\documentclass{article}
\usepackage{amsmath}
\usepackage{amssymb}
\usepackage{amscd}
\usepackage{cite}
\usepackage{threeparttable}
\usepackage{amsthm}

\theoremstyle{plain}

\newcommand{\vs}[1]{\vskip #1pt}
\usepackage{color}

 \newtheorem{The}{Theorem}[section]
 
 \newtheorem{Lem}[The]{Lemma}

 \newtheorem{defn}[The]{Definition}
 \newtheorem{Rem}[The]{Remark}
 
 \numberwithin{equation}{section}

\newcommand{\mc}{\mathcal}

\def\<{\langle}
\def\>{\rangle}

\def\be{\begin{equation}}
\def\ee{\end{equation}}
\def\p{\partial}

\catcode`\@=11

  \renewcommand{\section}%
  {\setcounter{equation}{0}\@startsection {section}{1}{\z@}{-3.5ex plus -1ex
   minus -.2ex}{2.3ex plus .2ex}{\Large\bf}}

\title{Long time stability for KAM tori of the derivative nonlinear Schr\"odinger equation\thanks{ Corresponding author: Xiaoping Yuan.}
\thanks{This work was partially supported by  National Natural Science Foundation of China (Grant No. 11871023 and No.12371189),  Guangdong Basic and Applied Basic Research Foundation (Grant No. 2021A1515111068) and  Shenzhen Science and Technology Innovation Program (Grant No. RCBS20210609103231040),.}}

\author{}

\author{ Shengqing Hu
\\
School of Science and Engineering,\\
The Chinese University of Hong Kong, Shenzhen\\
 Guangdong, 518172, China\\
{\tt e.mail: hushengqing@cuhk.edu.cn}
\\
\\
Huining Xue and  Xiaoping Yuan
\\
School of Mathematical Sciences, Fudan University\\
Shanghai, 200433, China\\
{\tt e.mail:21110180025@m.fudan.edu.cn;  xpyuan@fudan.edu.cn}
}

\date{}

\begin{document}
\maketitle
 \vs{12}

\tableofcontents

\centerline{\bf
Abstract}

\vs{12} {\small  This paper is concerned with the long time stability of KAM tori for a class of derivative nonlinear Schr\"odinger equations subjected to periodic boundary condition.
}

\vs{12} {\bf Keywords.} Stability;  Birkhoff normal form; KAM theory; Invariant tori; Hamiltonian system.
%\vs{12} \ni \vs{12} \ni
\section{Introduction and main results}
Let
 $H =H_{0}(I)+\varepsilon R(\theta,I)$ be a Hamiltonian function defined on some domain in phase space $\mathbb{T}^{n}\times \mathbb{R}^{n}$ endowed with the symplectic structure $dI\wedge d\theta$. Under suitable conditions, KAM theory \cite{kolmogorov1954conservation},\cite{arnold1963proof},\cite{moser1962} shows that there are many (in the sense of positive Lebesgue measure of the set of initial actions $I$'s) invariant tori which are today called KAM tori. On the other hand, there is such orbit $(I(t),\theta(t))$ satisfying that for any $A$ and $B$ with $0<A < B$, there are $t_{1}$ and $t_{2}$ such that  $|I(t_{1})|<A\ll B <I(t_2)|$. This kind of phenomenon is founded by Arnold\cite{Arn64} and called Arnold diffusion. The coexistence of both stable orbits (KAM tori) and unstable orbits (Arnold diffusion) leads to quite complicated dynamical behaviors. However, all orbits are long-time-stable in Nekhoroshev's sense \cite{Neh}:
\begin{equation*}
|I(t)-I(0)| \lesssim \varepsilon^{a}, \ for \ |t|\lesssim \exp{\varepsilon^{-b}},
\end{equation*}
where $I(0)$ is any initial action and $a, b>0$ are constants, provided that $H_0$ is steep. If $H_{0}(I)$ is quasi-convex (the steepest)\cite{Pos93},
\begin{equation}\label{yh}
a=b=\frac{1}{2n}.
\end{equation}
Consequently, the theories on KAM tori, Arnold diffusion and long time stability in Nekhoroshev's sense depict well the dynamical behaviours of nearly integrable Hamiltonian systems.

Since 1980's, many authors tried to generalize the theories mentioned above to  infinite dimensional Hamiltonian systems defined by some partial differential equations (PDEs).

\textbf{ KAM theory for PDEs.}

 Because the dimension of Hamiltonian system defined by PDEs is infinite, it is natural to seek the infinite dimensional
invariant tori (almost-periodic solutions) which are called full dimensional tori. In \cite{Bou96-GAFA} and \cite{P3} such kind of full dimensional  tori have been obtained for  nonlinear wave and Schr\"{o}dinger equations of spatial dimension $1$ subject to Dirichlet boundary conditions. However, the decay conditions  imposed on the initial sequence amplitudes are too strong to be natural. So far the problem on replacing the strong decay by a reasonable one is still open. See \cite{Bou96-GAFA} and \cite{cong}, for example. An option is to seek invariant tori (KAM tori) of lower (finite) dimension for PDEs.

In order to discuss KAM tori of lower dimension, let us consider the  Hamiltonian
\begin{equation}
 H =\langle w,y\rangle+ \<Az,\bar {z}\> +\varepsilon R(x,y,z,\bar{z}),
\end{equation}
with symplectic structure $dy\wedge dx +{\bf{i}}dz \wedge d\bar{z}$, $(x,y)\in {\mathbb T^{d}}\times \mathbb R^{d}$, $(z,\bar{z}) \in \mathcal{H}_{p} \times \mathcal{H}_{p}$, where $A$ is a linear operator with pure point spectra $\{\Omega_{j}, j \in \mathbb{Z}^{d}\}$ with $d \in \mathbb{N}$ being the spatial dimension of PDEs, and where $\mathcal{H}_{p}$ is the usual Sobolev space with weight $p>0$. Since the existence of KAM tori depends heavily on the structure of PDEs, we need to classify PDEs in terms of the resonances of linearized operator $A$ and unboundedness of nonlinear perturbation $R$. To the end, assume that the spectra of the linearized operator $A$ obey
\begin{equation*}
\Omega_{j}\approx|j|^{\kappa},\ \ \ \kappa \in \mathbb{R}.
\end{equation*}
Also assume the Hamiltonian vector field $X_{R}$ of the nonlinear perturbation $R$ maps $\mathcal{H}_{p}$ to $\mathcal{H}_{p'}$, that is,  $$X_{R}: \mathcal{H}_{p}\to \mathcal{H}_{p'}.$$ Let $\tilde\delta=p-{p'}$. We classify the vector field $X_{R}$ as follows (the classification is not complete):
\begin{enumerate}
  \item[Case 1.] bounded perturbation ($\kappa>0, \, \tilde\delta\le 0$).  Here are typical examples: Case(1.1).
  $\tilde\delta=0, \kappa=2,$
\be\label{NLS}{\bf i}
u_t+\Delta u+V(x) u+|u|^2u+\cdots=0,\ x\in \text{some compact subset of}\ \mathbb R^d,\text {(NLS)}.\ee

Case (1.2). $\tilde\delta=-1, \kappa=1,$
\be\label{NLW}u_{tt}-\Delta u+V(x) u+u^3+\cdots=0, \ x\in \text{some compact subset of}\,\, \mathbb R^d,\, \text {(NLW)}.\ee
  \item[Case 2.] unbounded perturbation ($\kappa>0,\, 0<\tilde{\delta}<\kappa-1$). Here two typical examples are perturbed KdV (pKdV) and Kadomtsev-Petviashvili (KP): $\tilde{\delta}=1,\kappa=3$,
  \be\label{kdv}  u_t+u_{xxx}+u_x u+\cdots=0, \; x\in\mathbb T, \, \text{ (pKdV)}.\ee
  \be\label{kp}  u_t+u_{xxx}+u_x u+\p_x^{-1}\p_{yy}u+\cdots=0, \, (x,y)\in\mathbb T^2, \; \text{(KP)}.\ee
  \item[Case 3.] unbounded critical perturbation ($\kappa>0\,\, \tilde\delta=\kappa-1$). Here a typical example is derivative nonlinear Schr\"{o}dinger equation: $\tilde\delta=1,\kappa=2$,
    \be\label{DNLS} {\bf i}
u_t+ u_{xx}+V(x) u+(|u|^2u)_x+\cdots=0,x\in \mathbb T,\; \text{(DNLS)}.\ee

  \item[Case 4.]  quasi-linear perturbation ($\kappa>0,\, \tilde\delta=\kappa$). A typical example is KdV with quasi-linear perturbation:
  \be\label{quasilinear-Kdv} u_t+u_{xxx}-6u u_x+F(x,u,u_x,u_{xx},u_{xxx})=0,\; x\in\mathbb T.\ee

  \item[Case 5.] fully nonlinear equation. Here is a typical example,the pure gravity water waves system (pgWW), with $\kappa=1/2$ and $\tilde\delta=1$,

  \begin{equation}\label{1.2xue}
\left\{ \begin{aligned}% requires amsmath; align* for no eq. number
 & \partial_{t}\Phi+\frac{1}{2}|\nabla \Phi|^{2}+g\eta=0  & at\ y=\eta(t,x),\\
 &\Delta\Phi=0   &in\ \mathcal{D}_{\eta},\\
 &\partial_{y}\Phi=0  &at\ y=-h,\\
 &\partial_{t}\eta=\partial_{y}\Phi-\partial_{x}\eta \cdot \partial_{x}\Phi  &at\ y=\eta(t,x),
\end{aligned} \right.
\end{equation}
where $g>0$ is the acceleration of gravity. The unknowns of the problem are the free surface $y=\eta(t,x)$ and the velocity potential $\Phi:\mathcal{D}\ \rightarrow \mathbb{R}$, i.e. the rotational velocity field $v=\nabla_{x,y}\Phi$ of the fluid.

  \item[Case 6.] $\kappa\le 0$. A typical example is the generalized Pochhammer-Chree equation (gPC): $\kappa=0$,
  \be\label{gPC} u_{tt}-\Delta u_{tt}-\Delta u+\Delta u^3+\cdots=0,\, x\in \mathbb T^d.  \ee

\end{enumerate}

 For the spatial dimension $d=1$ and bounded perturbations, there have been  too many results. Here are two original works \cite{K} and \cite{W} among them. For $d>1$ or unbounded perturbations, there are relatively less results. Here are a part of representative works. See \cite{K},\cite{Kuk00} \cite{KP},\cite{W}, \cite{P2},\cite{P1} for case (1); \cite{Kuk00} for case (2); \cite{LY11} for case (3); \cite{BBP13},\cite{baldi2016kam} for case (4);  \cite{baldi2018time} for case (5),\cite{yuan} for case (6).

\textbf{Arnold diffusion or instability for PDEs.}\  Consider the periodic defocusing cubic nonlinear Schr\"{o}dinger (NLS) equation
\begin{eqnarray*}
  &- {\bf i}{\partial _t}u + \Delta u ={\left| u \right|^2}u, \\
  & u\left( {0,x} \right)= {u_0}\left( x \right),\,  x \in \mathbb T^{2}=\mathbb{R}^{2}/(2\pi\mathbb{Z})^{2}.
\end{eqnarray*}
In \cite{CK} it is  proved that for given $K\gg 1$ and $0 < \delta\ll 1$, there exists global smooth solution $u(t,x)$ to the NLS and $T>0$ with
\begin{equation*}
\|u(0)\|_{\mathcal{H}_{p}} \leq \delta \ {\rm and} \ \|u(T)\|_{\mathcal{H}_{p}} \geq K.
\end{equation*}
See \cite{guardia2015growth} for the related works.

\textbf{Long time stability for PDEs.} \ Firstly, look at the long time stability for equilibrium points of PDEs. For PDEs, the freedom degree $n=+\infty$. From \eqref{yh}, we see that Nekhoroshev's theorem no longer valid for PDEs. Therefore, Bourgain \cite{Bou96-GAFA} proposes to restrict the initial action $I(0)$ to a sufficiently small neighborhood of some invariant set (equilibrium points and KAM tori, for example) and shorten the stable times $\mathrm{exp}(\epsilon^{-a})$ to $|t|\leq \varepsilon^{-M}$ for some large $M$. And in \cite{Bou96-GAFA} it is proved that  the solutions to
$${\bf i}u_{t}-u_{xx}+V(x)u+\varepsilon \frac{\partial H}{{\partial \bar{u}}}=0$$
and
$$y_{tt}-y_{xx}+\rho y +F'(y)=0$$ are $\varepsilon$-close to the unperturbed solution (linear flow) for times $|t|<T \thicksim \varepsilon^{-M},$
subject to Dirichlet boundary conditions, with smooth initial date of size $\varepsilon$, for ``typical" $V$ and $\rho$, respectively.  There are plenty references on this direction. Here we list part of them \cite{Bam03,BG,BFG88,imekraz2016long,baldi2019existence,ionescu2019long,biasco2020abstract,bernier2020long,berti2022birkhoff,yuan2014long,yuan2016averaging}. Especially, we mention Bambusi-Gr\'{e}bert's
work \cite{BG} where the tame property of the perturbation vector field $X_{R}$ is found and used to simplify the proof. Secondly, look at the long time stability of KAM tori for PDEs. It has been almost twenty  years since the study of long time stability of KAM tori was proposed by Eliasson \footnote{L. H. Eliasson, A talk in Fudan University, 2007} and there has been few results. In \cite{CLY}, the long time stability is investigated for NLS equation of space dimension $d=1$,
$${\bf i}u_{t}=u_{xx}+M_{\xi}u+\varepsilon |u^{2}|u+ \cdots =0,$$
subject to Dirichlet boundary conditions, where the potential $M_{\xi}$ is a Fourier multiplier. More exactly, for a given  KAM torus depending on a typical parameters $\xi$ and for any $0<\delta \ll 1$, there exist constants $p$ and $M$ with $M\gg 1$ and $p\gtrsim 8M^{4}$ such that any solution $u(t,x)$ with the initial data $u(0,x)$ $\delta$-close to the KAM torus $T_{\xi}$ in Sobolev norm $\|\cdot\|_{\mathcal{H}_{p}}$ obeys
$$\|u(t,x)-T_{\xi}\|_{\mathcal{H}_{p}} \leq 2\delta , \ for\  |t|\leq \delta^{-M}.$$
 See \cite{cong2015long} for a similar result for NLW with $d=1$.  Note that $\delta$ is independent of the perturbation size $\varepsilon$. For the nonlinear Schr\"odinger equations of space dimension $d>1$, the long time stability of KAM tori is investigated in \cite{he2021long}, where the perturbation does not involve the space variable $x$ explicitly. See \cite{maspero2018long} and \cite{kappeler2021stability} for some related results.

As for the nonlinear wave equation with $d>1$ and those PDEs in case (2)-(6), the long time stability of KAM tori is still open. See {\it Tables 1 and 2.}

\begin{table}\caption{$\kappa>0,d=1$ (Existence and stability of KAM Tori)}\label{table}
  \centering
  \begin{threeparttable}
 \begin{tabular}{|c||c|c|c|c|c||}
  \hline
 & NLS & NLW & pKdV
& DNLS & WW\\
 \hline
  \textbf{KAM} & $\surd$ & $\surd$ & $\surd$ & $\surd$& $\surd$\\
  %\hline
  %\textbf{Instability\tnote{$\ast$}} & ? $\ddagger$& ?& ? &? &?\\
  \hline
  \textbf{Stability\tnote{$\dagger$}} & $\surd$&$\surd$
&? & ?\tnote{\S} & ? \\
  \hline
  \end{tabular}
\begin{tablenotes}
        \footnotesize
        \item[$\ast$]  Arnold diffusion.
        \item[$\dagger$]  Long time stability of KAM tori.
        \item [$\ddagger$]? means unsolved.
        \item[$\S$] Solved in the present paper.
      \end{tablenotes}
\end{threeparttable}
\end{table}

\begin{table}\caption{$d>1$(Existence and stability of KAM Tori)}\label{table}
  \centering
  \begin{threeparttable}
 \begin{tabular}{|c||c|c|c|c||}
  \hline
  & NLS & NLW & KP
& gPC \\
 \hline
  \textbf{KAM} & $\surd$ & $\surd$ & ? & $\surd$\\
  %\hline
  %\textbf{Instability} & $\surd$\tnote{$\ast$} & ?& ? &? \\
  \hline
  \textbf{Stability} & ?&?
&? & ?  \\
  \hline
  \end{tabular}
%\begin{tablenotes}
    %    \footnotesize
    %    \item[$\ast$]   $d=2$.
      %        \end{tablenotes}
\end{threeparttable}

\end{table}

All the long time stability results on KAM tori mentioned above are for those PDEs of bounded nonlinear perturbation. In the present paper, we focus on the long time stability of KAM tori for nonlinear Schr\"odinger equations of unbounded nonlinear perturbation. More exactly, for
 a class of derivative nonlinear Schr\"odinger equations subject to periodic boundary condition
\begin{align}\label{main equation}
&{\bf i}u_t+u_{xx}-\Xi(x)*u+{\bf i}\epsilon (|u|^2u)_x=0, \quad(t,x)\in \mathbb{R}\times \mathbb{T},
%\\& u(t,0)=u(t,\pi)
\end{align}
where $\Xi(x)*u$ denotes the convolution function between the potential $\Xi: \mathbb{T}\rightarrow \mathbb{C}$ and the function $u$. Moreover, the function $\Xi$ is analytic and the Fourier coefficients are real, when expanding $\Xi$ into Fourier series $\Xi(x)=\sum_{j\in\mathbb{Z}\setminus\{0\}}\xi_je^{{\bf i} jx}$.
%where $M_\xi$ is a real Fourier multiplier,
%$$M_\xi e^{2\pi {\bf i} jx}=\xi_je^{2\pi {\bf i} jx},\quad \xi_j=\xi_{-j}\in\mathbb{R},\quad j\in \mathbb{Z},$$
Let the parameter set be
$$\Pi=\left\{\xi=(\xi_j)_{j\in \mathbb{Z}\setminus \{0\}}: \xi_j \in [1,2]/|j|\right\}.$$
%It has positive Lebesgue measure, i.e., ${\rm Meas}\ \Pi:=\prod_{j\in\mathbb{Z}\setminus\{0\}}|[1,2]|=1>0$.
%%%
Then, we have the following main result.
\begin{The}\label{theorem1.1}
Consider equation \eqref{main equation}. Given an integer $n\geq 1$ and a real number $p\geq 1$. If $\epsilon$ is sufficiently small, then there exists a positive measure set ${\Pi}_\eta\in\Pi$ with
$$\lim_{\eta\to 0}\frac{{\rm Meas}\  \Pi_\eta}{{\rm Meas}\ \Pi}=1,$$
%%%
where $\eta$ is some constant in $(0,1)$. For any $\xi\in {\Pi}_\eta$, the nonlinear Schr\"odinger equation \eqref{main equation} possesses a linearly stable $n$-dimensional KAM torus $\mathcal{T}_\xi$ in Sobolev space $\mathcal{H}^p(\mathbb{T})$.

Moreover,  for $0\le \mathcal{M}\le C(\epsilon)$ (where $C(\epsilon)$ is a constant depending on $\epsilon$ and $C(\epsilon)\rightarrow \infty$ as $\epsilon \rightarrow 0$) and $p\geq 8C_0(\mathcal{M}+7)^4+1$ (where $C_0$ is a constant depending on $n$), there exists small $\delta_0>0$ depending on $p,\mathcal{M}$ and $n$, such that  for  any $0<\delta<\delta_0$, there exists a subset $\Pi_{\acute \eta}\subset \Pi_\eta$ with
 $$\lim_{\acute{\eta}\to 0}\frac{{\rm Meas}\ \Pi_{\acute\eta}}{{\rm Meas}\ \Pi_\eta}=1,$$
where  $\acute\eta$ is some constant in $(0,1)$. For any $\xi\in\Pi_{\acute \eta}$, any solution $u(t,x)$ to equation \eqref{main equation} with the initial datum satisfying
  $$d_{\mathcal{H}^p(\mathbb{T})}(u(0,x),\mathcal{T}_\xi):=\inf_{w\in \mathcal{T}_\xi}\Vert u(0,x)-w\Vert_{\mathcal{H}^p(\mathbb{T})}\le \delta,$$
obeys
 $$d_{\mathcal{H}^{p/2}(\mathbb{T})}(u(t,x),\mathcal{T}_{\xi}):=\inf_{w\in \mathcal{T}_\xi}\Vert u(t,x)-w\Vert_{\mathcal{H}^{p/2}(\mathbb{T})}\le 2\delta,\quad \forall \ \ |t|\le \delta^{-\mathcal{M}/4}.$$
\end{The}

%%%%%%%%%%%????????????

\iffalse%%%%%%%%%%%%%%%%%
In the present paper, we are addressed to the long time stability of the KAM tori for the derivative nonlinear Schr\"{o}dinger equation
\begin{equation}
{\bf i}u_{t}+u_{xx}-v(x)\ast u+ {\bf i}\varepsilon(|u|^{2}u)_{x}=0,\ x\in \mathbb{T}.
\label{1.5+}
\end{equation}
%%%
\fi
%%

\begin{Rem}
The measure of the parameter set $\Pi\subset\mathbb{R}^{\mathbb{Z}\setminus\{0\}}$ in the above theorem is in the sense of Kolmogorov.
\end{Rem}

\begin{Rem}
Equation \eqref{main equation} is in case (3). Thus, the second question mark `` $?$" could be replaced by  check mark ``$\checkmark$" in Table 1.
\end{Rem}

%An abstract KAM Theorem is constructed in \cite{LY11} for case (3) and is used to prove the %existence of KAM tori of \eqref{main equation} for typical $\Xi(x)$. The long time stability of %the KAM tori has been open for \eqref{main equation} .

For the readers' convenience, here are the outline of the proof of Theorem \ref{theorem1.1}, as well as some remarks:

%%%%%%%%

{\bf Step 1}: Construct a KAM theorem in terms of Kolmogorov's original idea \cite{kolmogorov1954conservation}, that is, by a series of symplectic coordinate changes, to reduce Hamiltonian
$$H=\<w,y\>+\<Az,\bar{z}\>+\varepsilon R(x,y,z,\bar{z}),(x,y,z,\bar{z}) \in D(s_{0},r_{0})$$
to
\begin{equation}
H_{\infty}=\<\tilde{w},y\>+\<\tilde{A}(x)z,\bar{z}\>+\varepsilon R_{\infty}(x,y,z,\bar{z}),\ (x,y,z,\bar{z}) \in D\left(\frac{s_{0}}{{2}} ,\frac{r_{0}}{{2}}\right),
\label{1.6+}
\end{equation}
where $A=(\Omega_j)_{j\in\mathbb{Z}}$ is a constant-valued diagonal operator and $\tilde{A}(x)=(\tilde{\Omega}_j)_{j\in\mathbb{Z}}$ is $x$-dependent diagonal operator, and
\begin{equation}
R_{\infty}=O(|y|^{2}+|y|\|z\|_{p} +\|z\|_{p}^{3}),
\label{1.7+}
\end{equation}
and
$$D(s,r)=\{(x,y,z,\bar{z}) \in (\mathbb{C}/ 2\pi \mathbb{Z})^{n} \times \mathbb{C}^{n} \times \mathcal{H}_{p} \times \mathcal{H}_{p} \mid |\text{Im} x|<s, |y|<r^{2},|z|<r,|\bar{z}|<r\}.$$
This task is fulfilled in Sections 2 and 3.

Here are some remarks. It follows from \eqref{1.7+} that the set $\mathbb{T}^n\times\{y=0\} \times \{z=\bar{z}=0\}$ is an invariant torus (KAM torus) of the vector field defined by $H_{\infty}$. The key point is that the definition domain $D(\frac{s_{0}}{{2}} ,\frac{r_{0}}{{2}})$ of $R_{\infty}$ is a neighborhood of KAM torus $\mathbb{T}^n\times\{y=0\} \times \{z=\bar{z}=0\}$, which provides a working room for the long time stability of the KAM torus.  While the definition domain of $R_{\infty}$ is $D(\frac{s_{0}}{{2}} ,0)$ in the usual KAM theorem by Arnold's
idea \cite{arnold1963proof}. In \cite{LY11}, a KAM theorem has been established for Schr\"odinger equation \eqref{main equation} with the definition domain $D(0,0)$. Actually, we can expand $D(0,0)$ to $D(0,\frac{r_{0}}{{2}})$ in terms of Kolmogorov's original idea \cite{kolmogorov1954conservation} and that in \cite{LY11}. The  domain $D(0,\frac{r_{0}}{2})$ can provide the linear stability of the KAM torus. See \cite{yuan2013reduction}.  However, the expanded domain $D(0,\frac{r_{0}}{{2}})$ is not yet enough to establish the long time stability of the KAM torus. While expanding $D(0,\frac{r_{0}}{{2}})$ to $D(\frac{s_0}{2},\frac{r_{0}}{{2}})$, a basic ingredient is a suitable estimate for the so-called ``small denominator equation with large variable coefficient" by \cite{kuksin1997small}
\begin{equation}
-{\bf i} w\cdot \partial_ x u+(\lambda +\mu(x))u=r(x),\ |\text{Im} \, x|\leq s,
\label{1.8+}
\end{equation}
where ${\bf i}^2=-1$ and $\lambda \gg 1$ is constant, the coefficient $\mu(x)$ depends on the angle variable $x \in \mathbb{T}$, and it is large: $\mu(x) \thicksim \lambda ^{\theta}$.
When $0< \theta <1$, in \cite{kuksin1997small} Kuksin gives an estimate
\begin{equation}\label{1.9++}
\mathop {\sup }\limits_{\left| {{\mathop{\rm Im}\nolimits} x} \right| \le s - \sigma }|u(x) |\leq C_{1}\exp(C_{2}C_{3}^{\frac{1}{1-\theta}})\mathop {\sup }\limits_{\left| {{\mathop{\rm Im}\nolimits} x} \right| \le s} |r(x)|,\ \forall\ 0< s \ll 1.
\end{equation}
This estimate is suitable for case (2) which can be applied to the perturbed KdV equation among others.

When considering DNLS equation, we have that the vector field $X_{R}:\mathcal{H}_{p}\to \mathcal{H}_{p-1}$ is in case (3): critical unbounded perturbation.
At this time, $\theta=1$. One sees that the estimate \eqref{1.9++} is invalid for $\theta=1$.
 In \cite{liu2010spectrum}, for $\theta=1$, Liu and Yuan gives out an estimate
\begin{equation}\label{1.10+yuan}
\mathop {\sup }\limits_{\left| {{\mathop{\rm Im}\nolimits} x} \right| \le s - \sigma }| u(x) |\leq \sigma^{C_{4}}\exp(c\gamma s )\mathop {\sup }\limits_{\left| {{\mathop{\rm Im}\nolimits} x} \right| \le s} |r(x)|,\ \forall\ 0< s \ll 1.
\end{equation}
%%%%
Although Liu-Yuan's estimate is worse than Kukuk's one, it is enough to get $R_{\infty} =O(|y|^{2} +|y|\|z\|_{p}+ \|z\|_{p}^{3}), \ (x,y,z,\bar{z}) \in D(0,0))$,  which guarantees the existence of the KAM tori. As mentioned above, $D(0,0)=\mathbb{T}^n\times \{y=0\} \times \{z=\bar{z}=0\}$ is not a neighborhood of the torus $\mathbb{T}^n\times\{y=0\} \times \{z=\bar{z}=0\}$. This excludes the possibility of proving the long time stability of the torus. In the present paper, we improve \eqref{1.10+yuan} by absorbing the method by the Italian School's works \cite{baldi2018time},\cite{baldi2016kam},\cite{BBP13}, etc., which deal with   cases (4) and (5). More precisely, in Section 6, we show that for $0<\theta\le 1$
\begin{equation}\label{1.11+}
\mathop {\sup }\limits_{\left| {{\mathop{\rm Im}\nolimits} x} \right| \le s-\sigma}|u(x)| \le \sigma^{-C_5} \mathop {\sup }\limits_{\left| {{\mathop{\rm Im}\nolimits} x} \right| \le s }|r(x)|,\ \forall \ 0<\sigma<s.
\end{equation}
Once getting this estimate, we can extend the definition domain of $R_\infty$ to $D(s_0/2,r_0/2)$ by following Kolmogorov's iteration procedure\cite{kolmogorov1954conservation}. We believe that the improved estimate \eqref{1.11+} is of more use in the related theory on small divisors.

\medskip

{\bf Step 2}: Write $z=(\acute{z},\hat{z})$ , $\hat{z} = (z_{j}:|j|\leq \mathcal{N})$, $\acute{z} = z \ominus \hat{z}$, for some large $\mathcal{N}$. In Section 4,
by means of Birkhoff normal form, we reduce $H_{\infty}$ to
$$H_{\infty}^{\ast}=\<\tilde{w},y\>+\<\tilde{A}(x)z,\bar{z}\>+Z(x,y,z,\bar{z})+R_{\infty}^{\ast},$$
where
$$Z(x,y,z,\bar{z})=\sum_{4\le 2|\alpha|+|\beta|+|\mu|\le \mathcal{M}+2,|\mu|\le 1}Z^{\alpha\beta\mu}(x)y^{\alpha} |\acute{z}|^{2\beta} |\hat{z}|^{2\mu},$$
$$R_{\infty}^{\ast} = O(|\tilde{y}|^{\mathcal{M}+1}+\|\acute{z}\|_{p}^{3}),$$
 and $\tilde{y}=(y^{\frac{1}{2}},z_{j},\bar{z}_j, |j|\leq \mathcal{N})$.  This Birkhoff normal form is somehow different from that of bounded Hamiltonian. See \cite{CLY} for the bounded case, where the coefficients $\tilde{A}(x)$ and $Z^{\alpha\beta\mu}(x)$ are independent of the angle variable $x$. Since the Hamiltonian vector field is unbounded in our case, we could not find a bounded transformation to make those  terms independent of $x$. This leads to the Birkhoff normal form involves angle variable $x$, which is not integrable any more. See Step 4 below.

In addition, the Liu-Yuan's estimate \eqref{1.10+yuan} is still of its own use, although it is improved by \eqref{1.11+}. As the Kuksin's estimate, \eqref{1.11+} holds under condition $|u(x)| \leq C\lambda.$ While in constructing the Birkhoff normal $H_{\infty}^{*}$, the case $|u(x)| \geq C\lambda$ might appear. At this time, we have to appeal to \eqref{1.10+yuan} rather than \eqref{1.11+}. Recall that the definition domain of $H_{\infty}$ is $D(\frac{s_{0}}{2}, \frac{r_{0}}{2})$ for $r_{0},s_{0} >0.$  Once \eqref{1.10+yuan} is used, the domain $D(\frac{s_{0}}{2}, \frac{r_{0}}{2})$ will shrink rapidly to $D(0,0)$. Fortunately, we use \eqref{1.10+yuan}, at most, $\mathcal{M} < \infty$ times. Thus, the  choice of $\mathcal{M},\mathcal{N}$ and $\delta$ in Theorem \ref{theorem1.1} is very delicate. See Section 4.

\medskip

{\bf Step 3}: By absorbing the idea in  \cite{BG,CLY,yuan2014long}, we verify that the vector field $X_{R_{\infty}^{\ast}}$ has tame property, thus,
$$\|X_{R_{\infty}^{\ast}}\|_{p/2} \leq \delta ^{\mathcal{M}+1} +\| \acute{z}\| \lesssim \delta^{\mathcal{M}+1} + \frac{1}{\mathcal{N}} \lesssim \delta ^{\mathcal{M}+1}$$
 by taking $\mathcal{N} \thicksim \delta ^{-(\mathcal{M}+1)}$. Remark that this step is not trivial, since $R_{\infty}^{\ast}$ depends on the angle variable $x$ and the vector field $ X_{R^*_\infty}$ is unbounded. Instead of using the tame norm in \cite{CLY}, we use the sup norm of the coefficients of Hamiltonian. (See Definition \ref{defnd} for the detail). Under this norm, by a method similar to that in \cite{yuan2014long}, we can prove the vector field of the Hamilton functions with momentum conservation has the tame property. See Lemma \ref{guji}.  Finally, we estimate the norm of the Possion bracket in Lemmas \ref{brak} and \ref{brak111}.

 %In the present paper, we obtain the long time stability result with initial datum in $\mathcal{H}_p$ space and the solution in $\mathcal{H}_{p-1}$ space, where $p$ is only depend on $\epsilon$. Actually,  if one considers the term $\{ \Vert z\Vert_p^2, R_\infty^*\}$ as in \cite{yuan2014long}, one can also obtain the long time stability result with initial datum in $\mathcal{H}_p$ space and the solution in $\mathcal{H}_{p}$ space. However, the Sobolev index $p$ will depend on $\delta$ by using the method in this paper and the proof process will be more complex.

 \medskip

{\bf Step 4}: In step 2, we obtained a Birkhoff normal form $H_{\infty}^{\ast}$. Which can afford the long time stability estimate in $z$ and $\bar{z}$ direction. However, this is not enough for the long time stability estimate in $y$ direction, since the coefficients $\tilde{A}(x)$ and $Z^{\alpha\beta\mu}(x)$ depend on the angle variable $x$. In this step, we will further do some transformations to make these terms independent of $x$ or being small enough ($\le \delta^{\mathcal{M}+1}$) such that the long time stability can be derived. More precisely, we reduce  $H_{\infty}^{\ast}$ into
$$H_{\infty}^{**}=N+Z+B+R_{\infty}^{**},$$
with
$$N=\sum_{j=1}^n\omega^*_j y_j+\sum_{j\geq 1}[\breve{\Omega}_j]|z_j|^2+\sum_{j\geq 1}(\breve{\Omega}_j-[\breve{\Omega}_j])(x)(|z_j|^2-|z_j(0)|^2),$$
$$%Z=\sum_{4\le 2|\alpha|+|\beta|+|\mu|\le \mathcal{M}+2,|\mu|\le 1}\tilde{Z}^{\alpha\beta\mu}y^{\alpha} |\acute{z}|^{2\beta} |\hat{z}|^{2\mu}+
Z=\sum_{|\alpha|\le \mathcal{M}+2}\tilde Z^\alpha(\xi) y^\alpha+
\sum_{|\alpha|\geq \mathcal{M}+2}\tilde{Z}^\alpha(x;\xi) y^\alpha,$$
$$
B=\sum_{ |\beta|\geq 1}B^{\alpha\beta}(x,y;\xi)(|z|^{2}-|z(0)|^2)^{\beta},$$
$$R_{\infty}^{**}=O(|\tilde{y}|^{\mathcal{M}+1}+\|\acute{z}\|_{p}^{3})+O(\delta^{\mathcal{M}+1}),$$
where the coefficients $\tilde{Z}^{\alpha\beta\mu}$ is independent of $x$ and  $\Vert X_{R_{\infty}^{**}}\Vert_{p/2}\le \delta^{\mathcal{M}+1}$.
In this process, we have used KAM iterative method to eliminate the function of the form
$$\sum_{0\le |\alpha|\le \mathcal{M}+1,|\beta|\geq 1}g^{\alpha\beta}(x)y^{\alpha}|z(0)|^{2\beta}.$$
Since $z(0)$ is the initial datum, so the norm of this function  will be less than $\delta^2$. Thus, KAM iteration works.

Finally one can easily verify that the KAM torus $\mathbb{T}^n\times\{y=0\} \times \{z=\bar{z}=0\}$ of $H_{\infty}$ is stable for $|t|\lesssim \delta^{-\mathcal{M}/2}$. Consider the Hamiltonian $H_{\infty}^{**}$. The terms $\sum_{j\geq 1}[\breve{\Omega}_j]|z_j|^2$ and
$$%\sum_{4\le 2|\alpha|+|\beta|+|\mu|\le \mathcal{M}+2,|\mu|\le 1}\tilde{Z}^{\alpha\beta\mu}y^{\alpha} |\acute{z}|^{2\beta} |\hat{z}|^{2\mu}
\sum_{|\alpha|\le \mathcal{M}+2}\tilde Z^\alpha(\xi) y^\alpha,$$
are integrable terms which do not effect the stability. The terms
$$\copyright:=\sum_{j\geq 1}(\breve{\Omega}_j-[\breve{\Omega}_j])(x)(|z_j|^2-|z_j(0)|^2)$$ and $B$ depend on $|z|^{2}-|z(0)|^2$.  Since we have stability estimate in $z$ and $\bar{z}$ direction, it means that $\sum_{j}|j|^{2p/2}||z_j|^2-|z_j(0)|^2|\le \delta^{\mathcal{M}+1}|t|$. Substituting it into $\copyright$ and $B$, we immediately obtain the long time stability estimate with time length $|t|\le \delta^{-\mathcal{M}/2}$. See Section \ref{lthea}.

\begin{Rem}
Another difficulty we encountered is the measure estimate. It is under some proper non-resonant conditions that estimates \eqref{1.10+yuan} and \eqref{1.11+} hold true. In other words, one needs to remove some resonant sets consisting of ``bad" parameters. For example, in the process of obtaining the Birkhoff normal form, we need to eliminate the resonant set
$$\{\xi\in \Pi: \langle k, \omega(\xi)\rangle+\langle \acute{l}, \acute{\Omega}(\xi)\rangle+ \langle \hat{l},\hat{\Omega}(\xi)\rangle\, {\rm is}\,{\rm small}\},$$
where $k\in\mathbb{Z}^n, \acute{l}\in\mathbb{Z}^{2\mathbb{\mathcal{N}}}, \hat{l}\in\mathbb{Z}^{\mathbb{Z}_*}$ and $|k|+|\acute{l}|+|\hat{l}|\neq 0, |\acute{l}|+|\hat{l}|\le\mathcal{M}+2, |\hat{l}|\le 2$. Note that,  the normal frequencies $\acute{\Omega}(\xi)=(\Omega_{j}(\xi))_{|j|\le \mathcal{N}}$ and $\hat{\Omega}(\xi)=(\Omega_j(\xi))_{|j|>\mathcal{N}}$ satisfy $\Omega_j=j^2+\xi_j+O(|j|\epsilon)$.
We hope that the Lebesgue measure of the resonant set is small. For instance, we consider the case $k=0$ and
$$ \langle \hat{l},\hat{\Omega}(\xi)\rangle:=\Omega_i-\Omega_j\approx i^2-j^2> C\max\{|i|,|j|\}.$$
Thus, for $|i|,|j|>|\langle \acute{l}, \acute{\Omega}(\xi)\rangle|:=C_0(\mathcal{M},\mathcal{N})$, the divisor $|\langle \acute{l},  \acute{\Omega}(\xi)\rangle+ \langle \hat{l},\hat{\Omega}(\xi)\rangle|$ is not small. So we only need to discuss the case $|i|,|j|<C_0(\mathcal{M}, \mathcal{N})$. To obtain the estimate of Lebesgue measure, we need to verify that $ \langle k, \omega(\xi)\rangle+\langle \acute{l}, \acute{\Omega}(\xi)\rangle+ \langle \hat{l},\hat{\Omega}(\xi)\rangle$ is twisted with respect to $\xi\in\Pi$. In a standard way, we letting $|\acute{l}_b|=\max_{|j|\le \mathcal{N}}\{|\acute{l}_j|\}$, we have
\begin{align*}
&|\partial_{\xi_{b}}(\langle k,\breve{\omega}\rangle+\langle \acute{l},\acute{\Omega}\rangle+\langle \hat{l},\hat{\Omega}\rangle)|\\
\geq &|\acute {l}_b||\partial_{\xi_{b}}\acute{\Omega}_b|-|\partial_{\xi_{b}}(\langle \acute{l},\acute{\Omega}\rangle+\langle \hat{l},\hat{\Omega}\rangle-\acute{l}_p\acute{\Omega}_b)|\\
\geq& |\acute{l}_b|(1-O(|b|\epsilon))-\left(\sum_{j\neq b}|\acute{l}_j|O(|j|\epsilon)+\sum_{j\geq \mathcal{N}, |\hat{l}|\le 2}|\hat{l}_j|O(|j|\epsilon)\right) \\
\geq& |\acute{l}_b|-(\mathcal{M}+2)C_0(\mathcal{M},\mathcal{N})O(\epsilon)\\
> &0,
\end{align*}
provided that $(\mathcal{M}+2)C_0(\mathcal{M},\mathcal{N})O(\epsilon)<\frac{1}{2}$. However, as $\mathcal{N}$ is chosen to be a very large constant (It usually depends on the constant $\delta$ which is independent of $\epsilon$),  this condition can not hold. Thus, in this way, we fail to verify the twist property and to obtain the measure estimate.
Fortunately, we observe that the tangent and normal frequencies has a good property:
\begin{equation}\label{huru}
\partial_{\xi_a}\omega_i\approx \frac{1}{|a|}O(\epsilon),\quad \partial_{\xi_a}\Omega_i\approx \delta_{ja}+\frac{|j|}{|a|}O(\epsilon), \quad 1\le i\le n,a,j\in\mathbb{Z}^*,
\end{equation}
where $\delta_{ja}$ is the Dirac function. By this property, we can overcome the difficulty and obtain the measure estimate.
More exactly,  define
$$j_{max}=\max\{j: \acute{l}_j\, {\rm or}\, \hat{l}_j \neq 0\},\quad j_{min}=\min\{j: \acute{l}_j\, {\rm or}\, \hat{l}_j \neq 0\}
$$
$$J:=\left\{
\begin{array}{l}
j_{max},\quad {\rm if}\ j_{max}>-j_{min},\\
j_{min},\quad {\rm if}\ j_{max}\le-j_{min}.\\
\end{array}
\right.$$

Without loss of generality, we assume $\acute{l}_J\neq 0$. Then, from \eqref{huru}, we have
\begin{align*}
&|\partial_{\xi_{J}}(\langle k,\breve{\omega}\rangle+\langle \acute{l},[\acute{\Omega}]\rangle+\langle \hat{l},[\hat{\Omega}]\rangle)|\\
\geq &|\acute {l}_J||\partial_{\xi_{J}}[\acute{\Omega}]_J|-|\partial_{\xi_{J}}(\langle \acute{l},[\acute{\Omega}]\rangle+\langle \hat{l},[\hat{\Omega}]\rangle-\acute{l}_J[\acute{\Omega}]_J)|\\
\geq& |\breve{l}_J|(1-O(\epsilon))-\left(\sum_{j\neq J}|\acute{l}_j|\frac{|j|}{|J|}O(\epsilon)+\sum_{j\geq \mathcal{N}, |\hat{l}|\le 2}|\hat{l}_j|\frac{|j|}{|J|}O(\epsilon)\right) \\
\geq& |\breve{l}_J|-(\mathcal{M}+2)O(\epsilon)\\
> &0,
\end{align*}
if we let $(\mathcal{M}+2)O(\epsilon)<1$. Thus, property \eqref{huru} plays an important role in the measure estimate. See Section \ref{Mease}.
\end{Rem}

\begin{Rem}
In the present paper, the moment conservation plays an important role, as in \cite{LiuYuan2014}. When nonlinearity in \eqref{main equation} depends explicitly on $x$, the moment conservation is not valid. So the long time stability is still open for this case.
\end{Rem}

%\begin{Rem}
%Another critical point is to prove that the tame property of $X_{R}$ can be preserved through the symplectic coordinate change defined by the small denominator equation \eqref{1.8+} with large variable coefficients . The variable coefficients and unboundedness of $X_{R}$ lead to very complicated calculation in proving the tame property of $X_{R_{\infty}}$ and $X_{R_{\infty}^{\ast}}$. Finally, in the line of the procedure in \cite{CLY} with the modification as above, we have the following theorem:
%\end{Rem}

\section{KAM theorem}
To finish the proof of Theorem \ref{theorem1.1}, several abstract theorems about infinitely dimensional Hamiltonian system are given. We firstly give some notations as a preliminary.

\subsection{Preliminary}

\subsubsection{Hamiltonian formulation of DNLS equation}\label{2.1.1}
Set
$$e_j(x)=\frac{1}{\sqrt{2\pi}{\bf i}}e^{{\bf i}jx},\quad j\in\mathbb{Z}\setminus\{0\},$$
 and
$$u=\sum_{j\in \mathbb{Z}\setminus\{0\}}q_j e_j.$$

Then equation \eqref{main equation} can be written as
$$\dot{q}_j={\bf i}j \frac{\partial H}{\partial \bar{q}_j},\quad j\in  \mathbb{Z}\setminus\{0\},$$
with Hamiltonian
$$H(q,\bar{q})=\sum_{j\in  \mathbb{Z}\setminus\{0\}}\frac{\lambda_j}{j}|q_j|^2+ P(q,\bar{q}),$$
where
$$\lambda_j=j^2+\xi_j,\quad j\in \mathbb{Z}\setminus\{0\},$$
$$P(q,\bar{q})=\frac{\epsilon}{4}\sum_{i-j+k-l=0}P_{ijkl}q_i\bar{q}_jq_k\bar{q}_l,\quad P_{ijkl}=\int_0^{2\pi}e_ie_{-j}e_ke_{-l}dx=\frac{1}{2\pi}.$$
The symplectic form is $-{\bf i}\sum_{j\in\mathbb{Z}\setminus \{0\}}j^{-1}dq_j\wedge d\bar{q}_j$.

Fix $ j_1,\cdots,j_n \in \mathbb{Z}$ and let $J_n=\{j_1,\cdots,j_n\}$. Now we choose $n$ modes $e_{j_1},\cdots,e_{j_n}$ as tangent direction and the remained as normal direction. Let
$$\tilde{q}=(q_{j_1},\cdots,q_{j_n}),\quad z=(q_j)_{j\in\mathbb{Z}_*=\mathbb{Z}\setminus (\{0\}\cup J_n)},$$
be the tangent and normal variables, respectively. Then rewrite $P(q,\bar{q})$ in the multiple-index as
\begin{equation}\label{pqq}
P(q,\bar{q})=\epsilon\sum_{|\mu|+|\nu|+|\alpha|+|\beta|=4}P^{\mu\nu\alpha\beta}\tilde{q}^{\mu}\bar{\tilde{q}}^{\nu}z^{\alpha}\bar{z}^{\beta}.
\end{equation}

Introduce action-angle variables
$$q_{j_i}=\sqrt{j_i(\zeta_i+y_i)}e^{{\bf i}x_i},\quad \bar{q}_{j_i}=\sqrt{j_i(\zeta_i+y_i)}e^{-{\bf i}x_i}, \quad i\in\{1,\cdots,n\},$$
where each $\zeta_i\in [1,2]$ is the initial datum and will be considered as a constant. Then symplectic structure is
\begin{equation}\label{symplectic structure}
\sum_{1\le i\le n}dy_i\wedge dx_i-\sum_{j\in\mathbb{Z}_*}{\bf i}j^{-1}dz_j\wedge d\bar{z}_j.
\end{equation}

 Hence, Hamiltonian \eqref{pqq} becomes
\begin{align*}
P(q,\bar{q})=P(x,y,z,\bar{z})=\epsilon\sum_{|\alpha|+|\beta|\le 4}P^{\alpha\beta}(x,y)z^{\alpha}\bar{z}^{\beta},
\end{align*}
where
$$P^{\alpha\beta}(x,y)=\sum_{|\mu|+|\nu|=4-|\alpha|-|\beta|}P^{\mu\nu\alpha\beta} e^{{\bf i}\langle \mu-\nu,x\rangle}\sqrt{(\zeta+y)^{\mu+\nu}}\sqrt{\prod_{i=1}^nj_i^{\mu_i+\nu_i}}.$$

%Then, since $P^{\mu\nu\alpha\beta}=\frac{1}{8\pi}$ for all $|\mu|+|\nu|+|\alpha|+|\beta|=4$, we have the symmetric property:
%$$P^{\alpha\beta}(x,y)=P^{\beta\alpha}(x,y)\in\mathbb{R}.$$
%If we write $P^{\alpha\beta}(x,y)=\sum_{\gamma\in\mathbb{N}^n}P^{\gamma\beta\alpha}(x)y^\mu$, then it is also equivalent to
%$$P^{\gamma\alpha\beta}(x)=P^{\gamma\beta\alpha}(x)\in\mathbb{R}.$$
In addition, we have the momentum conservation: for the monomials in the Hamiltonian
$$e^{{\bf i}(k_1x_1+\cdots+k_nx_n)}y_1^{m_1}\cdots y_n^{m_n}\prod_{j\in \mathbb{Z}_*}z_j^{\alpha_j}\bar{z}_j^{\beta_j},$$
the following equality holds:
$$M(k,\alpha,\beta):=-\sum_{i=1}^nk_ij_i+\sum_{j\in \mathbb{Z}_*}j(\alpha_j-\beta_j)=0.$$

As consequence, we obtain a Hamiltonian $H(x,y, z,\bar{z};\xi)$ of the form
$$H(x,y, z,\bar{z};\xi)=N(y, z,\bar{z};\xi)+ P(x,y, z,\bar{z};\xi),$$
where
$$
N(y, z,\bar{z};\xi)=\sum_{i=1}^n\omega_i(\xi)y_i+\sum_{j\in\mathbb{Z}_*}\frac{\Omega_j(\xi)}{j}z_j\bar{z}_j,
$$
with the tangent frequency
\begin{equation}\label{oi}
\omega(\xi)=(\omega_{1}(\xi),\cdots,\omega_{n}(\xi))=\left(j_1^2+\xi_{j_1},\cdots,j_n^2+\xi_{j_n}\right),
\end{equation}
and the normal frequency $\left(\Omega_{j}(\xi)\right)_{j\in\mathbb{Z}_*}$ where
\begin{equation}\label{Oi}
(\Omega_{j}(\xi))_{j\in\mathbb{Z}_*}=\left(j^2+\xi_{j}\right)_{j\in\mathbb{Z}_*},
\end{equation}
and where the perturbation $P(x,y,z,\bar{z})$ is independent of $\xi$.

Introduce an infinite dimensional symplectic phase space
$$(x,y,z,\bar{z})\in\mathcal{P}^p:=\mathbb{T}^n\times \mathbb{C}^n\times l_p^2\times l_p^2.$$
where the space $l_p^2$ is defined in \eqref{lpp2}.Then it is easy to verify that the vector field $X_{P}$ of $P$ obeys
$$X_P: \mathcal{P}^p\rightarrow \mathcal{P}^{p-1}.$$
where $X_{P}=(\partial_yP,-\partial_xP,{\bf i}\kappa \partial_zP,-{\bf i}\kappa \partial_{\bar{z}}P)$ with ${\bf i}\kappa \partial_zP=({\bf i}j \partial_{z_j}P)_{j\in\mathbb{Z}_*}$.

\subsubsection{Some notations}
 Let $\Vert \cdot\Vert$ be an operator norm or $l^2$ norm. The notation $|\cdot |$ denotes sup norm. Let $l_p^2$ be the Hilbert space of all complex sequences $z=(z_j)_{j\in\mathbb{Z}_*}$ with
\begin{equation}\label{lpp2}
\Vert z\Vert_{p}^2=\sum_{j\in\mathbb{Z}_*}|j|^{2p}|z_j|^2<\infty.
\end{equation}

\begin{defn}
Let
$$D(s)=\{x\in\mathbb{T}^n|\, |{\rm Im}\,x|\le s\}.$$
 If $W(x;\xi):D(s)\times \Pi\rightarrow \mathbb{C}$ is analytic in $x\in D(s)$ and $C^1$ (in Whitney's sense) in $\xi\in\Pi$ with the Fourier series
$W(x;\xi)=\sum_{k\in\mathbb{Z}^n}W(k;\xi)e^{{\bf i}\langle k,x\rangle}$,
define the norm
$$\Vert W(x;\xi)\Vert_{D(s)\times \Pi}=\sup_{\xi\in\Pi}\sum_{k\in\mathbb{Z}^n}(|W(k;\xi)|+|\partial_\xi W(k;\xi)|)e^{|k|s}.$$
\end{defn}

\begin{defn}
Let
$$D(s,r,r)=\{(x,y)\in\mathcal{P}^p|\ |{\rm Im}\,x|\le s,\, |y|\le r^2, \Vert z\Vert_p+\Vert \bar{z}\Vert_p\le r\}.$$
Consider a function
$$ W_h(x,y,z,\bar{z};\xi)=\sum_{\substack{\alpha\in \mathbb{N}^n,\mu,\gamma\in\mathbb{N}^{\mathbb{Z}_*}\\2|\alpha|+|\mu+\gamma|=h}} W^{\alpha\mu\gamma}(x;\xi)y^\alpha z^\mu\bar{z}^\gamma: D(s,r,r)\times\Pi\rightarrow \mathbb{R},$$
which is analytic in $(x,y,z,\bar{z})\in D(s,r,r)$ and $C^1$ in $\xi\in\Pi$. Then we define the norms by
$$||| W_h(x,y,z,\bar{z};\xi)|||_{D(s,r,r)\times \Pi}:%=|||\lfloor W_h(x,y,z,\bar{z};\xi)\rceil_{D(s)\times \Pi}(y,z,\bar{z})|||$$
%$$
=r^h\sup_{\substack{\alpha\in \mathbb{N}^n,\mu,\gamma\in\mathbb{N}^{\mathbb{Z}_*}\\2|\alpha|+|\mu+\gamma|=h}}\Vert W^{\alpha\mu\gamma}(x;\xi)\Vert_{D(s)\times \Pi},$$
$$||| W_h(x,y,z,\bar{z};\xi)|||^*_{D(s,r,r)\times \Pi}:%=|||\lfloor W(x,y,z,\bar{z};\xi)\rceil_{D(s)\times \Pi}(z,\bar{z})|||^*$$
%$$
=r^h\sup_{\substack{\alpha\in\mathbb{N}^n,\mu,\gamma\in\mathbb{N}^{\mathbb{Z}_*}\\2|\alpha|+|\mu+\gamma|=h}}M_{\mu\gamma}\Vert W^{\alpha\mu\gamma}(x;\xi)\Vert_{D(s)\times \Pi},$$
\end{defn}
where $$M_{\mu\gamma}=\sup\left\{|j|\,|\, \mu_j\neq 0 \,{\rm or}\, \gamma_j\neq 0\right\}.$$

\begin{defn}\label{defnd}
Consider a Hamiltonian $W(x,y,z,\bar{z};\xi)=\sum_{h\geq 0}W_h(x,y,z,\bar{z};\xi)$ is analytic in $(x,y,z,\bar{z})\in D(s,r,r)$ and $C^1$ in $\xi\in\Pi$, where
$$ W_h(x,y,z,\bar{z};\xi)=\sum_{\substack{\alpha\in \mathbb{N}^n,\mu,\gamma\in\mathbb{N}^{\mathbb{Z}_*}\\2|\alpha|+|\mu+\gamma|=h}} W^{\alpha\mu\gamma}(x;\xi)y^\alpha z^\mu\bar{z}^\gamma.$$
Then we define the norms of the function $W$ by
$$||| W(x,y,z,\bar{z};\xi)|||_{D(s,r,r)\times \Pi}:=\sum_{h\geq 0}||| W_h(x,y,z,\bar{z};\xi)|||_{D(s,r,r)\times \Pi},$$
$$||| W(x,y,z,\bar{z};\xi)|||^*_{D(s,r,r)\times \Pi}:=\sum_{h\geq 0}||| W_h(x,y,z,\bar{z};\xi)|||^*_{D(s,r,r)\times \Pi}.$$
We say $W$ belongs to $\mathcal{H}$ if $||| W(x,y,z,\bar{z};\xi)|||_{D(s,r,r)\times \Pi}<\infty$. Similarly, we say $W$ belongs to $\mathcal{H}^*$ if $||| W(x,y,z,\bar{z};\xi)|||^*_{D(s,r,r)\times \Pi}<\infty$.
\end{defn}

\begin{defn}
For a function $W(x,y,z,\bar{z};\xi)$ is analytic in $(x,y,z,\bar{z})$ on $D(s,r,r)$ and $C^0$ in $\xi\in\Pi$.
Define $$|||W|||^0_{D(s,r,r)\times\Pi}=:\sum_{h\geq 0}||| W_h(x,y,z,\bar{z};\xi)|||^0_{D(s,r,r)\times \Pi},$$
$$|||W|||^{0*}_{D(s,r,r)\times\Pi}=:\sum_{h\geq 0}||| W_h(x,y,z,\bar{z};\xi)|||^{0*}_{D(s,r,r)\times \Pi},$$
with $$||| W_h(x,y,z,\bar{z};\xi)|||^{0}_{D(s,r,r)\times \Pi}=r^h\sup_{\substack{\alpha\in \mathbb{N}^n,\mu,\gamma\in\mathbb{N}^{\mathbb{Z}_*}\\2|\alpha|+|\mu+\gamma|=h}}\Vert W^{\alpha\mu\gamma}(x;\xi)\Vert^{0}_{D(s)\times \Pi},$$
$$||| W_h(x,y,z,\bar{z};\xi)|||^{0*}_{D(s,r,r)\times \Pi}=r^h\sup_{\substack{\alpha\in \mathbb{N}^n,\mu,\gamma\in\mathbb{N}^{\mathbb{Z}_*}\\2|\alpha|+|\mu+\gamma|=h}}M_{\mu\gamma}\Vert W^{\alpha\mu\gamma}(x;\xi)\Vert^{0}_{D(s)\times \Pi},$$
where
$$\Vert W^{\alpha\mu\gamma}(x;\xi)\Vert^0_{D(s)\times \Pi}=\sup_{\xi\in\Pi}\sum_{k\in\mathbb{Z}^n}|W(k;\xi)|e^{|k|s}.$$
\end{defn}

\begin{defn}\label{vf}
Consider a function
$$ W(x,y,z,\bar{z};\xi)=\sum_{\mu,\gamma\in\mathbb{N}^{\mathbb{Z}_*}} W^{\mu\gamma}(x,y;\xi)z^\mu\bar{z}^\gamma: D(s,r,r)\times\Pi\rightarrow \mathbb{R}.$$
%belong to the space $H^*$.
The corresponding vector field for the function $W$ is defined by
$$X_W=(\partial_yW,-\partial_xW,{\bf i}\kappa \partial_zW,-{\bf i}\kappa \partial_{\bar{z}}W),$$
 where ${\bf i}\kappa \partial_zW=({\bf i}j \partial_{z_j}W)_{j\in\mathbb{Z}_*}$. The norm for the vector field is defined by
$$\Vert X_W\Vert_{p,D(s,r,r)\times \Pi}=\Vert \partial_yW\Vert+\frac{1}{r^2}\Vert \partial_xW\Vert+\frac{1}{r}\Vert {\bf i}\kappa  \partial_zW\Vert_{p}+\frac{1}{r}\Vert {\bf i}\kappa \partial_{\bar{z}}W \Vert_{p}.$$
\end{defn}

\subsection{The abstract results}
Consider the nearly integrable Hamiltonian
\begin{equation}\label{hhh}
H=N(y,z,\bar{z};\xi)+P(x,y,z,\bar{z};\xi),
\end{equation}
$$N(y,z,\bar{z};\xi)=\sum_{1\le j\le n}\omega_j(\xi)y_j+\sum_{j\in \mathbb{Z}_*}\frac{\Omega_j(\xi)}{j}z_j\bar{z}_j,$$
where
$$\omega(\xi)=(\omega_1(\xi),\cdots,\omega_n(\xi)),\quad \left(\Omega_j(\xi)\right)_{j\in \mathbb{Z}_*},$$
are called tangent and normal frequencies, respectively. Denote $\Omega(\xi)=(\Omega_j(\xi))_{j\in\mathbb{Z}_*}$, and where the function
$$P(x,y,z,\bar{z};\xi)%=\sum_{\beta,\gamma\in \mathbb{N}^{\mathbb{Z}_*}}P^{\beta\gamma}(x,y;\xi)z^\beta\bar{z}^\gamma
%$$ $$
=\sum_{\alpha\in\mathbb{N}^n,\beta,\gamma\in \mathbb{N}^{\mathbb{Z}_*}}P^{\alpha\beta\gamma}(x;\xi)y^{\alpha}z^\beta\bar{z}^\gamma$$
is the perturbation. Hamiltonian \eqref{hhh} is equipped with the symplectic structure
$$\sum_{1\le i\le n}dy_i\wedge dx_i-{\bf i}\sum_{j\in \mathbb{Z}_*} j^{-1} d\bar{z}_j\wedge dz_j.$$

Suppose the following conditions are satisfied:

{\bf Assumption A: Frequency Asymptotics}. The tangent frequencies $\omega_{i}'s$ are of the form
$$\omega_i(\xi)=j_i^2+\xi_{j_i},\quad j_{i} \in J_n=\{j_{1},\cdots,j_{n}\},\ \ i=1,\cdots,n $$
and  the normal frequencies $\Omega_{j}'s$ are of the form
$$\Omega_j(\xi)=j^2+\xi_j,\quad j\in\mathbb{Z}_*.$$
%{\bf Assumption B: Twist conditions}.
%$$\partial_{\xi_j}\omega_i(\xi)=\delta_{ij},\quad \partial_{\xi_j}(\Omega_{j^\prime}(\xi))=\delta_{jj^{\prime}},\quad 1\le i\le n,\ j,j^\prime\in\mathbb{Z}_*.$$

{\bf Assumption B: Momentum conservation}. The Hamiltonian $H$ has the momentum conservation , i.e., for each monomial
$$e^{{\bf i}\langle k,x\rangle}\prod_{1\le j\le n}y_j^{\alpha_j}\prod_{j\in\mathbb{Z}_*}z_j^{\beta_j}\bar{z}_j^{\gamma_j},\quad k\in\mathbb{Z}^n,\,\alpha_j, \beta_j,\gamma_j\in\mathbb{N},$$
of $H(x,y,z,\bar{z};\xi)$,
$$M(k,\beta,\gamma)=-\sum_{1\le i\le n}k_i j_i+\sum_{j\in\mathbb{Z}_*}(\beta_j-\gamma_j)j=0.$$

%{\bf Assumption C: Real index symmetry}. The coefficients of the perturbation $P$ satisfy: for real variables $x$
%$$P^{\alpha\beta\gamma}(x)=P^{\alpha\gamma\beta}(x)\in\mathbb{R},\quad \forall \beta,\gamma\in \mathbb{N}^{\mathbb{Z}_*}.$$

{\bf Assumption C: The independence on $\xi$ of P}. In the initial step, the perturbation $P$ is independent of $\xi$.

\begin{The}[KAM theorem or Normal form of order $2$]\label{normaltheorem}
Around KAM tori, suppose that $H=N+P$
is a perturbation of the integrable Hamiltonian $N(y,z,\bar{z};\xi)$ on the domain $D(s_0,r_0,r_0)$ with some $s_0, r_0\in(0,1]$ depending on parameter $\xi\in\Pi$. In addition, the tangent frequency and normal frequency satisfy the {\bf assumptions  A}, {\bf B} and {\bf C}. Moreover, the perturbation belongs to $\mathcal{H}$ with the estimate
$$\epsilon:=|||P|||_{D(s_0,r_0,r_0)\times \Pi}\le \eta^{12}\varepsilon, \quad {\rm for}\,\  {\rm some}\,\  \eta\in(0,1),$$
where $\varepsilon>0$ is a sufficiently small  constant depending on $s_0, r_0$ and $n$. In addition, for fixed $\xi \in \Pi, P=P(x,y,z,\overline{z}; \xi)$ is analytic in $(x,y,z,\overline{z}) \in D(s_{0},r_{0},r_{0})$, and it is real when $x,y$ is real and $\overline{z}$ is regarded as the complex conjugate of $z$. (We briefly call $P$ is real analytic in $(x,y,z,\overline{z}) \in D(s_{0},r_{0},r_{0})$.)

Then there is a subset $\Pi_\eta\in\Pi$ with the estimate
$${\rm Meas}\ \Pi_\eta\geq ({\rm Meas}\ \Pi)(1-O(\eta)).$$
For each $\xi\in\Pi_\eta$, there is a symplectic map
$$\Psi:D(s_0/2,r_0/2,r_0/2)\rightarrow D(s_0,r_0,r_0),$$
such that
$$H\circ \Psi=\breve{N}(x,y,z,\bar{z};\xi)+\breve{P}(x,y,z,\bar{z};\xi),$$
where
$$\breve{N}(x,y,z,\bar{z};\xi)=\sum_{j=1}^n\breve{\omega}_j(\xi)y_j+\sum_{j\in \mathbb{Z}_*}\frac{\breve{\Omega}_j(x;\xi)}{j}z_j\bar{z}_j,$$
%$$=\sum_{j=1}^n\breve{\omega}_j(\xi)y_j+\sum_{|j|\le \mathcal{N}}\frac{\breve{\Omega}_j(\xi)}{j}z_j\bar{z}_j+\sum_{j\in \mathbb{Z}_*,|j|\geq \mathcal{N}+1}\frac{\breve{\Omega}_j(x;\xi)}{j}z_j\bar{z}_j$$
and
$$\breve{P}(x,y,z,\bar{z};\xi)=\sum_{\substack{\alpha\in\mathbb{N}^n,\beta,\gamma\in \mathbb{N}^{\mathbb{Z}^*},\\ 2|\alpha|+|\beta+\gamma|\geq 3}}\breve{P}^{\alpha\beta\gamma}(x;\xi)y^{\alpha}z^{\beta}\bar{z}^{\gamma}.$$
In addition, the following estimate holds:

(1) for each $\xi\in\Pi_\eta$, the symplectic map $\Psi$ satisfies
\begin{equation}\label{lho1}
\Vert \Psi-id\Vert_{p,D(s_0/2,r_0/2,r_0/2)}\le c\varepsilon \eta^6,
\end{equation}
where
$$\Vert \Psi-id\Vert_{p,D(s_0/2,r_0/2,r_0/2)}=\sup_{w\in D(s_0/2,r_0/2,r_0/2)}\Vert (\Psi-id)w\Vert_{\mathcal{P}^p,D(s_0/2,r_0/2,r_0/2)},$$
and
$$\Vert w\Vert_{\mathcal{P}^p,D(s_0,r_0,r_0)}=\Vert w_x\Vert+\frac{1}{r_0^2}\Vert w_y\Vert+\frac{1}{r_0}\Vert w_z\Vert_p+\frac{1}{r_0}\Vert w_{\bar{z}}\Vert_p,$$
for each $w=(w_x,w_y,w_z,w_{\bar{z}})\in D(s_0,r_0,r_0)$; moreover,
$$|| D\Psi-id||_{p,D(s_0/2,r_0/2,r_0/2)}\le c\varepsilon \eta^6,$$
where
$$|| D\Psi-id||_{p,D(s_0/2,r_0/2,r_0/2)}=\sup_{0\neq w\in D(s_0/2,r_0/2,r_0/2)}\frac{\Vert (D\Psi-Id)w\Vert_{\mathcal{P}^p,D(s_0/2,r_0/2,r_0/2)}}{\Vert w\Vert_{\mathcal{P}^p,D(s_0/2,r_0/2,r_0/2)}};$$
(2) the frequencies $\breve{\omega}$ and $\breve{\Omega}$ satisfy
$$\Vert \breve{\omega}-\omega\Vert_{\Pi_\eta}\le  c\varepsilon \eta^6,$$
and
\begin{equation}\label{huomom}
\Vert \breve{\Omega}-\Omega\Vert_{-1,D(s_0/2)\times\Pi_\eta}:=\sup_{j\in\mathbb{Z}^*}\Vert j^{-1}  (\breve{\Omega}_j-\Omega_j)\Vert_{D(s_0/2)\times\Pi_\eta}\le  c\varepsilon \eta^6;
\end{equation}
(3) the new Hamiltonian $\breve{P}$ satisfies
$$|||\breve{P}|||_{D(s_0/2,r_0/2,r_0/2)\times \Pi_\eta}\le \epsilon(1+c\eta^6\varepsilon).$$
\end{The}
\begin{Rem}\label{rem2.7}
From Theorem \ref{normaltheorem}, it is easy to see that for each $\xi\in \Pi_\eta$ there is an analytic embedding invariant torus $\mathcal{T}^0_\xi=\mathbb{T}^n\times \{0\}\times\{0\}\times\{0\}$ with frequency $\breve{\omega}(\xi)$ for Hamiltonian $\breve H$ and $\mathcal{T}_\xi=\Psi^{-1}\mathcal{T}^0_\xi$ is an analytic embedding invariant torus for the original Hamiltonian $H$. In addition, the invariant torus is linearly stable. In Theorem \ref{lthe} below, it is shown that the invariant torus is long time stable.
\end{Rem}

\begin{The}[The long time stability of KAM tori]\label{lthe}
Given any $0\le \mathcal{M}\le (2c\eta^6 \epsilon)^{-1/2}$ and $p\geq 8C_0(\mathcal{M}+7)^4+1$, there exists a small positive $\delta_0$ depending on $s_0, r_0, n$ and $\mathcal{M}$ and a subset $\Pi_{\acute{\eta}}\subset \Pi_{\eta}$ satisfying
$${\rm Meas}\ \Pi_{\acute{\eta}}\geq ({\rm Meas}\ \Pi_{\eta})(1-O(\acute{\eta})),$$
where $\acute{\eta}$ is some constant in $(0,1)$. For each $\xi\in\Pi_{\acute{\eta}}$, the KAM tori $\mathcal{T}_\xi$ is stable in long time, i.e., if $w(t)$ is a solution of Hamiltonian vector field $X_H$ with initial datum $w(0)$
satisfying
$$d_p(w(0),\mathcal{T}_\xi)\le \delta,$$
then
$$d_{p/2}(w(t),\mathcal{T}_\xi)\le 2\delta,\quad \forall \  |t|\le \delta^{-\mathcal{M}/4}.$$
%where the distance $d_p(w,v)$ between any two points
\end{The}

\begin{Rem}
In Theorem \ref{normaltheorem} and \ref{lthe}, the measure is in the sense of Kolmogorov, since the parameter set $\Pi\in \mathbb{R}^{\mathbb{Z}\setminus\{0\}}$ is infinite dimension. However, in the present paper we can understand it as Lebesgue measure. Indeed, if we write
$$\xi=(\xi^n,\xi^{C(\mathcal{M},\mathcal{N})},\xi^{\mathbb{Z}\setminus\{0\}})\in \mathbb{R}^n\times \mathbb{R}^{C(\mathcal{M},\mathcal{N})}\times l^2$$
where $C(\mathcal{M},\mathcal{N})$ is a constant depending on some large $\mathcal{M},\mathcal{N}$ which are defined in Theorem \ref{lthe}. In the proof of Theorem \ref{normaltheorem} and \ref{lthe}, it is enough to treat $\xi^n$ and $(\xi^n,\xi^{C(\mathcal{M},\mathcal{N})})$ as parameters, respectively. Therefore, we only need finite parameters in the whole proof of Theorem \ref{theorem1.1}.
\end{Rem}

\section{Proof of Theorem \ref{normaltheorem}}
\subsection{Derivation of homological equations}

We suppose at $\nu$-th step of Newton iteration, the Hamiltonian reads
$$H^{\nu}(x,y,z,\bar{z};\xi)=N^{\nu}(x,y,z,\bar{z};\xi)+P^{\nu}(x,y,z,\bar{z};\xi),$$
with
$$N^{\nu}(x,y,z,\bar{z};\xi)=\sum_{1\le j\le n}\omega^{\nu}_j(\xi)y_j+\sum_{j\in\mathbb{Z}_* }\frac{\Omega_j^{\nu}(x;\xi)}{j}z_j\bar{z}_j.$$
%$$=\sum_{j=1}^n\breve{\omega}_j(\xi)y_j+\sum_{|j|\le \mathcal{N}}\frac{\breve{\Omega}_j(\xi)}{j}z_j\bar{z}_j+\sum_{j\in \mathbb{Z}_*,|j|\geq \mathcal{N}+1}\frac{\breve{\Omega}_j(x;\xi)}{j}z_j\bar{z}_j.$$
Denote $P^{\nu}(x,y,z,\bar{z};\xi)=P^{\nu,low}(x,y,z,\bar{z};\xi)+P^{\nu,high}(x,y,z,\bar{z};\xi)$ with
\begin{equation}\label{*1}
P^{\nu,low}(x,y,z,\bar{z};\xi)=\sum_{\substack{\alpha\in\mathbb{N}^n,\beta,\gamma\in \mathbb{N}^{\mathbb{Z}_*},\\ 2|\alpha|+|\beta+\gamma|\le 2}}P^{\nu,\alpha\beta\gamma}(x;\xi)y^{\alpha}z^{\beta}\bar{z}^{\gamma},
\end{equation}
\begin{equation}\label{*2}
P^{\nu,high}(x,y,z,\bar{z};\xi)=\sum_{\substack{\alpha\in\mathbb{N}^n,\beta,\gamma\in \mathbb{N}^{\mathbb{Z}_*},\\ 2|\alpha|+|\beta+\gamma|\geq 3}}P^{\nu,\alpha\beta\gamma}(x;\xi)y^{\alpha}z^{\beta}\bar{z}^{\gamma}.
\end{equation}
Suppose the following estimates hold
\begin{equation}\label{hupl}
|||P^{\nu,low}|||_{D(s,r,r)\times\Pi}\le \epsilon,\quad |||\partial_{\xi_a}P^{\nu,low}|||^0_{D(s,r,r)\times \Pi}\le \frac{\epsilon}{|a|},
\end{equation}
and
\begin{equation}
|||P^{\nu,high}|||_{D(s,r,r)\times\Pi}\le 1,\quad |||\partial_{\xi_a}P^{\nu,high}|||^0_{D(s,r,r)\times\Pi}\le \frac{1}{|a|},\label{huph}
\end{equation}
where $a\in\mathbb{Z}_*$.

Write
$$\Omega^{\nu}(x;\xi)=[\Omega]^{\nu}(\xi)+\tilde{\Omega}^{\nu}(x;\xi),$$
where $[\Omega]^{\nu}$ is the mean value of the variable $x$ over $\mathbb{T}^n$.
For the tangent and normal frequencies, we suppose the following non-resonant conditions hold:
\begin{align}
|\langle k,\omega^{\nu}(\xi)\rangle+\langle l,[\Omega]^{\nu}(\xi)\rangle|\geq\frac{M_{l}\eta}{(1+|k|)^\tau},\quad &|k|+|l|\neq 0,\,|l|\le2,%\\
 %|\langle l,[\Omega]^{\nu}(\xi)\rangle|\geq\frac{M_l\eta}{|k|^\tau},\quad &|l|\le2,
%|\tilde{\Omega}|_{D(s)}\le \gamma_0| j|,\quad &j\in \mathbb{Z}_*.\label{to}
\end{align}
where $M_l:=\max\{|j|\,|\, |l_j|\neq 0\}$, and
$$\tilde{\Omega}^{\nu}_j(x;\xi)=\sum_{i=0}^{\nu} \tilde{\Omega}^{\nu,i}_j(x;\xi),\quad j\in\mathbb{Z}_*$$
with
\begin{equation}\label{omegai}
\Vert \tilde{\Omega}^{\nu,i}(x;\xi)\Vert_{-1,D(s_i)\times \Pi_i}:=\sup_{j\in\mathbb{Z}_*}\Vert j^{-1}\tilde{\Omega}^{\nu,i}_j(x;\xi)\Vert_{D(s_i)\times \Pi_i}\le \epsilon_i^{1-}.
\end{equation}
In addition, we assume,
$$\left|\frac{\partial \omega^\nu_j}{\partial \xi_a}\right|\le \frac{1}{|a|}\epsilon_0,\quad \left|\frac{\partial \Omega^\nu_j}{\partial \xi_a}\right|\le \delta_{ja}+\frac{|j|}{|a|}\epsilon_0,$$
for $a\in \mathbb{Z}_*$.

In the following, in order to avoid a fluid of constants, we will use the notations ``$\lessdot$" to denote ``$<c$" with constant $c$ independent of iterations. Moreover, the notations without subscript ``$+$" represent the quantities at the $\nu$-th step, and those with subscript ``$+$" represent the corresponding ones at $(\nu+1)$-the step.

The symplectic transformation $\Phi$ is obtained as the time-$1$-map $X_F^t|_{t=1}$ of a Hamiltonian vector field $X_F=(\partial_yF,-\partial_xF,{\bf i}\kappa \partial_zF,-{\bf i}\kappa \partial_{\bar{z}}F)$ with ${\bf i}\kappa \partial_zF=({\bf i}j \partial_{z_j}F)_{j\in\mathbb{Z}_*}$,  where $F(x,y,z,\bar{z};\xi)$ is of the form:
$$F(x,y,z,\bar{z};\xi)=\sum_{2|\alpha|+|\beta+\gamma|\le 2}F^{\alpha\beta\gamma}(x;\xi)y^{\alpha}z^{\beta}\bar{z}^{\gamma}.$$
For simplicity, here and in the following, we omit $\alpha\in\mathbb{N}^n,\beta,\gamma\in \mathbb{N}^{\mathbb{Z}_*}$ in the summation notation. From Taylor's formula, we have
\begin{align*}
H_+=&H\circ \Phi=H\circ X_F^1\\
=&N+\{N,F\}+\int_0^1(1-t)\{\{N,F\},F\}\circ X_F^tdt\\
&+P^{low}+\int_0^1 \{P^{low},F\}\circ X_F^tdt\\
&+P^{high}+\{P^{high},F\}+\int_0^1(1-t)\{\{P^{high},F\},F\}\circ X_F^tdt.%\\
%=&N+\{N,F\}+\int_0^1(1-t)\{\{N,F\},F\}\circ X_F^tdt\\
%&+P_{\le }^{low}+(P^{low}-P_{\le }^{low})+\int_0^1 \{P^{low},F\}\circ X_F^tdt\\
%&+P^{high}+\{P^{high},F\}^{low}_{\le }+(\{P^{high},F\}-\{P^{high},F\}^{low}_{\le })\\
%&+\int_0^1(1-t)\{\{P^{high},F\},F\}\circ X_F^tdt,
\end{align*}
%where $P_{\le }^{low}$ denotes the Fourier series of the function $P^{low}$ truncating at $K$-th order, where the constant $K$ will be determined later. The definition for $\{P^{high},F\}_{\le }$ is similar.
We need to solve the following homological equation
\begin{equation}\label{homological}
\{N,F\}+P^{low}+\{P^{high},F\}^{low}=\hat{N}+\hat{P}.
\end{equation}
where $\{P^{hign},F\}^{low}$ is defined similarly as in \eqref{*1}.
After we solve this homological equation, the new normal form $N_+$ and new perturbation $P_+$ have the following form, respectively,
\begin{align}
N_+=&N+\hat{N},\nonumber\\
P_+=&\hat{P}+\int_0^1(1-t)\{\{N,F\},F\}\circ X_F^tdt\label{phigh}\\
%&+(P^{low}-P_{\le }^{low})+(\{P^{high},F\}^{low}-\{P^{high},F\}^{low}_{\le })\\
&+ P^{high}+\{P^{high},F\}^{high}+\int_0^1 \{P^{low},F\}\circ X_F^tdt\nonumber\\
&+\int_0^1(1-t)\{\{P^{high},F\},F\}\circ X_F^tdt.\nonumber
\end{align}
where $\{P^{high},F\}^{high}$ is defined similarly as in \eqref{*2}.

\subsection{Solving the homological equations}
We firstly rewrite $P^{low}$ as
$$P^{low}=P^x+P^y+P^1+P^2,$$
with
\begin{align*}
P^x=&P^x(x;\xi)=\sum_{2|\alpha|+|\beta+\gamma|=0}P^{\alpha\beta\gamma}(x;\xi),\\
P^y=&\sum_{1\le j\le n}P^{y_j}(x;\xi)y_j=\sum_{|\alpha|=1,|\beta+\gamma|=0}P^{\alpha\beta\gamma}(x;\xi)y^{\alpha},\\
P^1=&\sum_{j\in\mathbb{Z}_*}\left(P^{z_j}(x;\xi)z_j+P^{\bar{z}_j}(x;\xi)\bar{z}_j\right)
=\sum_{ |\alpha|=0,|\beta+\gamma|=1}P^{\alpha\beta\gamma}(x;\xi)z^{\beta}\bar{z}^{\gamma},\\
P^2=&\sum_{i,j\in\mathbb{Z}_*}\left(P^{z_iz_j}(x;\xi)z_iz_j+P^{z_i\bar{z}_j}(x;\xi)z_i\bar{z}_j+P^{\bar{z}_i\bar{z}_j}(x;\xi)\bar{z}_i\bar{z}_j\right)\\
=&\sum_{ |\alpha|=0,|\beta+\gamma|=2}P^{\alpha\beta\gamma}(x;\xi)z^{\beta}\bar{z}^{\gamma}.
\end{align*}
%where $P^{\alpha\beta\gamma}(x;\xi)=\sum_{|k|\le K}P^{\alpha\beta\gamma}(k;\xi)e^{{\bf i}\langle k,x\rangle}$, the constant $K$ will be defined later.
The function $F$ is written as the same form as $P^{low}$:
$$F=F^x+F^y+F^1+F^2,$$
where $F^x,F^y,F^1,F^2$ are similar to $P^x,P^y,P^1,P^2$. In order to write the term $\{P^{high},F\}^{low}$ explicitly, we rewrite
$$P^{high}=\sum_{j=0}^{4}P^{(j)},$$
where
\begin{align*}
P^{(0)}=&\sum_{ |\alpha|=2,|\beta+\gamma|=0}P^{\alpha\beta\gamma}(x;\xi)y^{\alpha},\\
P^{(1)}=&\sum_{ |\alpha|=1,|\beta+\gamma|=1}P^{\alpha\beta\gamma}(x;\xi)y^{\alpha}z^{\beta}\bar{z}^{\gamma},\\
P^{(2)}=&\sum_{ |\alpha|=1,|\beta+\gamma|=2}P^{\alpha\beta\gamma}(x;\xi)y^{\alpha}z^{\beta}\bar{z}^{\gamma},\\
P^{(3)}=&\sum_{ |\alpha|=0,|\beta+\gamma|=3}P^{\alpha\beta\gamma}(x;\xi)y^{\alpha}z^{\beta}\bar{z}^{\gamma},\\
P^{(4)}=&\sum_{ 2|\alpha|+|\beta+\gamma|\geq 5\, {\rm or}\, |\beta+\gamma|\geq 4}P^{\alpha\beta\gamma}(x;\xi)y^{\alpha}z^{\beta}\bar{z}^{\gamma}.
\end{align*}
By the definition of Poisson bracket with the symplectic structure \eqref{symplectic structure}, we have
$$\{P^{high},F\}^{low}=\{P^{high},F\}^{y}+\{P^{high},F\}^{1}+\{P^{high},F\}^{2},$$
where $\{P^{high},F\}^{y},\{P^{high},F\}^{1},\{P^{high},F\}^{2}$ are similar to $P^y,P^1,P^2$ and taken on the following explicit form
\begin{align}
\{P^{high},F\}^{y}=&\sum_{1\le j\le n}P_{y_j}^{(0)}F^{x}_{x_j}+{\bf i}\sum_{j\in\mathbb{Z}_*}j(P^{(1)}_{\bar{q_j}}F^1_{q_j}-P^{(1)}_{q_j}F^1_{\bar{q}_j}),\quad
\{P^{high},F\}^{1}=\sum_{1\le j\le n}P_{y_j}^{(1)}F^{x}_{x_j},\\
\{P^{high},F\}^{2}=&\sum_{1\le j\le n}(P_{y_j}^{(1)}F^{1}_{x_j}+P_{y_j}^{(2)}F^{x}_{x_j})+{\bf i}\sum_{j\in\mathbb{Z}_*}j(P^{(3)}_{\bar{q_j}}F^1_{q_j}-P^{(3)}_{{q_j}}F^1_{\bar{q}_j}).\label{high2}
\end{align}
Denote $\partial_{\omega}=\omega\cdot \partial_x$. The homological equation \eqref{homological} becomes
\begin{align}
&\partial_\omega F^x=P^x-\hat{N}^x-\hat{P}^x,\label{2.9}\\
&\partial_\omega F^{y_j}=P^{y_j}+\{P^{high},F\}^{y_j}-\hat{N}^{y_j}-\hat{P}^{y_j},\quad 1\le j\le n,\\
&(\partial_\omega +{\bf i}\Omega_j )F^{z_j}=P^{z_j}+\{P^{high},F\}^{z_j}-\hat{P}^{z_j},\quad j\in\mathbb{Z}_*,\\
&(\partial_\omega-{\bf i}\Omega_j) F^{\bar{z}_j}=P^{\bar{z}_j}+\{P^{high},F\}^{\bar{z}_j}-\hat{P}^{\bar{z}_j},\quad j\in\mathbb{Z}_*,\\
&(\partial_\omega+{\bf i}(\Omega_i+\Omega_j)) F^{z_iz_j}=P^{z_iz_j}+\{P^{high},F\}^{z_iz_j}-\hat{P}^{z_iz_j},\quad i,j\in\mathbb{Z}_*,\\
&(\partial_\omega-{\bf i}(\Omega_i+\Omega_j)) F^{\bar{z}_i\bar{z}_j}=P^{\bar{z}_i\bar{z}_j}+\{P^{high},F\}^{\bar{z}_i\bar{z}_j}-\hat{P}^{\bar{z}_i\bar{z}_j},\quad i,j\in\mathbb{Z}_*,\\
&(\partial_\omega+{\bf i}(\Omega_i-\Omega_j)) F^{z_i\bar{z}_j}=P^{z_i\bar{z}_j}+\{P^{high},F\}^{z_i\bar{z}_j}-\hat{P}^{z_i\bar{z}_j},\quad i,j\in\mathbb{Z}_*,\,\,i\neq  j,\label{zbarz}
%&\partial_\omega F^{z_j\bar{z}_j}=P^{z_j\bar{z}_j}+\{P^{high},F\}^{z_j\bar{z}_j}-\hat{N}^{z_j\bar{z}_j},\quad |j|\le\mathcal{N},
\end{align}
where
\begin{align*}
&\hat{N}^x=[P^x],\quad  \hat{N}^{y_j}=[P^{y_j}+\{P^{high},F\}^{y_j}],\\
&\hat{P}^x=\partial_\omega (1-\Gamma_K)G^x,\quad \hat{P}^{y_j}=\partial_\omega (1-\Gamma_K)G^{y_j},\\
&\hat{P}^{z_j}=(\partial_\omega +{\bf i}\Omega_j )(1-\Gamma_K)G^{z_j},\\
&\hat{P}^{z_j}=(\partial_\omega-{\bf i}\Omega_j) (1-\Gamma_K)G^{\bar{z}_j},\\
&\hat{P}^{z_iz_j}=(\partial_\omega+{\bf i}(\Omega_i+\Omega_j)) (1-\Gamma_K)G^{z_iz_j},\\
&\hat{P}^{\bar{z}_i\bar{z}_j}=(\partial_\omega-{\bf i}(\Omega_i+\Omega_j))(1-\Gamma_K) G^{\bar{z}_i\bar{z}_j},\\
&\hat{P}^{z_i\bar{z}_j}=(\partial_\omega+{\bf i}(\Omega_i-\Omega_j)) (1-\Gamma_K)G^{z_i\bar{z}_j}.
\end{align*}
%$$\hat{N}^{z_j\bar{z}_j}=[P^{z_j\bar{z}_j}+\{P^{high},F\}^{z_j\bar{z}_j}],\quad |j|\le \mathcal{N}$$
Here, $G$ is the solution to the homological equation \eqref{homological} without the term $\hat{P}$, and $\Gamma_K G$ denotes the Fourier series of the function $G$ truncated at $K$-th order, where the constant $K$ will be specified later.  Therefore, we obtain
$$\hat{P}=\{N,(1-\Gamma_K) G\}.$$

The terms $P^{z_j\bar{z}_j}+\{P^{high},F\}^{z_j\bar{z}_j}$ ($j\in\mathbb{Z}_*$) %and $P^{z_j\bar{z}_{-j}}+\{P^{high},F\}^{z_j\bar{z}_{-j}}$
can not be eliminated, %Since $z_j=-z_{-j}$, thus we have
%$$P^{z_j\bar{z}_{-j}}z_j\bar{z}_{-j}+\{P^{high},F\}^{z_j\bar{z}_{-j}}z_j\bar{z}_{-j}=-P^{z_j\bar{z}_{-j}}z_j\bar{z}_{j}-\{P^{high},F\}^{z_j\bar{z}_{-j}}z_j\bar{z}_{j}$$
which can be put into the new normal form.
Finally, we have
$$N_+=N+\hat{N}=\sum_{1\le j\le n}\omega_{+j}y_j+\sum_{\in\mathbb{Z}_*}\frac{\Omega_{+j}(x;\xi)}{j}z_j\bar{z}_j$$
with
\begin{align}
\omega_{+j}(\xi)=&\omega_j+[P^{y_j}+\{P^{high},F\}^{y_j}],\quad 1\le j\le n,\label{o+}\\
%\Omega_{+j}(\xi)=&\Omega_j+j[P^{z_j\bar{z}_j}+\{P^{high},F\}^{z_j\bar{z}_j}],\quad |j|\le\mathcal{N},\\
\Omega_{+j}(x;\xi)=&\Omega_j+j\big(P^{z_j\bar{z}_j}+\{P^{high},F\}^{z_j\bar{z}_j}+\left\langle \partial_{x}\Omega_j/j,F^{y}\right\rangle\big),\quad j\in\mathbb{Z}_*.%\nonumber\\
%&-P^{z_j\bar{z}_{-j}}-\{P^{high},F\}^{z_j\bar{z}_{-j}}]\label{O+}
\end{align}
It means
\begin{align*}
\hat{N}=&\sum_{1\le j\le n}[P^{y_j}+\{P^{high},F\}^{y_j}]y_j\\
%&+\sum_{|j|\le\mathcal{N}}[P^{z_j\bar{z}_j}+\{P^{high},F\}^{z_j\bar{z}_j}]z_j\bar{z}_j\\
&+\sum_{j\in\mathbb{Z}_*}(P^{z_j\bar{z}_j}+\{P^{high},F\}^{z_j\bar{z}_j}+\left\langle \partial_{x}\Omega_j/j,F^{y}\right\rangle%\\
%&-P^{z_j\bar{z}_{-j}}-\{P^{high},F\}^{z_j\bar{z}_{-j}}
)z_j\bar{z}_j.
\end{align*}

\subsection{The solution of Homological equation \eqref{homological}}
%\begin{The}
%Consider an nearly integrable Hamiltonian
%$$H(x,y,z,\bar{z};\xi)=N(x,y,z,\bar{z};\xi)+P(x,y,z,\bar{z};\xi)$$
%with
%$$N(x,y,z,\bar{z};\xi)=\sum_{1\le j\le n}\omega_j(\xi)y_j+\sum_{j\in\mathbb{Z}_* }\frac{\Omega_j(x;\xi)}{j}z_j\bar{z}_j.$$
%and $$P(x,y,z,\bar{z};\xi)=P^{low}(x,y,z,\bar{z};\xi)+P^{high}(x,y,z,\bar{z};\xi).$$
%Suppose the following estimates holds for the perturbation
%$$|||X_{P^{low}}|||^T_{p,p-1,D(s,r,r)\times\Pi}\le \epsilon,$$
%and
%$$|||X_{P^{high}}|||^T_{p,p-1,D(s,r,r)\times\Pi}\le 1.$$
%Write
%$$\Omega(x;\xi)=[\Omega](\xi)+\tilde{\Omega}(x;\xi)$$
%where $[\Omega]$ is the mean value of the variable $x$ over $\mathbb{T}^n$.
%Suppose
%$$\tilde{\Omega}^{\nu}_j(x;\xi)=\sum_{i=0}^{\nu} \tilde{\Omega}^i_j(x;\xi),\quad j\in\mathbb{Z}_*$$
%with
%\begin{equation}\label{omegai}
%\Vert \tilde{\Omega}^i(x;\xi)\Vert_{-1,D(s_i)\times \Pi_i}:=\sup_{j\in\mathbb{Z}_*}\Vert j^{-1}\tilde{\Omega}^i_j(x;\xi)\Vert_{D(s_i)\times \Pi_i}\le \epsilon_i^{1-}.
%\end{equation}
%\end{The}
 \begin{defn}
We say an analytic function $f(x)=\sum\limits_{k \in \mathbb{Z}^{n}}\widehat{f}(x)e^{ik\cdot x}, x\in\mathbb{T}^n$ belongs to $\mathcal{F}(b)$ where $b$ is an integer, if $k\in\mathbb{Z}^n$ with $\widehat{f}(k) \neq 0$ satisfying
 $$\sum_{1\le i\le n}k_ij_i=b, \ \ j_{i} \in J_n.$$
\end{defn}

The following lemma is useful for us to obtain the estimate of the solutions of homological equation.
\begin{Lem}\label{lemma1}
Consider the homological equation
\begin{equation*}
\left( {{\bf i}\omega  \cdot {\partial _\varphi } + \lambda \left( {1 + a\left( \varphi  \right)} \right)} \right)x\left( \varphi  \right) = R\left( \varphi  \right),\quad\varphi  \in D(s_{\nu}).
\end{equation*}
 Here, $(\omega,\lambda)$ satisfies the Diophantine conditions
 \begin{align*}
&|\langle k,\omega\rangle|\geq \eta_0/|k|^\tau,\quad \forall k\in\mathbb{Z}^n\setminus\{0\},\\
&|\langle k,\omega\rangle+\lambda|\geq \eta/|k|^\tau,\quad \forall k\in\mathbb{Z}^n,
\end{align*}
where $\eta_{0}$ and $\eta$ with $0<\eta\le\eta_0\ll 1$ are constants. The function
 $a\left( \theta  \right): D(s_{\nu}) \to \mathbb C$ obeys
  \begin{equation*}
  a\left( \theta  \right) = \sum\limits_{j = 0}^\nu {{a_j}\left( \theta  \right)},
  \end{equation*}
where ${a_j}\left( \theta  \right): D(s_j) \to \mathbb C$ is real analytic, $a_j\in \mathcal{F}(0)$ and
 \begin{equation*}
\begin{array}{l}
{\left\| {{a_j}} \right\|_{D({s_j})}} \le \epsilon _j^{1 - },
\end{array}
\end{equation*}
where $\epsilon _j^{1 - }$ represents $\epsilon _j^{1 -\rho }$ for some constant $\rho$ with $0<\rho\ll 1$. Moreover, the function $R(\varphi)$ is analytic in $D(s_{\nu})$ with $R\in\mathcal{F}(b)$. Then the homological equation has a unique solution $x=x(\varphi): D(s_{\nu}) \to \mathbb C$ which is real analytic, $x(\varphi)\in\mathcal{F}(b)$ and
\begin{equation*}
{\left\| x\right\|_{D(s_{\nu}-\sigma_\nu)}} \le \frac{c(n,\tau)}{\eta \sigma_\nu^{20\left( {n + \tau } \right)}}{\left\| R \right\|_{D({s_\nu})}},
\end{equation*}
where $c(n,\tau)$ is a constant depending on $n$ and $\tau$.
\end{Lem}

The existence, uniqueness and related estimates of solutions of the homological equation is listed in section \ref{slarge}. The proof of $x(\varphi)\in\mathcal{F}(b)$ can be obtained by the uniqueness of solutions. And the proof is the same as that of Lemma 4.5 in \cite{LiuYuan2014}, here we omit it.

From this lemma, we know that the equations \eqref{2.9}-\eqref{zbarz} have solutions obeying the following estimates
\begin{align}
&|||{F^x}|||^*_{D(s-\sigma,r,r)\times\Pi}\lessdot \frac{\epsilon}{\eta^2\sigma^{C}},\quad |||\partial_{\xi_a}{F^x}|||^{0*}_{D(s-\sigma,r,r)\times\Pi}\lessdot \frac{\epsilon}{|a|\eta^2\sigma^{C}},\label{1}\\
&|||{F^1}|||^*_{D(s-4\sigma,r-4\sigma,r-4\sigma)\times\Pi}\lessdot \frac{\epsilon K^n}{\eta^4\sigma^{C}},\quad |||\partial_{\xi_a}{F^1}|||^{0*}_{D(s-4\sigma,r-4\sigma,r-4\sigma)\times\Pi}\lessdot \frac{\epsilon K^n}{|a|\eta^4\sigma^{C}},\label{3}\\
&|||{F^y}|||^*_{D(s-7\sigma,r-7\sigma,r-7\sigma)\times\Pi}\lessdot \frac{\epsilon K^{2n}}{\eta^6\sigma^{C}},\quad |||\partial_{\xi_a}{F^y}|||^{0*}_{D(s-7\sigma,r-7\sigma,r-7\sigma)\times\Pi}\lessdot \frac{\epsilon K^{2n}}{|a|\eta^6\sigma^{C}},\label{2}\\
&|||{F^2}|||^*_{D(s-7\sigma,r-7\sigma,r-7\sigma)\times\Pi}\lessdot \frac{\epsilon K^{2n}}{\eta^6\sigma^{C}},\quad |||\partial_{\xi_a}{F^2}|||^{0*}_{D(s-7\sigma,r-7\sigma,r-7\sigma)\times\Pi}\lessdot \frac{\epsilon K^{2n}}{|a|\eta^6\sigma^{C}},\label{4}
\end{align}
where $a\in\mathbb{Z}_*$ and $C$ is a constant depending on $n$ and $\tau$. In addition, we have
$$F^x\in \mathcal{F}(0),\quad F^y_j\in \mathcal{F}(0),\quad 1\le j\le n,$$
$$F_j^z\in \mathcal{F}(j),\quad F_j^{\bar{z}}\in \mathcal{F}(-j),\quad j\in \mathbb{Z}_*,$$
$$F_{ij}^{zz}\in \mathcal{F}(i+j),\quad F_{ij}^{\bar{z}\bar{z}}\in \mathcal{F}(-i-j),\quad F_{ij}^{z\bar{z}}\in \mathcal{F}(i-j),\quad  i,j\in \mathbb{Z}_*,$$
which implies the momentum conservation for the solutions of homological equations.

\medskip

These estimates can be obtained by the following order $\eqref{1}\rightarrow \eqref{3}\rightarrow\eqref{2}\rightarrow\eqref{4}$.
Here we only consider estimate \eqref{4}, which is more difficult than the others. For this purpose, let's look at \eqref{zbarz}. Let $F^{z_i\bar{z}_j}=\Gamma_K G^{z_i\bar{z}_j}$, then equation \eqref{zbarz} becomes
\begin{align*}
&\partial_\omega  G^{z_i\bar{z}_j}+{\bf i}(\Omega_i-\Omega_j) G^{z_i\bar{z}_j}=R^{z_i\bar{z}_j}+\{P^{high},F\}^{z_i\bar{z}_j},\quad i,j\in\mathbb{Z}_*,\,\,i\neq  j.
\end{align*}
Applying $\partial_{\xi_{a}}$ to both sides of this equation, we arrive at
\begin{align*}
&\partial_\omega (\partial_{\xi_a} G^{z_i\bar{z}_j})+{\bf i}(\Omega_i-\Omega_j) \partial_{\xi_a} G^{z_i\bar{z}_j}\\
=&\partial_{\xi_a} R^{z_i\bar{z}_j}+\partial_{\xi_a} \{P^{high},F\}^{z_i\bar{z}_j}-\partial_{\xi_a}\omega\cdot \partial_x G^{z_i\bar{z}_j}-{\bf i} \partial_{\xi_a}(\Omega_i-\Omega_j) G^{z_i\bar{z}_j},\quad i,j,a\in\mathbb{Z}_*,\,\,i\neq j.
\end{align*}
From \eqref{huph} \eqref{high2}, \eqref{1}, \eqref{3} and Lemma \ref{brak}, one has the estimate:
\begin{align}\label{Fengg}
&|||{\{P^{high},F\}^{2}}|||_{D(s-5\sigma,r-5\sigma,r-5\sigma)\times\Pi}\nonumber\\
=&|||{\{P^{high},F^x+F^1\}}|||_{D(s-5\sigma,r-5\sigma,r-5\sigma)\times\Pi}\lessdot  \frac{\epsilon K^{2n}}{\eta^4(r-4\sigma)^2\sigma^{C}}\lessdot \frac{\epsilon K^{2n}}{\eta^4\sigma^{C}},
\end{align}
and
\begin{align}\label{Fengg1}
&|||\partial_{\xi_a}{\{P^{high},F\}^{2}}|||^0_{D(s-5\sigma,r-5\sigma,r-5\sigma)\times\Pi}\nonumber\\
=&|||\partial_{\xi_a}{\{P^{high},F^x+F^1\}}|||^0_{D(s-5\sigma,r-5\sigma,r-5\sigma)\times\Pi}\lessdot \frac{\epsilon K^{2n}}{|a|\eta^4\sigma^{C}},
\end{align}
where we use $r-4\sigma>\sigma$ and put $(r-4\sigma)^2$ into $\sigma^{C}$.
By applying Lemma \ref{lemma1}, one immediately gets
$$|G^{z_i\bar{z}_j}|_{D(s-6\sigma)\times \Pi}\lessdot \frac{1}{\max\{|i|,|j|\}\eta \sigma^{C}}\Vert R^{z_i\bar{z}_j}+\{P^{high},F\}^{z_i\bar{z}_j}\Vert_{D(s-5\sigma)\times \Pi},$$
and
$$|\partial_{\xi_a} G^{z_i\bar{z}_j}|_{D(s-7\sigma)\times \Pi}\lessdot \frac{1}{\max\{|i|,|j|\}\eta \sigma^{C}}\Bigg\{\Vert \partial_{\xi_a}R^{z_i\bar{z}_j}+\partial_{\xi_a}\{P^{high},F\}^{z_i\bar{z}_j}\Vert_{D(s-6\sigma)\times \Pi}$$
$$+| \partial_{\xi_a}\omega|_{\Pi}\cdot \Vert\partial_x G^{z_i\bar{z}_j}\Vert_{D(s-6\sigma)\times \Pi}+(|i|+|j|)| \partial_{\xi_a}\Omega|_{-1,D(s-6\sigma)\times \Pi}\cdot \Vert G^{z_i\bar{z}_j}\Vert_{D(s-6\sigma)\times \Pi}\Bigg\}.$$
By Cauchy's estimate and \eqref{hupl}, \eqref{Fengg}, \eqref{Fengg1}, it immediately implies,
\begin{align*}
\Vert G^{z_i\bar{z}_j}\Vert_{D(s-7\sigma)\times \Pi}\lessdot &\frac{1}{\max\{|i|,|j|\}\eta^2 \sigma^{C}}(1+|\partial_\xi \omega|_{\Pi}+|\partial_\xi\Omega|_{-1,D(s-6\sigma)\times \Pi})\\
&\cdot\Vert R^{z_i\bar{z}_j}+\{P^{high},F\}^{z_i\bar{z}_j}\Vert_{D(s-6\sigma)\times \Pi},\\
 \lessdot& \frac{\epsilon K^{2n}}{\max\{|i|,|j|\}\eta^6 \sigma^{C}},
 \end{align*}
and
$$\Vert\partial_{\xi_a} G^{z_i\bar{z}_j}\Vert_{D(s-7\sigma)\times \Pi}\lessdot \frac{\epsilon K^{2n}}{|a|\max\{|i|,|j|\}\eta^6 \sigma^{C}}. $$
Since $F^{z_i\bar{z}_j}=\Gamma_K G^{z_i\bar{z}_j}$, we also have the following estimate for $F$:
$$\Vert F^{z_i\bar{z}_j}\Vert_{D(s-7\sigma)\times \Pi}
 \lessdot \frac{\epsilon K^{2n}}{\max\{|i|,|j|\}\eta^6  \sigma^{C}},$$
$$\Vert\partial_{\xi_a} F^{z_i\bar{z}_j}\Vert_{D(s-7\sigma)\times \Pi}\lessdot \frac{\epsilon K^{2n}}{|a|\max\{|i|,|j|\}\eta^6  \sigma^{C}}. $$

Similar estimates holds for $F^{z_iz_j}$ and $F^{\bar{z}_i\bar{z}_j}$. Therefore, we get
\begin{align*}
&|||{F^2}|||^*_{D(s-7\sigma,r-7\sigma,r-7\sigma)\times\Pi}
\lessdot \frac{\epsilon K^{2n}}{\eta^6\sigma^{C}},\quad |||\partial_{\xi_a}{F^2}|||^*_{D(s-7\sigma,r-7\sigma,r-7\sigma)\times\Pi}
\lessdot \frac{\epsilon K^{2n}}{|a|\eta^6\sigma^{C}}.
\end{align*}
%Then from Cauchy's estimate in Lemma \ref{}, one has
%\begin{equation}\label{df}
%|||DX_{F}|||_{D(s-6\sigma,r-6\sigma,r-6\sigma)\times\Pi}\lessdot \frac{\epsilon}{\gamma^6\sigma^{C(n,\tau)}}.
%\end{equation}
Note $\Phi=X_F^1$. Thus, by Lemma \ref{youy} and Lemma $A.4$ in \cite{P1},  we have
$$||\Phi-id||_{p,D(s-8\sigma,r-8\sigma,r-8\sigma)\times\Pi}\le \Vert X_F\Vert_{p,D(s-8\sigma,r-8\sigma,r-8\sigma)\times\Pi}$$
$$\lessdot \frac{1}{(r-8\sigma)^2\sigma^p} |||F|||^*_{D(s-7\sigma,r-7\sigma,r-7\sigma)\times\Pi}\lessdot \frac{\epsilon K^{2n}}{\eta^6\sigma^{C}},$$
where we suppose $\sigma<r-8\sigma$. It implies that $\Phi$ is well defined and maps the domain $D(s-10\sigma,r-10\sigma,r-10\sigma)$ into $D(s-9\sigma,r-9\sigma,r-9\sigma)$. Moreover, by Cauchy's inequality in Lemma A.3 in \cite{P1}, we obtain
$$||D\Phi-Id||_{D(s-10\sigma,r-10\sigma,r-10\sigma)\times\Pi}\lessdot \frac{\epsilon K^{2n}}{\eta^6\sigma^{C}}.$$
\subsection{The new normal form}
From \eqref{o+}, %and \eqref{O+}
we have
$N_+=N+\hat{N}$
with
$$\hat{N}=\sum_{1\le j\le n}\hat{\omega}_{j}(\xi)y_j+\sum_{j\in\mathbb{Z}_*}\hat{\Omega}_{j}(x;\xi)z_j\bar{z}_j,$$
where
\begin{align}
&\hat{\omega}_j(\xi)=[P^{y_j}+\{P^{high},F\}^{y_j}], \ \ 1 \leq j \leq n, \label{om}\\
&\hat{\Omega}_{j}(x;\xi)=j\big(P^{z_j\bar{z}_j}+\{P^{high},F\}^{z_j\bar{z}_j}+\langle \partial_{x}\Omega_j/j,F^{y}\rangle
\big),\quad j\in\mathbb{Z}_*.\label{Om}
\end{align}
From \eqref{om}, one easily finds that $\hat{\omega} = (\hat{w}_{j}: 1 \leq j \leq n)$ obeys
\begin{align*}
&\Vert \hat{\omega}\Vert_{D(s-7\sigma)\times \Pi}\lessdot  \Vert P^{y_j}\Vert_{D(s)\times\Pi}+\Vert\{P^{high},F\}^{y_j}\Vert_{D(s-7\sigma)\times\Pi}
\lessdot \frac{\epsilon K^{n}}{\eta^6\sigma^{C}},\\
&\Vert \partial_{\xi_a}\hat{\omega}\Vert_{D(s-7\sigma)\times \Pi}\lessdot \Vert\partial_{\xi_a}P^{y_j}\Vert_{D(s)\times\Pi}+\Vert\partial_{\xi_a}\{P^{high},F\}^{y_j}\Vert_{D(s-7\sigma)\times\Pi}
\lessdot \frac{\epsilon K^{n}}{|a|\eta^6\sigma^{C}}.
\end{align*}
In the following, we estimate $\hat{\Omega} = (\hat{\Omega}_{j}:j \in \mathbb{Z}_{*})$. From \eqref{2} and \eqref{omegai}, one has that for $a\in\mathbb{Z}_*$
$$
||\langle \partial_{x}\Omega_j/j,F^{y_j}\rangle||_{D(s-7\sigma)\times\Pi}\le \sigma^{-1}\Vert\tilde{\Omega}_j/j\Vert_{D(s-7\sigma)\times\Pi}\cdot ||{F^{y_j}}||_{D(s-7\sigma)\times\Pi}\le  \frac{\epsilon K^{2n}}{\eta^6\sigma^{C}},
$$
$$||\partial_{\xi_a}\langle \partial_{x}\Omega_j/j,F^{y_j}\rangle||_{D(s-7\sigma)\times\Pi}\le \sigma^{-1}\Vert\partial_{\xi_a}\tilde{\Omega}_j/j\Vert_{D(s-7\sigma)\times\Pi}\cdot ||{F^{y_j}}||_{D(s-7\sigma)\times\Pi}$$
$$+ \sigma^{-1}\Vert\tilde{\Omega}_j/j\Vert_{D(s-7\sigma)\times\Pi}\cdot ||\partial_{\xi_a}{F^{y_j}}||_{D(s-7\sigma)\times\Pi}\le  \frac{\epsilon K^{2n}}{|a|\eta^6\sigma^{C}}.$$
Together with
$$\Vert P^{z_j\bar{z}_j}\Vert_{D(s)\times \Pi}\le \frac{\epsilon}{r^2}\le \frac{\epsilon}{\sigma^2},\quad \Vert \partial_{\xi_a}P^{z_j\bar{z}_{j}}\Vert_{D(s)\times \Pi}
\le \frac{\epsilon}{|a|\sigma^2},$$

$$\Vert\{P^{high},F\}^{z_j\bar{z}_j}\Vert_{D(s-7\sigma)\times\Pi}
\le\frac{\epsilon K^{n}}{\eta^6\sigma^{C}},$$
and
$$\Vert\partial_{\xi_a}\{P^{high},F\}^{z_j\bar{z}_j}\Vert_{D(s-7\sigma)\times\Pi}
\le\frac{\epsilon K^{n}}{|a|\eta^6\sigma^{C}},$$
we obtain
$$||\hat{\Omega}||_{D(s-7\sigma)\times\Pi}\lessdot \frac{\epsilon K^{2n}}{\eta^6\sigma^{C}},\quad ||\partial_{\xi_a}\hat{\Omega}||_{D(s-7\sigma)\times\Pi}\lessdot \frac{\epsilon K^{2n}}{|a|\eta^6\sigma^{C}}.$$
%Note
%\begin{align*}
%\partial_\xi \hat{\Omega}=\partial_\xi  P^{z_j\bar{z}_j}+\partial_\xi \{P^{high},F\}^{z_j\bar{z}_j}+\langle \partial_{x}\partial_\xi\Omega_j,F^{y}\rangle+\langle \partial_{x}\Omega_j,\partial_\xi F^{y}\rangle
%\end{align*}
%Similar to the above discussion, we have
%$$||\hat{\Omega}||_{D(s-\sigma)\times \Pi}\lessdot j^\delta \sigma^{-1} K^C\epsilon.$$

\subsection{Estimate the new high order term and the perturbation}
Firstly, by Lemma \ref{brak}, we have
\begin{align*}
&|||\{P^{high},F\}^{high}|||_{D(s-8\sigma,r-8\sigma,r-8\sigma)\times\Pi},\quad |||\{P^{high},F\}^{low}|||_{D(s-8\sigma,r-8\sigma,r-8\sigma)\times\Pi}\\
\le &|||\{P^{high},F\}|||_{D(s-8\sigma,r-8\sigma,r-8\sigma)\times\Pi}\\
\lessdot &\max\left\{\frac{K^n}{(r-7\sigma)\sigma},\frac{1}{(r-7\sigma)^2\sigma}\right\}|||P^{high}|||_{D(s-7\sigma,r-7\sigma,r-7\sigma)\times\Pi} \\
&\cdot|||F|||^*_{D(s-7\sigma,r-7\sigma,r-7\sigma)\times\Pi}
\lessdot  \frac{\epsilon K^{3n}}{\eta^6\sigma^{C}},
\end{align*}
and
\begin{align*}
&|||\partial_{\xi_a}\{P^{high},F\}^{high}|||^0_{D(s-8\sigma,r-8\sigma,r-8\sigma)\times\Pi},\\
&\quad |||\partial_{\xi_a}\{P^{high},F\}^{low}|||^0_{D(s-8\sigma,r-8\sigma,r-8\sigma)\times\Pi}
\lessdot   \frac{\epsilon K^{3n}}{|a|\eta^6\sigma^{C}}.
\end{align*}
Secondly, we estimate the term $\int_0^1 \{P^{low},F\}\circ X_F^tdt$. Actually, from Lemma \ref{lemma5.3}, we derive
\begin{align*}
&|||{\int_0^1 \{P^{low},F\}\circ X_F^tdt}|||_{D(s-9\sigma,r-9\sigma,r-9\sigma)\times\Pi}\\
\le&|||\{P^{low},F\}\circ X_F^t|||_{D(s-9\sigma,r-9\sigma,r-9\sigma)\times\Pi}\\
\lessdot & |||\{P^{low},F\}|||_{D(s-8\sigma,r-8\sigma,r-8\sigma)\times\Pi}\\
 \le&\frac{\epsilon^2 K^{3n}}{\eta^6\sigma^{C}},
\end{align*}
and
\begin{align*}
|||\partial_{\xi_a}{\int_0^1 \{P^{low},F\}\circ X_F^tdt}|||^0_{D(s-9\sigma,r-9\sigma,r-9\sigma)\times\Pi}
 \le& \frac{\epsilon^2 K^{3n}}{|a|\eta^6\sigma^{C}}, \ a \in \mathbb{Z}_{*}.
\end{align*}
Then, we consider the term $\int_0^1(1-t)\{\{P^{high},F\},F\}\circ X_F^tdt$
\begin{align*}
&|||{\int_0^1(1-t)\{\{P^{high},F\},F\}\circ X_F^tdt}|||_{D(s-10\sigma,r-10\sigma,r-10\sigma)\times\Pi}\\
\le & |||\{\{P^{high},F\},F\}\circ X_F^t|||_{D(s-10\sigma,r-10\sigma,r-10\sigma)\times\Pi}\\
\lessdot & |||\{\{P^{high},F\},F\}|||_{D(s-9\sigma,r-9\sigma,r-9\sigma)\times\Pi}\\
 \le & \left(\frac{\epsilon K^{3n}}{\eta^6\sigma^{C}}\right)^2,
\end{align*}
and for $a\in\mathbb{Z}_*$,
$$|||\partial_{\xi_a}{\int_0^1(1-t)\{\{P^{high},F\},F\}\circ X_F^tdt}|||^0_{D(s-10\sigma,r-10\sigma,r-10\sigma)\times\Pi}\le \frac{1}{|a|}\left(\frac{\epsilon K^{3n}}{\eta^6\sigma^{C}}\right)^2.$$

Next, we are going to estimate $\hat{P}=\{N,(1-\Gamma_K)G\}$. Since
$$\{N,G\}=-P^{low}-\{P^{high},F\}^{low}+\hat{N},$$
we have
$$\{[N],G\}=-P^{low}-\{P^{high},F\}^{low}+\hat{N}-\{N-[N],G\},$$
where $[N]$ denotes the mean value of the function $N$ with respect to the variable $x$ over $\mathbb{T}^n$. Therefore,
\begin{equation}\label{ngng}
\{[N], (1-\Gamma_K)G\}=(1-\Gamma_K)(-P^{low}-\{P^{high},F\}^{low}+\hat{N}-\{N-[N],G\}).
\end{equation}
Let $K=\frac{1}{\sigma}|\log \frac{1}{\epsilon}|$. From \eqref{ngng}, it is easy to obtain
\begin{equation}\label{hatp1}
|||\{[N], (1-\Gamma_K)G\}|||_{D(s-8\sigma,r-8\sigma,r-8\sigma)\times\Pi}\le \frac{\epsilon^2 K^{3n}}{\eta^6\sigma^{C}}.
\end{equation}
In addition, we have
\begin{equation}\label{hatp2}
|||\{N-[N], (1-\Gamma_K)G\}|||_{D(s-8\sigma,r-8\sigma,r-8\sigma)\times\Pi}\le \frac{\epsilon^2 K^{3n}}{\eta^6\sigma^{C }}.
\end{equation}
Combining \eqref{hatp1} and \eqref{hatp2}, one obtains
$$|||\hat{P}|||_{D(s-8\sigma,r-8\sigma,r-8\sigma)\times\Pi}=|||\{N, (1-\Gamma_K)G\}|||_{D(s-8\sigma,r-8\sigma,r-8\sigma)\times\Pi}\le \frac{\epsilon^2 K^{3n}}{\eta^6\sigma^{C}}.
$$
Furthermore, we can get that for $a\in\mathbb{Z}_*$,
$$|||\partial_{\xi_a}\hat{P}|||^0_{D(s-8\sigma,r-8\sigma,r-8\sigma)\times\Pi}\le \frac{\epsilon^2 K^{3n}}{|a|\eta^6\sigma^{C}}.
$$

Finally, we study the term $\int_0^1(1-t)\{\{N,F\},F\}\circ X_F^tdt$. Note
$$\{N,F\}=-P^{low}-\{P^{high},F\}^{low}+\hat{N}+\hat{P}.$$
Thus, we have that, for $a\in\mathbb{Z}_*$,
$$|||{\{N,F\}}|||_{D(s-8\sigma,r-8\sigma,r-8\sigma)\times\Pi} \le \frac{\epsilon K^{3n}}{\eta^6\sigma^{C}},$$
$$|||\partial_{\xi_t}{\{N,F\}}|||^0_{D(s-8\sigma,r-8\sigma,r-8\sigma)\times\Pi} \le \frac{\epsilon K^{3n}}{|a|\eta^6\sigma^{C}}.$$
It immediately follows that
\begin{align*}
&|||{\int_0^1(1-t)\{\{N,F\},F\}\circ X_F^tdt}|||_{D(s-10\sigma,r-10\sigma,r-10\sigma)\times\Pi}\\
=& |||\{\{N,F\},F\}\circ X_F^t|||_{D(s-10\sigma,r-10\sigma,r-10\sigma)\times\Pi}\\
\le & |||\{\{N,F\},F\}|||_{D(s-9\sigma,r-9\sigma,r-9\sigma)\times\Pi}\\
 \le& \left(\frac{\epsilon K^{3n}}{\eta^6\sigma^{C}}\right)^2,
\end{align*}
and, for $a\in\mathbb{Z}_*$,
\begin{align*}
&|||\partial_{\xi_a}{\int_0^1(1-t)\{\{N,F\},F\}\circ X_F^tdt}|||^0_{D(s-10\sigma,r-10\sigma,r-10\sigma)\times\Pi} \le \frac{1}{|a|} \left(\frac{\epsilon K^{3n}}{\eta^6\sigma^{C}}\right)^2.
\end{align*}
From \eqref{phigh}, we let
\begin{align*}
P_+^{low}=&(\int_0^1(1-t)\{\{N,F\},F\}\circ X_F^tdt)^{low}\\
&+\hat{P}+(\int_0^1 \{P^{low},F\}\circ X_F^tdt)^{low}\\
&+(\int_0^1(1-t)\{\{P^{high},F\},F\}\circ X_F^tdt)^{low}.
\end{align*}
and
\begin{align*}
P_+^{high}=&(\int_0^1(1-t)\{\{N,F\},F\}\circ X_F^tdt)^{high}\\
&+ P^{high}+\{P^{high},F\}^{high}+(\int_0^1 \{P^{low},F\}\circ X_F^tdt)^{high}\\
&+(\int_0^1(1-t)\{\{P^{high},F\},F\}\circ X_F^tdt)^{high}.
\end{align*}
Combining the above estimates together, one gets that
\begin{align*}
&|||{P_+^{low}}|||_{D(s-10\sigma,r-10\sigma,r-10\sigma)\times\Pi}\le \left(\frac{\epsilon K^{3n}}{\eta^6\sigma^{C}}\right)^2,\\
& |||\partial_{\xi_a}{P_+^{low}}|||^0_{D(s-10\sigma,r-10\sigma,r-10\sigma)\times\Pi}\le \frac{1}{|a|}\left(\frac{\epsilon K^{3n}}{\eta^6\sigma^{C}}\right)^2,
\end{align*}
and that, for $a\in\mathbb{Z}_*$,
\begin{align*}
&|||{P_+^{high}}|||_{D(s-10\sigma,r-10\sigma,r-10\sigma)\times\Pi}\le 1+\frac{\epsilon K^{3n}}{\eta^6\sigma^{C}}+\left(\frac{\epsilon K^{3n}}{\eta^6\sigma^{C}}\right)^2,\\
&|||\partial_{\xi_a}{P_+^{high}}|||^0_{D(s-10\sigma,r-10\sigma,r-10\sigma)\times\Pi} \le \frac{1}{|a|}\left(1+\frac{\epsilon K^{3n}}{\eta^6\sigma^{C}}+\left(\frac{\epsilon K^{3n}}{\eta^6\sigma^{C}}\right)^2\right).
\end{align*}

\subsection{Iterative and convergence}
Firstly, we define some iterative constants:
\begin{itemize}
\item $\eta_m=\eta 2^{-m}, m=0,1,2,\cdots$;
\item $\epsilon_0=\epsilon=\eta^{12}\varepsilon$, $\epsilon_m=\eta^{12}\varepsilon^{(4/3)^m}, m=0,1,2,\cdots$;
\item $\tau_0=0, \tau_m=(1^{-2}+\cdots+m^{-2})/2\sum_{j=1}^\infty j^{-2}, m=1,2,\cdots$;
\item Given $0<s_0, r_0\le 1$. Let $0<\chi\le \min\{s_0,r_0\}$ and $s_m=(1-\tau_m)\chi, r_m=(1-\tau_m)\chi$, $m=1,2\cdots$;
%\item $\sigma_m=\frac{s_m-s_{m+1}}{10}=\frac{r_m-r_{m+1}}{10}$,  $m=0,1,2\cdots$.
\item $K_m=|\frac{1}{r_{m}-r_{m+1}}\log\frac{1}{\epsilon_m}|$, $m=0,1,2\cdots$.
\end{itemize}
\begin{Lem}[Iterative Lemma]
Consider the nearly integrable Hamiltonian
$$H_m(x,y,z,\bar{z};\xi)=N_m(x,y,z,\bar{z};\xi)+P_m(x,y,z,\bar{z};\xi),$$
where
$$N_m(x,y,z,\bar{z};\xi)=\sum_{1\le j\le n}\omega_{mj}(\xi)y_j+\sum_{j\in\mathbb{Z}_* }\frac{\Omega_{mj}(x;\xi)}{j}z_j\bar{z}_j,$$
%$$=\sum_{1\le j\le n}\omega_{mj}(\xi)y_j+\sum_{|j|\le \mathcal{N} }\frac{\Omega_{mj}(\xi)}{j}z_j\bar{z}_j+\sum_{|j|\geq \mathcal{N} +1}\frac{\Omega_{mj}(x;\xi)}{j}z_j\bar{z}_j$$
$P_m(x,y,z,\bar{z};\xi)=P_m^{low}(x,y,z,\bar{z};\xi)+P_m^{high}(x,y,z,\bar{z};\xi)$ is the perturbation with
$$P_m^{low}(x,y,z,\bar{z};\xi)=\sum_{\substack{\alpha\in\mathbb{N}^n,\beta,\gamma\in \mathbb{N}^{\mathbb{Z}_*},\\ 2|\alpha|+|\beta+\gamma|\le 2}}P_m^{\alpha\beta\gamma}(x;\xi)y^{\alpha}z^{\beta}\bar{z}^{\gamma},$$
and
$$P_m^{high}(x,y,z,\bar{z};\xi)=\sum_{\substack{\alpha\in\mathbb{N}^n,\beta,\gamma\in \mathbb{N}^{\mathbb{Z}_*},\\ 2|\alpha|+|\beta+\gamma|\geq 3}}P_m^{\alpha\beta\gamma}(x;\xi)y^{\alpha}z^{\beta}\bar{z}^{\gamma}.$$
Suppose the assumptions \textbf{A} and \textbf{B} in Section $2.2$ are fulfilled for $\omega_0$ and $\Omega_0$. In addition, $\omega_m$ obeys
$$\Vert \omega_m(\xi)-\omega_0(\xi)\Vert_{\Pi_m}\le \sum_{i=1}^m\epsilon_{i-1}^{2/3},\quad \Vert \partial_{\xi_a}\omega_m(\xi)\Vert_{\Pi_m}\le \sum_{i=1}^m\frac{\epsilon_{i-1}^{2/3}}{|a|},$$
and $\Omega_{m}$ satifies,
$$\Vert \Omega_m(x;\xi)-\Omega_0(\xi)\Vert_{-1,D(s_m)\times\Pi_m}\le \sum_{i=1}^m\epsilon_{i-1}^{2/3},$$
$$\Vert \partial_{\xi_a}(\Omega_m(x;\xi)-\Omega_0(\xi))\Vert_{-1,D(s_m)\times\Pi_m}\le\sum_{i=1}^m\frac{\epsilon_{i-1}^{2/3}}{|a|},$$
where $a\in\mathbb{Z}_*$.
Suppose the assumption \textbf C in Section $2.2$ is fulfilled for the initial perturbation $P_0$. Furthermore, the perturbations $P_m$ is real analytic and $P_{m}^{low}$ and $P_{m}^{high}$ obey the following estimates:
$$|||{P_m^{low}}|||_{D(s_m,r_m,r_m)\times\Pi_m}\le \epsilon_m,\quad |||\partial_{\xi_a}{P_m^{low}}|||^0_{D(s_m,r_m,r_m)\times\Pi_m}\le \frac{\epsilon_m}{|a|},$$
and
$$|||{P_m^{high}}|||_{D(s_m,r_m,r_m)\times\Pi_m}\le \epsilon_0+\sum_{i=1}^m\epsilon_i^{2/3},$$
$$ |||\partial_{\xi_a}{P_m^{high}}|||^0_{D(s_m,r_m,r_m)\times\Pi_m}\le \frac{\epsilon_0}{|a|}+\sum_{i=1}^m\frac{\epsilon_i^{2/3}}{|a|},$$
where $a\in\mathbb{Z}_*$.

Write
$$\Omega_m(x;\xi)=[\Omega]_m(\xi)+\tilde{\Omega}_m(x;\xi),$$
where $[\Omega]$ is the mean value of the variable $x$ over $\mathbb{T}^n$. Let
$$\mathcal{R}_{kl}^m=\left\{\xi\in\Pi:|\langle k,\omega_m(\xi)+\langle l,[\Omega]_m(\xi)\rangle|\le \frac{\eta_m  M_{l}}{(|k|+1)^\tau}\right\},\quad k\in\mathbb{Z}^n, l\in \mathbb{Z}^{\mathbb{Z}_*},$$
and let
$$\Pi_{m+1}=\Pi_m\setminus\bigcup_{\substack{k\in\mathbb{Z}^n, l\in \mathbb{Z}^{\mathbb{Z}_*}\\|l|\le 2, |k|+|l|\neq 0}}\mathcal{R}_{kl}^m.$$
Then for each $\xi\in\Pi_{m+1}$, the homological equation
$$\{N_m,F_m\}+P_m^{low}+\{P_m^{high},F_m\}^{low}=\hat{N}_m+\hat{P}_m,$$
has a solution $F_m(x,y,z,\bar{z};\xi)$ with the estimate
$$|||{F_m}|||^*_{D(s_{m+1},r_{m+1},r_{m+1})\times\Pi_{m+1}}\le \epsilon_m^{2/3},$$
$$ |||\partial_{\xi_a}{F_m}|||^{0*}_{D(s_{m+1},r_{m+1},r_{m+1})\times\Pi_{m+1}}\le \frac{\epsilon_m^{2/3}}{|a|},$$
and
$$|||\hat{N}_m|||_{D(s_{m+1},r_{m+1},r_{m+1})\times\Pi_{m+1}}\le \epsilon_m^{2/3},$$
$$ |||\partial_{\xi_a}\hat{N}_m|||^0_{D(s_{m+1},r_{m+1},r_{m+1})\times\Pi_{m+1}}\le \frac{\epsilon_m^{2/3}}{|a|},$$
where $a\in\mathbb{Z}_*$.
Moreover, under the transformation $\Psi_m=X_{F_m}^t|_{t=1}$, the Hamiltonian reads
$$H_{m+1}=H_m\circ\Psi_m=N_{m+1}(x,y,z,\bar{z};\xi)+P_{m+1}(x,y,z,\bar{z};\xi),$$
where
$$N_{m+1}(x,y,z,\bar{z};\xi)=N_m(x,y,z,\bar{z};\xi)+\hat{N}_m(x,y,z,\bar{z};\xi)$$
$$=\sum_{1\le j\le n}\omega_{m+1j}(\xi)y_j+\sum_{j\in\mathbb{Z}_* }\frac{\Omega_{m+1j}(x;\xi)}{j}z_j\bar{z}_j,$$
and
\begin{align*}
P_{m+1}=&\hat{P}_m+\int_0^1(1-t)\{\{N_m,F_m\},F\}\circ X_{F_m}^tdt\\
&+ P^{high}_m+\{P_m^{high},F_m\}^{high}+\int_0^1 \{P_m^{low},F_m\}\circ X_{F_m}^tdt\\
&+\int_0^1(1-t)\{\{P_m^{high},F_m\},F_m\}\circ X_{F_m}^tdt.
\end{align*}
Moreover the following estimates hold:

\noindent(1) for each $\xi\in \Pi_{m+1}$, the symplectic map $\Psi_m=X_{F_m}^t|_{t=1}$ satisfies
$$\Vert \Psi_m-id\Vert_{p,D(s_{m+1},r_{m+1},r_{m+1})}\le \epsilon_m^{2/3},$$
and
$$\Vert D\Psi_m-Id\Vert_{p,D(s_{m+1},r_{m+1},r_{m+1})}\le \epsilon_m^{2/3};$$
(2) the frequencies $\omega_{m+1}(\xi)$ and $\Omega_{m+1}(\xi)$ satisfy, for $a\in\mathbb{Z}_*$
$$\Vert \omega_{m+1}(\xi)-\omega_0(\xi)\Vert_{\Pi_{m+1}}\le \sum_{i=1}^{m+1}\epsilon_{i-1}^{2/3},\quad \Vert \partial_{\xi_a}\omega_{m+1}(\xi)\Vert_{\Pi_{m+1}}\le \sum_{i=1}^{m+1}\frac{\epsilon_{i-1}^{2/3}}{|a|},$$
and
$$\Vert \Omega_{m+1}(x;\xi)-\Omega_0(\xi)\Vert_{-1,D(s_{m+1})\times\Pi_{m+1}}\le \sum_{i=1}^{m+1}\epsilon_{i-1}^{2/3},$$
$$ \Vert \partial_{\xi_a}(\Omega_{m+1}(x;\xi)-\Omega_0(\xi))\Vert_{-1,D(s_{m+1})\times\Pi_{m+1}}\le \sum_{i=1}^{m+1}\frac{\epsilon_{i-1}^{2/3}}{|a|};$$
(3) the perturbation $P_{m}$ is real analytic and satisfies
$$|||{P_{m+1}^{low}}|||_{D(s_{m+1},r_{m+1},r_{m+1})\times\Pi_{m+1}}\le \epsilon_{m+1},$$
$$ |||\partial_{\xi_a}{P_{m+1}^{low}}|||^0_{D(s_{m+1},r_{m+1},r_{m+1})\times\Pi_{m+1}}\le \frac{\epsilon_{m+1}}{|a|},$$
and
$$|||{P_{m+1}^{high}}|||_{D(s_{m+1},r_{m+1},r_{m+1})\times\Pi_{m+1}}\le \epsilon_0+\sum_{i=1}^{m+1}\epsilon_i^{2/3},$$
$$|||\partial_{\xi_a}{P_{m+1}^{high}}|||^0_{D(s_{m+1},r_{m+1},r_{m+1})\times\Pi_{m+1}}\le \frac{1}{|a|}\left(\epsilon_0+\sum_{i=1}^{m+1}\epsilon_i^{2/3}\right);$$
(4) the measure of the subset $\Pi_{m+1}$ of $\Pi_m$ satisfies
$${\rm Meas}\ \Pi_{m+1}\geq ({\rm Meas} \ \Pi_m)(1-O(\eta_m)).$$
\end{Lem}

Based on the previous iterative process, the iterative lemma follows immediately. Here we omit it. For the measure estimates, we can follow the proof of measure of $\Pi_{\acute{\eta}}$ in subsection \ref{Mease}.

\subsection{Proof of Theorem \ref{normaltheorem}}
Let $\Pi_\eta=\bigcap_{m=0}^\infty\Pi_m$, $D(s_0/2,r_0/2,r_0/2)\subset \bigcap_{m=0}^\infty D(s_m,r_m,r_m)$ and $\Psi=\Pi_{m=0}^\infty\Psi_m$. By standard argument in KAM theory, we conclude that $\Psi, D\Psi$, and $H_m$ converge uniformly on the domain $D(s_0/2,r_0/2,r_0/2)\times \Pi_\eta$. Let
$$\breve{H}(x,y,z,\bar{z};\xi):=\lim_{m\rightarrow\infty}H_m=\breve{N}(x,y,z,\bar{z};\xi)+\breve{P}(x,y,z,\bar{z};\xi),$$
where
$$\breve{N}(x,y,z,\bar{z};\xi)=\sum_{j=1}^n\breve{\omega}_j(\xi)y_j+\sum_{j\in\mathbb{Z}_*}\frac{\breve{\Omega}_j(x;\xi)}{j}z_j\bar{z}_j,$$
%$$=\sum_{j=1}^n\breve{\omega}_j(\xi)y_j+\sum_{|j|\le\mathcal{N}}\frac{\breve{\Omega}_j(\xi)}{j}z_j\bar{z}_j+\sum_{|j|\geq \mathcal{N}+1}\frac{\breve{\Omega}_j(x;\xi)}{j}z_j\bar{z}_j,$$
and
$$\breve{P}(x,y,z,\bar{z};\xi)=\sum_{\substack{\alpha\in\mathbb{N}^n,\beta,\gamma\in \mathbb{N}^{\mathbb{Z}_*}\\2|\alpha|+|\beta+\gamma|\geq 3}}\breve{P}^{\alpha\beta\gamma}(x,\xi)y^{\alpha}z^\beta\bar{z}^\gamma.$$
Moreover, by the standard KAM proof, we obtain the following estimates:

\noindent(1) for each $\xi\in \Pi_{\eta}$, the symplectic map $\Psi$ satisfies
$$\Vert \Psi-id\Vert_{p, D(s_{0}/2,r_{0}/2,r_{0}/2)}\le c\eta^{6}\varepsilon,$$
and
$$\Vert D\Psi-Id\Vert_{p, D(s_{0}/2,r_{0}/2,r_{0}/2)}\le c\eta^{6}\varepsilon;$$
(2) the frequencies $\breve{\omega}(\xi)$ and $\breve{\Omega}(\xi)$ satisfy that for $a\in \mathbb{Z}_*$
\begin{equation}\label{pgyy}
\Vert \breve{\omega}(\xi)-\omega_0(\xi)\Vert_{\Pi_\eta}\le c\eta^{6}\varepsilon, \quad \Vert \partial_{\xi_a}(\breve{\omega}(\xi)-\omega_0(\xi))\Vert_{\Pi_\eta}\le \frac{c\eta^{6}\varepsilon}{|a|},
\end{equation}
and
\begin{align}
&\Vert \breve{\Omega}(x;\xi)-\Omega_0(\xi)\Vert_{-1,D(s_{0}/2)\times\Pi_\eta}\le c\eta^{6}\varepsilon,\\
& \Vert \partial_{\xi_a}(\breve{\Omega}(x;\xi)-\Omega_0(\xi))\Vert_{-1,D(s_{0}/2)\times\Pi_\eta}\le \frac{c\eta^{6}\varepsilon}{|a|}.\label{iii}
\end{align}
(3) the perturbation $\breve{P}$ satisfies that, for $a\in \mathbb{Z}_*$,
\begin{equation}\label{brevep}
|||{\breve{P}}|||_{D(s_{0}/2,r_{0}/2,r_{0}/2)\times\Pi_\eta}\le \epsilon(1+c\eta^{6}\varepsilon),
\end{equation}
$$|||\partial_{\xi_a}{\breve{P}}|||^0_{D(s_{0}/2,r_{0}/2,r_{0}/2)\times\Pi_\eta}\le \frac{\epsilon(1+c\eta^{6}\varepsilon)}{|a|},$$
where $c$ is a constant depending on $s_0, r_0, n$ and $\tau$;

(4) the measure of the subset $\Pi_{\eta}$ satisfies
$${\rm Meas}\,\Pi_{\eta}\geq ({\rm Meas}\ \Pi)(1-O(\eta)).$$

\section{Long time stability theorem}
Given a large number $\mathcal{N}$, split the normal frequency $\breve{\Omega}(x;\xi)$ and normal variable $(z,\bar{z})$ into two parts, respectively, i.e.,
$$\breve{\Omega}(x;\xi)=(\acute{\Omega}(x,\xi),\hat{\Omega}(x;\xi)),\quad z=(\acute{z},\hat{z}),\quad\bar{z}=(\bar{\acute{z}},\bar{\hat{z}}),$$
where
$$\acute{\Omega}(\xi)=(\breve{\Omega}_j(x,\xi))_{|j|\le\mathcal{N}},\quad \acute{z}=(z_j)_{|j|\le\mathcal{N}},\quad \bar{\acute{z}}=(\bar{z}_j)_{|j|\le\mathcal{N}},$$
$$\hat{\Omega}(x;\xi)=(\breve{\Omega}_{j}(x;\xi),\cdots))_{|j|\geq\mathcal{N}+1},\quad\hat{z}=(z_j)_{|j|\geq\mathcal{N}+1},\quad \bar{\hat{z}}=(\bar{z}_j)_{|j|\geq\mathcal{N}+1}.$$

Moreover, we write
$$\breve{\Omega}(x;\xi)=[\breve{\Omega}](\xi)+\tilde{\breve{\Omega}}(x;\xi),$$
where $[\breve{\Omega}]$ is the mean value of the variable $x$ over $\mathbb{T}^n$.

Given $0<\acute{\eta}<1$ and $\tau>0$. If the frequencies $\breve{\omega}$ and $\breve{\Omega}(x;\xi)$ satisfy the following inequality
$$\big|\langle k,\breve{\omega}\rangle+\langle \acute{l},[\acute{\Omega}]\rangle+\langle \hat{l},[\hat{\Omega}]\rangle\big|\geq \frac{\acute{\eta}M_{\acute{l},\hat{l}}}{4^{\mathcal{M}}(|k|+1)^{\tau}C(\mathcal{N},\acute{l})},$$
for any $k\in\mathbb{Z}^n$, $\acute{l}\in \mathbb{Z}^{2\mathcal{N}}$ and $\hat{l}\in\mathbb{Z}^{\mathbb{Z}_*}$ with
$$|k|+|\acute{l}|+|\hat{l}|\neq 0,\quad |\acute{l}|+|\hat{l}|\le \mathcal{M}+2,\quad |\hat{l}|\le 2,$$
$$M_{\acute{l},\hat{l}}:=\max\{|j|: \acute{l}_j\neq 0\ {\rm or}\ \hat{l}_j\neq 0\},\quad C(\mathcal{N},\acute{l})=\mathcal{N}^{(|\acute{l}|+4)^2},$$
then we call the frequencies $\breve{\omega}$ and $\breve{\Omega}(x;\xi)$ are $(\acute{\eta},\mathcal{N},\mathcal{M})$-non-resonant.
Define the resonant sets $\mathcal{R}_{k\acute{l}\hat{l}}$ by
\begin{equation}\label{rlll}
\mathcal{R}_{k\acute{l}\hat{l}}=\left\{\xi\in\Pi_\eta: |\langle k,\breve{\omega}\rangle+\langle \acute{l},\acute{\Omega}\rangle+\langle \hat{l},[\hat{\Omega}]\rangle|\le \frac{\acute{\eta}M_{\acute{l},\hat{l}}}{4^{\mathcal{M}}(|k|+1)^{\tau}C(\mathcal{N},\acute{l})}\right\}.
\end{equation}
Let $$\mathcal{R}=\bigcup_{\substack{|k|+|\acute{l}|+|\hat{l}|\neq 0\\|\acute{l}|+|\hat{l}|\le \mathcal{M}+2,|\hat{l}|\le 2}}\mathcal{R}_{k\acute{l}\hat{l}},$$
and
$$\Pi_{\acute{\eta}}=\Pi_{\eta}\setminus{\mathcal{R}},$$
where $\Pi_{\eta}$ is obtained by Theorem \ref{normaltheorem}. Then it follows that for each $\xi\in \Pi_{\acute{\eta}}$, the frequencies $\breve{\omega}$ and $\breve{\Omega}(x;\xi)$ are $(\acute{\eta},\mathcal{N},\mathcal{M})$-non-resonant.

\begin{The}[Partial normal form of order $\mathcal{M}+2$]\label{parth}
Consider the normal form of order $2$
$$\breve{H}(x,y,z,\bar{z};\xi)=\breve{N}(x,y,z,\bar{z};\xi)+\breve{P}(x,y,z,\bar{z};\xi),$$
obtained in Theorem \ref{normaltheorem}. Given any positive integer $\mathcal{M}$ and $0<\acute{\eta}<1$, there exists
a small $\rho_0>0$ and
%a large positive integer $\mathcal{N}_0$
 depending on $s_0,r_0, n$ and $\mathcal{M}$.  %{\color{red}For each $\mathcal{N}>\mathcal{N}_0$, there exists a small $\rho$ satisfying
%$$0<\rho<\frac{\acute{\eta}^2}{2}\mathcal{N}^{-2C_0(\mathcal{M}+7)^2},$$}
For each $0<\rho<\rho_0$ and any integer $\mathcal{N}$ satisfying
$$\mathcal{N}_0<\mathcal{N}<(\frac{\acute{\eta}}{2\rho})^{1/2C_0(\mathcal{M}+7)^2},$$
(where $C_0$ is a constant depending on $n$ and $\tau$) and for each $\xi\in \Pi_{\acute{\eta}}$, there is a symplectic map
$$\Phi:D(4\varrho,4\rho,4\rho)\rightarrow D(5\varrho, 5\rho,5\rho),$$
(where $\varrho=\frac{\mathcal{M}}{\mathcal{N}^2}$) such that
\begin{align}\label{pnf}
\breve{\breve{H}}(x,y,z,\bar{z};\xi)=\breve{H}\circ\Phi=&N(x,y,z,\bar{z};\xi)+Z(x,y,z,\bar{z};\xi)\nonumber\\
&+Q(x,y,z,\bar{z};\xi)+T(x,y,z,\bar{z};\xi)
\end{align}
where
$$N(x,y,z,\bar{z};\xi)=\sum_{j=1}^n\breve{\omega}_j(\xi)y_j+\sum_{j\in\mathbb{Z}_*}\frac{\breve{\Omega}_j(x;\xi)}{j}z_j\bar{z}_j,$$
$$Z(x,y,z,\bar{z};\xi)=\sum_{4\le j\le\mathcal{M}+2}Z_j(x,y,z,\bar{z};\xi),$$
%$$A(x,y,z,\bar{z};\xi)=\sum_{4\le j\le\mathcal{M}+2}A_j(x,y,z,\bar{z};\xi)$$
%$$R(x,y,z,\bar{z};\xi)=\sum_{ j\geq\mathcal{M}+3}R_j(x,y,z,\bar{z};\xi),$$
$$Q(x,y,z,\bar{z};\xi)=\sum_{3\le j\le \mathcal{M}+2}Q_j(x,y,z,\bar{z};\xi),$$
and
\begin{equation}\label{T}
T(x,y,z,\bar{z};\xi)=\sum_{j\geq 3}T_j(x,y,z,\bar{z};\xi),
\end{equation}
with
$$Z_j(x,y,z,\bar{z};\xi)=\sum_{2|\alpha|+2|\beta|+2|\mu|=j, |\mu|=0 }Z^{\alpha\beta\beta\mu\mu}(\xi)y^{\alpha}\acute{z}^{\beta}\bar{\acute{z}}^{\beta}\hat{z}^{\mu}\bar{\hat{z}}^{\mu},$$
$$+\sum_{2|\alpha|+2|\beta|+2|\mu|=j, |\mu|=1 }Z^{\alpha\beta\beta\mu\mu}(x;\xi)y^{\alpha}\acute{z}^{\beta}\bar{\acute{z}}^{\beta}\hat{z}^{\mu}\bar{\hat{z}}^{\mu},\quad 4\le j\le\mathcal{M}+2,$$
%$$A_j(x,y,z,\bar{z};\xi)=\sum_{2|\alpha|+|\beta|+|\gamma|+|\mu|+|\nu|=j, |\mu|+|\nu|\le2 }A^{\alpha\beta\gamma\mu\nu}(x;\xi)y^{\alpha}\acute{z}^{\beta}\bar{\acute{z}}^{\gamma}\hat{z}^{\mu}\bar{\hat{z}}^{\nu},$$
%$$R_j(x,y,z,\bar{z};\xi)=\sum_{ 2|\alpha|+|\beta|+|\gamma|+|\mu|+|\nu|=j, |\mu|+|\nu|\le2 }R^{\alpha\beta\gamma\mu\nu}(x;\xi)y^{\alpha}\acute{z}^{\beta}\bar{\acute{z}}^{\gamma}\hat{z}^{\mu}\bar{\hat{z}}^{\nu},$$
$$Q_j(x,y,z,\bar{z};\xi)=\sum_{2|\alpha|+|\beta|+|\gamma|+|\mu|+|\nu|=j,|\mu|+|\nu|\geq3 }Q^{\alpha\beta\gamma\mu\nu}(x;\xi)y^{\alpha}\acute{z}^{\beta}\bar{\acute{z}}^{\gamma}\hat{z}^{\mu}\bar{\hat{z}}^{\nu},\quad 3\le j \le \mathcal M+2,$$
$$T_j(x,y,z,\bar{z};\xi)=\sum_{2|\alpha|+|\beta|+|\gamma|+|\mu|+|\nu|=j}T^{\alpha\beta\gamma\mu\nu}(x;\xi)y^{\alpha}\acute{z}^{\beta}\bar{\acute{z}}^{\gamma}\hat{z}^{\mu}\bar{\hat{z}}^{\nu},\quad j\geq\mathcal{M}+3,$$
and
$$T_j(x,y,z,\bar{z};\xi)=\sum_{\substack{2|\alpha|+|\beta|+|\gamma|+|\mu|+|\nu|=j\\ \beta\neq \gamma,\mu\neq \nu,|\mu|+|\nu|\le 2}}T^{\alpha\beta\gamma\mu\nu}(x;\xi)y^{\alpha}\acute{z}^{\beta}\bar{\acute{z}}^{\gamma}\hat{z}^{\mu}\bar{\hat{z}}^{\nu},\quad 3\le j\le \mathcal{M}+2.$$

Moreover, the following estimates hold:

(1) the symplectic map $\Psi$ satisfies
\begin{equation}\label{lho2}
\Vert \Psi-id\Vert_{p,D(4\varrho,4\rho,4\rho)}\le \frac{c\rho \mathcal{N}^C}{\acute{\eta}^2},
\end{equation}
$$\Vert D\Psi-Id\Vert_{p,D(4\varrho,4\rho,4\rho)}\le \frac{c \mathcal{N}^C}{\acute{\eta}^2},$$
where $C$ is a constant depending on $n$ and $\tau$.

(2) The functions ${Z}$, $Q$ and $T$ satisfy
$$|||Q|||_{D(4\varrho,4\rho,4\rho)\times \Pi_{\acute{\eta}}}\le c\rho^3,$$
$$%|||Z|||_{D(4\varrho,4\rho,4\rho)\times \Pi_{\acute{\eta}}}
\Vert X_Z\Vert_{p-1, D(4\varrho,4\rho,4\rho)\times \Pi_{\acute{\eta}}}\le c\rho^3\left(\frac{\mathcal{N}^{2C_0(\mathcal{M}+6)^2}}{\acute{\eta}^2}\rho\right),$$
%$$|||R|||_{D(4\varrho,4\rho,4\rho)\times \Pi_\eta}\le c\rho^3\left(\frac{\mathcal{N}^{2(\mathcal{M}+6)^2+2C_0}}{\acute{\eta}^2}\rho\right)^{\mathcal{M}},$$
and
\begin{equation}\label{TT}
\Vert X_T\Vert_{p-1, D(4\varrho,4\rho,4\rho)\times \Pi_{\acute{\eta}}}\le c\rho^3\left(\frac{\mathcal{N}^{2C_0(\mathcal{M}+6)^2}}{\acute{\eta}^2}\rho\right)^{\mathcal{M}},
\end{equation}
%(2)The functions ${Z}_j, R_j$ and $Q_j$ satisfy
%$$|||Z_j|||_{D(s_0/2)\times \Pi_\eta}\le c\rho\left(\frac{\mathcal{N}^{2(\mathcal{M}+6)^2+2C(\tau+n)}}{\acute{\eta}^2}\rho\right)^{j-3},$$
%$$|||R_j|||_{D(s_0/2)\times \Pi_\eta}\le c\rho\left(\frac{\mathcal{N}^{2(\mathcal{M}+6)^2+2C(\tau+n)}}{\acute{\eta}^2}\rho\right)^{j-3},$$
%$$|||Q_j|||_{D(s_0/2)\times \Pi_\eta}\le c\rho\left(\frac{\mathcal{N}^{2(\mathcal{M}+6)^2+2C(\tau+n)}}{\acute{\eta}^2}\rho\right)^{j-3}.$$
where $c$ is a constant depending on $s_0,r_0, n$, and $\mathcal{M}$.
\end{The}

\begin{Rem}
From \eqref{T}, it is easy to see that the function $T$ mainly consists of two parts. One part has the monomials of the form $T^{\alpha\beta\gamma\mu\nu}(x;\xi)y^{\alpha}\acute{z}^{\beta}\bar{\acute{z}}^{\gamma}\hat{z}^{\mu}\bar{\hat{z}}^{\nu}$ with $2|\alpha|+|\beta|+|\gamma|+|\mu|+|\nu|\geq \mathcal{M}+3$, we denote it by $\mathcal A_1$. Another part contains the monomials of the form $T^{\alpha\beta\gamma\mu\nu}(x;\xi)y^{\alpha}\acute{z}^{\beta}\bar{\acute{z}}^{\gamma}\hat{z}^{\mu}\bar{\hat{z}}^{\nu}$ with $3\le 2|\alpha|+|\beta|+|\gamma|+|\mu|+|\nu|\le \mathcal{M}+2$, we denote it by $\mathcal A_2$. In addition, the norm of the part $\mathcal A_1$ and $\mathcal A_2$ has the estimate \eqref{TT}. It means that $T$ is small enough for the long time stability estimate.
\end{Rem}

\subsection{Iterative Lemma}
Set
$$\varrho^{\prime}=\frac{\varrho}{12\mathcal{M}}:=\frac{1}{12\mathcal{N}^2},\quad \rho^\prime=\frac{\rho}{12\mathcal{M}}.$$
Let $2\le j_0\le \mathcal{M}+2$ and denote
$$D_{j_0}=D(5\varrho-2(j_0-2)\varrho^\prime, 5\rho-2(j_0-2)\rho^\prime, 5\rho-2(j_0-2)\rho^\prime),$$
$$D^\prime_{j_0}=D(5\varrho-(2(j_0-2)+1)\varrho^\prime, 5\rho-(2(j_0-2)+1)\rho^\prime, 5\rho-(2(j_0-2)+1)\rho^\prime),$$
$$D^{\prime\prime}_{j_0}=D(5\varrho-(2(j_0-2)+2)\varrho^\prime, 5\rho-(2(j_0-2)+2)\rho^\prime, 5\rho-(2(j_0-2)+2)\rho^\prime),$$
and
$$D_{j_0+1}\subset D^{\prime\prime}_{j_0}\subset D^\prime_{j_0}\subset D_{j_0}.$$
\begin{Lem}\label{iterm}
Consider the partial normal form of order $j_0\ (2\le j_0\le\mathcal{M}+1)$
\begin{equation*}
  \begin{split}
  H_{j_0}(x,y,z,\bar{z};\xi) =&\breve{N}(x,y,z,\bar{z};\xi)+Z_{j_0}(y,z,\bar{z};\xi) +R_{j_0}(x,y,z,\bar{z};\xi)\\
       & +Q_{j_0}(x,y,z,\bar{z};\xi)+T_{j_0}(x,y,z,\bar{z};\xi),
   \end{split}
\end{equation*}

where
\begin{align}
&\breve{N}(x,y,z,\bar{z};\xi)=\sum_{j=1}^n\breve{\omega}_j(\xi)y_j+\sum_{j\in\mathbb{Z}_*}\frac{\breve{\Omega}_j(x;\xi)}{j}z_j\bar{z}_j,\nonumber\\
%&=\sum_{j=1}^n\breve{\omega}_j(\xi)y_j+\sum_{|j|\le\mathcal{N}}\frac{\breve{\Omega}_j(\xi)}{j}z_j\bar{z}_j+\sum_{|j|\geq\mathcal{N}+1}\frac{\breve{\Omega}_j(x;\xi)}{j}z_j\bar{z}_j,\nonumber\\
&Z_{j_0}(x,y,z,\bar{z};\xi)=\sum_{3\le j\le j_0}Z_{j_0j}(x,y,z,\bar{z};\xi),\label{zj0}\\
&R_{j_0}(x,y,z,\bar{z};\xi)=\sum_{ j_0+1\le j\le \mathcal{M}+2}R_{j_0j}(x,y,z,\bar{z};\xi),\nonumber\\
&Q_{j_0}(x,y,z,\bar{z};\xi)=\sum_{3\le j\le\mathcal{M}+2}Q_{j_0j}(x,y,z,\bar{z};\xi)\nonumber\\
&T_{j_0}(x,y,z,\bar{z};\xi)=\sum_{j\geq 3}T_{j_0j}(x,y,z,\bar{z};\xi),\label{zj01}
\end{align}
with
\begin{equation*}
  \begin{split}
     Z_{j_0j}(y,z,\bar{z};\xi) &=\sum_{ 2|\alpha|+2|\beta|+2|\mu|=j, |\mu|=0}Z_{j_0}^{\alpha\beta\beta\mu\mu}(\xi)y^{\alpha}\acute{z}^{\beta}\bar{\acute{z}}^{\beta}\hat{z}^{\mu}\bar{\hat{z}}^{\mu}\\
       & =\sum_{ 2|\alpha|+2|\beta|+2|\mu|=j, |\mu|=1 }Z_{j_0}^{\alpha\beta\beta\mu\mu}(x;\xi)y^{\alpha}\acute{z}^{\beta}\bar{\acute{z}}^{\beta}\hat{z}^{\mu}\bar{\hat{z}}^{\mu},
   \end{split}
\end{equation*}

$$R_{j_0j}(x,y,z,\bar{z};\xi)=\sum_{2|\alpha|+|\beta|+|\gamma|+|\mu|+|\nu|=j, |\mu|+|\nu|\le2 }R_{j_0}^{\alpha\beta\gamma\mu\nu}(x;\xi)y^{\alpha}\acute{z}^{\beta}\bar{\acute{z}}^{\gamma}\hat{z}^{\mu}\bar{\hat{z}}^{\nu},$$
$$Q_{j_0j}(x,y,z,\bar{z};\xi)=\sum_{2|\alpha|+|\beta|+|\gamma|+|\mu|+|\nu|=j,|\mu|+|\nu|\geq3 }Q_{j_0}^{\alpha\beta\gamma\mu\nu}(x;\xi)y^{\alpha}\acute{z}^{\beta}\bar{\acute{z}}^{\gamma}\hat{z}^{\mu}\bar{\hat{z}}^{\nu},$$
$$T_{j_0j}(x,y,z,\bar{z};\xi)=\sum_{2|\alpha|+|\beta|+|\gamma|+|\mu|+|\nu|=j }T_{j_0}^{\alpha\beta\gamma\mu\nu}(x;\xi)y^{\alpha}\acute{z}^{\beta}\bar{\acute{z}}^{\gamma}\hat{z}^{\mu}\bar{\hat{z}}^{\nu},\quad j\geq \mathcal{M}+3,$$
and
$$T_{j_0j}(x,y,z,\bar{z};\xi)=\sum_{\substack{2|\alpha|+|\beta|+|\gamma|+|\mu|+|\nu|=j \\\beta\neq \gamma, \mu\neq\nu, |\mu|+|\nu|\le 2}}T_{j_0}^{\alpha\beta\gamma\mu\nu}(x;\xi)y^{\alpha}\acute{z}^{\beta}\bar{\acute{z}}^{\gamma}\hat{z}^{\mu}\bar{\hat{z}}^{\nu},\quad 3\le j\le \mathcal{M}+2.$$
Suppose $Z_{j_0j}(x,y,z,\bar{z};\xi)$, $R_{j_0j}(x,y,z,\bar{z};\xi)$ and $Q_{j_0j}(x,y,z,\bar{z};\xi)$ satisfy the following estimates
\begin{align}
&\Vert X_{Z_{j_0j}}\Vert_{p-1, D_{j_0}\times \Pi_{\acute{\eta}}}, \ \frac{1}{\rho^2}|||{Z_{j_0j}}|||_{D_{j_0}\times \Pi_{\acute{\eta}}}\lessdot \rho\left(\frac{\mathcal{N}^{2C_0(j_0+4)^2}}{\acute{\eta}^2}\rho\right)^{j-3},\label{xz}\\
 &\Vert X_{R_{j_0j}}\Vert_{p-1, D_{j_0}\times \Pi_{\acute{\eta}}}, \ \frac{1}{\rho^2}|||{R_{j_0j}}|||_{D_{j_0}\times \Pi_{\acute{\eta}}}\lessdot \rho\left(\frac{\mathcal{N}^{2C_0(j_0+5)^2}}{\acute{\eta}^2}\rho\right)^{j-3},\label{xr}\\
 &\Vert X_{Q_{j_0j}}\Vert_{p-1, D_{j_0}\times \Pi_{\acute{\eta}}},\ \frac{1}{\rho^2}|||{Q_{j_0j}}|||_{D_{j_0}\times \Pi_{\acute{\eta}}}\lessdot \rho\left(\frac{\mathcal{N}^{2C_0(j_0+5)^2}}{\acute{\eta}^2}\rho\right)^{j-3},\label{xy}
% &\Vert A_{j_0}\Vert_{p,D(4\alpha,4\rho,4\rho)\times \Pi_\eta}\lessdot \rho^3\left(\frac{\mathcal{N}^{2(j_0+5)^2+2C(\tau+n)}}{\acute{\eta}^2}\rho\right)^{\mathcal{M}}.
 \end{align}
 where $C_0>40(\tau+n)$ is a constant depending on $n$ and $ \tau$.
 The function $T_{j_0}(x,y,z,\bar{z};\xi)$ is already small enough for the long time stability estimate
 \begin{equation}\label{4.6+}
\Vert X_{T_{j_0}}\Vert_{p-1,D_{j_0}\times \Pi_{\acute{\eta}}}\lessdot \rho \left(\frac{\mathcal{N}^{2C_0(j_0+5)^2}}{\acute{\eta}^2}\rho\right)^{\mathcal{M}}.
\end{equation}
Then there exists a symplectic map $\Psi_{j_0}:D_{j_0+1}\rightarrow D_{j_0}$  satisfying the estimate
 \begin{equation}\label{psi1}
\Vert \Psi_{j_0}-id\Vert_{p,D_{j_0}^\prime}\lessdot \left(\frac{\mathcal{N}^{2C_0(j_0+5)^2}}{\acute{\eta}^2}\rho\right)^{j_0-1},
\end{equation}
\begin{equation}\label{psi2}
|| D\Psi_{j_0}-Id||_{p,D_{j_0+1}}\lessdot \rho^{-1}\left(\frac{\mathcal{N}^{2C_0(j_0+5)^2}}{\acute{\eta}^2}\rho\right)^{j_0-1},
\end{equation}
such that
\begin{equation*}
  \begin{split}
     \breve{\breve{H}}_{j_0+1}(x,y,z,\bar{z};\xi) =&H_{j_0}\circ\Psi_{j_0}(x,y,z,\bar{z};\xi)=\breve{N}(x,y,z,\bar{z};\xi)+Z_{j_0+1}(x,y,z,\bar{z};\xi) \\
     &+R_{j_0+1}(x,y,z,\bar{z};\xi)+Q_{j_0+1}(x,y,z,\bar{z};\xi)+T_{j_0+1}(x,y,z,\bar{z};\xi),
   \end{split}
\end{equation*}
where $Z_{j_0+1}, R_{j_0+1}$, $Q_{j_0+1}$ and $T_{j_0+1}$ satisfy \eqref{xz}, \eqref{xr},\eqref{xy} and \eqref{4.6+} by replacing $j_0$ by $j_{0}+1$ there.
\end{Lem}

\subsection{The derivation of homological equations}
Let
$$F_{j_0}(x,y,z,\bar{z};\xi)=\sum_{2|\alpha|+|\beta|+|\gamma|+|\mu|+|\nu|=j_0+1, |\mu|+|\nu|\le2 }F_{j_0}^{\alpha\beta\gamma\mu\nu}(x;\xi)y^{\alpha}\acute{z}^{\beta}\bar{\acute{z}}^{\gamma}\hat{z}^{\mu}\bar{\hat{z}}^{\nu},$$
and $\Psi_{j_0}=X_{F_{j_0}}^t|_{t=1}$ be the time-$1$-map of the Hamiltonian vector field $X_{F_{j_0}}$.

According to Taylor's formula, we have
\begin{align*}
H_{j_0+1}=&H_{j_0}\circ \Psi_{j_0}\\
=&(\breve{N}+Z_{j_0}+R_{j_0}+Q_{j_0}+T_{j_0})\circ X_{F_{j_0}}^t|_{t=1}\\
=&\breve{N}+\{\breve{N},F_{j_0}\}+\int_0^1(1-t)\{\{\breve{N},F_{j_0}\},F_{j_0}\}\circ X_{F_{j_0}}^tdt\\
&+R_{j_0(j_0+1)}+\int_0^1\{R_{j_0(j_0+1)},F_{j_0}\}\circ X_{F_{j_0}}^tdt\\
&+(Z_{j_0}+R_{j_0}-R_{j_0(j_0+1)}+Q_{j_0}+T_{j_0})\circ X_{F_{j_0}}^t|_{t=1}.
\end{align*}
We need to solve the homological equation
\begin{equation}\label{hh}
\{\breve{N},F_{j_0}\}+R_{j_0(j_0+1)}=\widehat{Z_{j_0}}+\widehat{T_{j_0}},
\end{equation}
where
\begin{equation}\label{hatz0}
\widehat{T_{j_0}}=\{\breve{N},(1-\Gamma_K)G_{j_0}\},
\end{equation}
\begin{align}
\widehat{Z_{j_0}}=&\sum_{2|\alpha|+2|\beta|+2|\mu|=j_{0}+1, |\mu|=0}R_{j_0}^{\alpha\beta\beta\mu\mu}(0;\xi)y^{\alpha}\acute{z}^{\beta}\bar{\acute{z}}^{\beta}\hat{z}^{\mu}\bar{\hat{z}}^{\mu}\label{hatz}\\
&+\sum_{2|\alpha|+2|\beta|+2|\mu|=j_{0}+1, |\mu|=1 }R_{j_0}^{\alpha\beta\beta\mu\mu}(x;\xi)y^{\alpha}\acute{z}^{\beta}\bar{\acute{z}}^{\beta}\hat{z}^{\mu}\bar{\hat{z}}^{\mu},\label{hatz1}
\end{align}
 with $R_{j_0}^{\alpha\beta\beta\mu\mu}(0;\xi)$ being the $0$-th Fourier series of $R_{j_0}^{\alpha\beta\beta\mu\mu}(x;\xi)$ and $G_{j_0}$ being the solution of equation \eqref{hh} without the term $\widehat{T_{j_0}}$.
After solving the homological equation \eqref{hh}, we have
$$H_{j_0+1}=\breve{N}+Z_{j_0+1}+R_{j_0+1}+Q_{j_0+1}+T_{j_0+1},$$
with
\begin{equation}\label{zj0+}
Z_{j_0+1}=Z_{j_0}+\widehat{Z_{j_0}},
\end{equation}
%\begin{equation}\label{zj0+1}
%A_{j_0+1}=A_{j_0}+\widehat{A_{j_0}},
%\end{equation}
and
\begin{align}
&R_{j_0+1}+Q_{j_0+1}+T_{j_0+1}\\
=&\widehat{T_{j_0}}+\int_0^1(1-t)\{\{\breve{N},F_{j_0}\},F_{j_0}\}\circ X_{F_{j_0}}^tdt\nonumber\\
&+\int_0^1\{R_{j_0(j_0+1)}+Z_{j_0},F_{j_0}\}\circ X_{F_{j_0}}^tdt\label{2u}\\
&+(R_{j_0}-R_{j_0(j_0+1)}+Q_{j_0}+T_{j_0})\circ X_{F_{j_0}}^t|_{t=1}.\nonumber
\end{align}
\subsection{The solution of homological equation}
The following lemma plays a key role in solving the homological equation.
\begin{Lem}{(\cite{liu2010spectrum})}\label{ylcp}
Consider the first-order partial differential equation \eqref{hh}
\begin{equation}\label{hoylcp}
{\bf i} \partial_{\omega} u+\lambda u+\mu(\phi)u=p(\phi),\quad \phi\in D(s), \lambda\in \mathbb{C},
\end{equation}
for the unknown function $u$ defined on the torus $D(s)$. Assume the following conditions hold.

(i) There are constants $\eta_0, \eta>0$ and $\tau>n$ such that
\begin{align}
&|\langle k,\omega\rangle|\geq \frac{\eta_0}{|k|^\tau},\quad \forall k\in\mathbb{Z}^n\setminus\{0\},\nonumber\\
&|\langle k,\omega\rangle+\lambda|\geq \frac{\eta\gamma }{1+|k|^\tau},\quad \forall k\in\mathbb{Z}^n.\label{sdi}
\end{align}

(ii) $\mu(\phi):D(s)\rightarrow \mathbb{C}$ is real analytic and zero average, $\int_{\mathbb{T}^n}\mu(\phi)d\phi=0.$ Moreover, assume
\begin{equation}\label{yuan1}|\mu|_{D(s),\tau+1}:=\sum_{k\in\mathbb{Z}^n}|\hat{\mu}_k||k|^{\tau+1}e^{|k|s}\le C\alpha,\end{equation}
where $\hat{\mu}_k$ is the $k$-Fourier coefficient of $\mu$.

(iii) $p$ is real analytic in $D(s)$.

Then homological equation \eqref{hoylcp} has a unique solution $u(\phi)$ defined in a narrower domain $D(s-\sigma)$ with $0<\sigma<s$, which satisfies
\begin{equation}\label{yuan3}\Vert u\Vert_{D(s-\sigma)}\le \frac{c(n,\tau)}{\eta\gamma  \sigma^{n+\tau}}e^{2C\alpha s/\eta_0}\Vert p\Vert_{D(s)}.\end{equation}
\end{Lem}

\begin{Rem}
The present lemma is a variant version of  Liu-Yuan's result (Theorem 1.4 of \cite{liu2010spectrum}). In \cite{liu2010spectrum}, the parameter $\alpha$ in both \eqref{yuan1} and \eqref{yuan3} is restricted to $\alpha=\gamma$ where the parameter $\gamma$ is that in the small divisor conditions \eqref{sdi}. By checking the proof of Theorem 1.4 of \cite{liu2010spectrum}, we find that the restrict condition $\alpha=\gamma$ can be dropped.
\end{Rem}

Now, we return to the homological equation \eqref{hh}. Writing $F_{j_0}=\Gamma_KG_{j_0}$, by \eqref{hatz0}, equation \eqref{hh} becomes
$$\{\breve{N},G_{j_0}\}+R_{j_0(j_0+1)}=\widehat{Z_{j_0}}.$$

Comparing the coefficients, equation \eqref{hh} reads
 \begin{equation}\label{hurg}
 (\partial_\omega+{\bf i}\langle \beta-\gamma, \acute{\Omega}\rangle +{\bf i}\langle \mu-\nu,\hat{\Omega}\rangle) G_{j_0}^{\alpha\beta\gamma\mu\nu}=R_{j_0}^{\alpha\beta\gamma\mu\nu},
 \end{equation}
 where $2|\alpha|+|\beta|+|\gamma|+|\mu|+|\nu|=j_0+1\le \mathcal{M}+2, |\mu|+|\nu|\le2$.
%By \eqref{iii}, we have
%$$|\langle \beta-\gamma, [\breve{\Omega}]\rangle|\le (|\beta|+|\gamma|)|[\breve{\Omega}]|\le (\mathcal{M}+2)|[\breve{\Omega}]|\le 2(\mathcal{M}+2)^3.$$
%Consider the term $\langle \mu-\nu,\hat{\Omega}\rangle$, it will have the following cases: $\hat{\Omega}_j$, $\hat{\Omega}_i-\hat{\Omega}_j$ and $\hat{\Omega}_i+\hat{\Omega}_j$ with $i,j\geq \mathcal{M}+3$.

%For case 1: $\hat{\Omega}_j$, if $|j|\geq 2(\mathcal{M}+2)^3$, then
%$$|\langle \beta-\gamma, [\breve{\Omega}]\rangle\pm[\hat{\Omega}_j]|\geq |j(j-1)|-2(\mathcal{M}+2)^3>|j(j-2)|>|j|.$$
%Therefore
%$$|\langle \beta-\gamma, [\breve{\Omega}]\rangle\pm [\hat{\Omega}_j]|<1/|j|,$$
%which implies
%$$\Vert(\langle \beta-\gamma, [\breve{\Omega}]\rangle\pm[\hat{\Omega}_j])^{-1}(\langle \beta-\gamma, \tilde{\breve{\Omega}}\rangle \pm\tilde{\hat{\Omega}}_j)\Vert \le c\gamma^6\epsilon.$$
%If $|j|<  2(\mathcal{M}+2)^3$, then by non-resonant condition, we have
%$$\Vert(\langle \beta-\gamma, [\breve{\Omega}]\rangle\pm[\hat{\Omega}_j])^{-1}(\langle \beta-\gamma, \tilde{\breve{\Omega}}\rangle \pm\tilde{\hat{\Omega}}_j)\Vert \le c\gamma^6\epsilon.$$
%By non-resonant condition, we have
%$$|\langle \beta-\gamma,[\acute{\Omega}]\rangle+\langle \mu-\nu,[\hat{\Omega}]\rangle|\geq \frac{\acute{\eta}M_{\beta-\gamma,\mu-\nu}}{4^{\mathcal{M}}C(\mathcal{N},\beta-\gamma)}.$$
%Therefore, we have
%$$\Vert\big(\langle \beta-\gamma,\acute{\Omega}\rangle+\langle \mu-\nu,[\hat{\Omega}]\rangle\big)^{-1}\langle \mu-\nu,\hat{\Omega}\rangle\Vert\le  \frac{4^{\mathcal{M}}C(\mathcal{N},\beta-\gamma)}{\acute{\eta}}c\eta^6\varepsilon.$$
In order to solve this homological equation, we split the term $\langle \mu-\nu,[\hat{\Omega}]\rangle$ into the following three cases:

\noindent {\bf Case 1}: $|\mu-\nu|=1$ and $|\langle \mu-\nu,[\hat{\Omega}]\rangle|=|[\breve{\Omega}_i]|$, for some $i\geq \mathcal{N}+1.$
Then, by assumption {\bf{A}} and \eqref{huomom}, there exists some constant $c_1>0$ such that
$$|\langle \mu-\nu,[\hat{\Omega}]\rangle|\geq \frac{c_1}{2}i^2\geq \frac{c_1}{2}i,$$
for some constant $c_1>0$.

\noindent{\bf Case 2}: $|\mu-\nu|=2$ and $|\langle \mu-\nu,[\hat{\Omega}]\rangle|=|[\breve{\Omega}_i]+[\breve{\Omega}_j]|$ or $|\langle \mu-\nu,[\hat{\Omega}]\rangle|=|[\breve{\Omega}_i]-[\breve{\Omega}_j]|$, ($ i\neq \pm j$) for some $i, j\geq \mathcal{N}+1.$
Then, we obtain
$$|\langle \mu-\nu,[\hat{\Omega}]\rangle|\geq c_1(i^2\pm j^2)-(i+j)c\eta^6 \epsilon\geq \frac{c_1}{2}|i+j||i-j|\geq \frac{c_1}{2}\max\{|i|,|j|\}.$$
In Cases 1 and 2, we assume $|\acute{\Omega}_j(x,\xi)|\le c_2|j|^2\le c_2 \mathcal{N}^2$. Thus, if
$$\max\{|i|,|j|\}\geq \frac{8}{c_1}c_2(\mathcal{M}+2)\mathcal{N}^2,$$
 then
\begin{align*}
|\langle \beta-\gamma,[\acute{\Omega}]\rangle+\langle \mu-\nu,[\hat{\Omega}]\rangle|\geq& |\langle \mu-\nu,[\hat{\Omega}]\rangle|-|\langle \beta-\gamma,[\acute{\Omega}]\rangle|\\
\geq &\frac{c_1}{4}\max\{|i|,|j|\}+ 2c_2(\mathcal{M}+2)\mathcal{N}^2\\
&-(1+c\eta^6\varepsilon)\cdot c_2|\acute{l}|\mathcal{N}^2\\
\geq& \frac{c_1}{4}\max\{|i|,|j|\}.
\end{align*}
Thus, in cases 1 and 2, we have the following two subcases:

{\bf Subcase 1}: If $\max\{|i|,|j|\}\geq \frac{8}{c_1}c_2(\mathcal{M}+2)\mathcal{N}^2$, then we have that, for $(\mathcal{M}+2)c\eta^6\le \varepsilon^{-1/2}$,
$$\Vert a(x)\Vert_{D(s_0/2)}:=\Vert\lambda^{-1}(\langle \beta-\gamma,\acute{\Omega}\rangle+\langle \mu-\nu,\hat{\Omega}\rangle-\lambda)\Vert_{D(s_0/2)}\le (\mathcal{M}+2)c\eta^6\varepsilon\le \varepsilon^{1/2},$$
where $\lambda:=\langle \beta-\gamma,[\acute{\Omega}]\rangle+\langle \mu-\nu,[\hat{\Omega}]\rangle$. Thus, equation \eqref{hurg} can be written as
$${\bf i} \partial_{\omega}G_{j_0}^{\alpha\beta\gamma\mu\nu} +\lambda (1+a(x))G_{j_0}^{\alpha\beta\gamma\mu\nu}=R_{j_0}^{\alpha\beta\gamma\mu\nu}.$$
 By Lemma \ref{lemma1}, one gets
$$| G_{j_0}^{\alpha\beta\gamma\mu\nu}|_{D^\prime_{j_0}}\le \frac{\mathcal{N}^{(j_0+5)^2}}{\acute{\eta}M_{\beta-\gamma,\mu-\nu}(\varrho^\prime)^{20(\tau+n)}}\Vert R_{j_0}^{\alpha\beta\gamma\mu\nu}\Vert_{D_{j_0}}$$
\begin{equation}\label{6y}
\le \frac{\mathcal{N}^{(j_0+5)^2+20(\tau+n)}}{\acute{\eta}M_{\beta-\gamma,\mu-\nu}}\Vert R_{j_0}^{\alpha\beta\gamma\mu\nu}\Vert_{D_{j_0}}.
\end{equation}
%$$\le \frac{1}{M_{\beta-\gamma,\mu-\nu}}\left(\frac{\mathcal{N}^{2(j_0+5)^2+2C(\tau+n)}}{\acute{\eta}^2}\right)^{j_0-1}.$$
In order to obtain the estimate for $\partial_{\xi}F_{j_0}^{\alpha\beta\gamma\mu\nu}$, by taking the partial derivative with respect to $\xi$ in both sides of the homological equation \eqref{hurg}, we have
\begin{align*}
&(\partial_\omega+{\bf i}\langle \beta-\gamma, \acute{\Omega}\rangle +{\bf i}\langle \mu-\nu,\hat{\Omega}\rangle) \partial_{\xi}G_{j_0}^{\alpha\beta\gamma\mu\nu}\\
=&\partial_{\xi}R_{j_0}^{\alpha\beta\gamma\mu\nu}-\partial_\xi \omega \cdot \partial_x G_{j_0}^{\alpha\beta\gamma\mu\nu}-({\bf i}\langle \beta-\gamma, \partial_\xi\acute{\Omega}\rangle +{\bf i}\langle \mu-\nu,\partial_\xi\hat{\Omega}\rangle)\cdot G_{j_0}^{\alpha\beta\gamma\mu\nu}.
\end{align*}
Similar to \eqref{6y}, one can obtain
\begin{equation}\label{hufj1}
| \partial_\xi G_{j_0}^{\alpha\beta\gamma\mu\nu}|_{D^{\prime\prime}_{j_0}}\lessdot \frac{\mathcal{N}^{2(j_0+5)^2+40(\tau+n)}}{\acute{\eta}^2M_{\beta-\gamma,\mu-\nu}}\Vert R_{j_0}^{\alpha\beta\gamma\mu\nu}\Vert_{D_{j_0}}.
\end{equation}
From \eqref{6y} and \eqref{hufj1}, one immediately obtains
$$\Vert G_{j_0}^{\alpha\beta\gamma\mu\nu}\Vert_{D^{\prime\prime}_{j_0}\times \Pi_{\acute{\eta}}}\lessdot \frac{\mathcal{N}^{2(j_0+5)^2+40(\tau+n)}}{\acute{\eta}^2M_{\beta-\gamma,\mu-\nu}}\Vert R_{j_0}^{\alpha\beta\gamma\mu\nu}\Vert_{D_{j_0}\times \Pi_{\acute{\eta}}}.$$
Since $F^{\alpha\beta\gamma\mu\nu}_{j_0}=\Gamma_KG^{\alpha\beta\gamma\mu\nu}_{j_0}$, we also have
\begin{equation}\label{hou1}
\Vert F_{j_0}^{\alpha\beta\gamma\mu\nu}\Vert_{D^{\prime\prime}_{j_0}\times \Pi_{\acute{\eta}}}\lessdot \frac{\mathcal{N}^{2(j_0+5)^2+40(\tau+n)}}{\acute{\eta}^2M_{\beta-\gamma,\mu-\nu}}\Vert R_{j_0}^{\alpha\beta\gamma\mu\nu}\Vert_{D_{j_0}\times \Pi_{\acute{\eta}}}.
\end{equation}

{\bf Subcase 2}: If $\max\{|i|,|j|\}< \frac{8}{c_1}c_2(\mathcal{M}+2)\mathcal{N}^2$, we can use Lemma \ref{ylcp} to solve equation \eqref{hurg}. Similar to subcase 1, we let $F^{\alpha\beta\gamma\mu\nu}_{j_0}=\Gamma_KG^{\alpha\beta\gamma\mu\nu}_{j_0}$. Letting $$\lambda={\bf i}\langle \beta-\gamma, [\acute{\Omega}]\rangle +{\bf i}\langle \mu-\nu,[\hat{\Omega}]\rangle,\quad\mu={\bf i}\langle \beta-\gamma, \acute{\Omega}-[\acute{\Omega}]\rangle +{\bf i}\langle \mu-\nu,\hat{\Omega}-[\hat{\Omega}]\rangle,$$
$u=G_{j_0}^{\alpha\beta\gamma\mu\nu}$ and $p=R_{j_0}^{\alpha\beta\gamma\mu\nu}$ in Lemma \ref{ylcp}, we obtain
$$\Vert \mu\Vert_{D(s_{j_0}),\tau+1}\le\Vert \mu\Vert_{D(s_0/2),\tau+1}\le C\mathcal{M}\mathcal{N}^2\varepsilon.$$
It immediately follows that
\begin{align*}
|G_{j_0}^{\alpha\beta\gamma\mu\nu}|_{D_{j_0}^\prime\times \Pi_{\acute{\eta}}}\lessdot& \frac{\mathcal{N}^{(j_0+5)^2}}{\acute{\eta}M_{\beta-\gamma,\mu-\nu}(\varrho^\prime)^{\tau+n}} e^{C\mathcal{M}\mathcal{N}^2\varepsilon s_{j_0}/\eta}\Vert R_{j_0}^{\alpha\beta\gamma\mu\nu}\Vert_{D_{j_0}\times \Pi_{\acute{\eta}}}\\
\le &\frac{\mathcal{N}^{(j_0+5)^2+20(\tau+n)}}{\acute{\eta} M_{\beta-\gamma,\mu-\nu}}\Vert R_{j_0}^{\alpha\beta\gamma\mu\nu}\Vert_{D_{j_0}\times \Pi_{\acute{\eta}}},
\end{align*}
since $C\mathcal{M}\mathcal{N}^2\varepsilon s_{j_0}/\eta\lessdot1$.
Similar to subcase 1, one can obtain the estimate for $ \partial_\xi G_{j_0}^{\alpha\beta\gamma\mu\nu}$ and prove \eqref{hou1}.

{\bf Case 3}. Let $|\mu-\nu|=2$ and $\langle \mu-\nu,[\hat{\Omega}]\rangle=[\hat{\Omega}]_i-[\hat{\Omega}]_j$ with $i=- j$.

%For $i=j$, the homological equation \eqref{hurg} can be solved by the similar subcase 2.

For $i=-j$, by the momentum conservation for $\breve{\Omega}_j-[\breve{\Omega}_j]=\sum_{k\neq 0}\breve{\Omega}_{j}(k;\xi)e^{{\bf i}\langle k,x\rangle}$, we have
$$ -\sum_{b=1}^nj_bk_b-2j+\sum_{i} i(\beta_i-\gamma_i)=0.$$
It implies that
$$|j|\le C_{J_n}|k|+\frac{1}{2}\mathcal{M}\mathcal{N}.$$
Therefore,
$$|\breve{\Omega}_j-[\breve{\Omega}_j]|_{D(s_0/2),\tau+1}=\sum_{0\neq k\in\mathbb{Z}^n}|\breve{\Omega}_{j}(k;\xi)||k|^{\tau+1}e^{|k|s_0/2}$$
$$\le C\sum_{0\neq k\in\mathbb{Z}^n}|j^{-1}\breve{\Omega}_{j}(k;\xi)||k|^{\tau+1}(|k|+\mathcal{M}\mathcal{N})e^{|k|s_0/2}\le C\epsilon\mathcal{M}\mathcal{N}. $$
Then, the homological equation \eqref{hurg} can be solved by the argument similar to subcase 2. In addition, the solution has the estimates similar to \eqref{hou1}.

{\bf Case 4}. $|\mu-\nu|=0$. For this case, we can solve the homological equation \eqref{hurg} as in subcase 2, and obtain the similar estimate for its solution. Since $|\mu-\nu|=0$, the homological equation becomes
$$ (\partial_\omega+{\bf i}\langle \beta-\gamma, \acute{\Omega}\rangle ) G_{j_0}^{\alpha\beta\gamma\mu\nu}=R_{j_0}^{\alpha\beta\gamma\mu\nu}.$$
 Setting $\lambda={\bf i}\langle \beta-\gamma, [\acute{\Omega}]\rangle$,  $\mu={\bf i}\langle \beta-\gamma, \acute{\Omega}-[\acute{\Omega}]\rangle$, $u=G_{j_0}^{\alpha\beta\gamma\mu\nu}$ and $p=R_{j_0}^{\alpha\beta\gamma\mu\nu}$ in Lemma \ref{ylcp}, we obtain $\Vert \mu\Vert_{D(s_{j_0}),\tau+1}\le\Vert \mu\Vert_{D(s_0/2),\tau+1}\le C\mathcal{M}\mathcal{N}\varepsilon$. It follows that
\begin{align*}
|G_{j_0}^{\alpha\beta\gamma\mu\nu}|_{D_{j_0}^\prime\times \Pi_{\acute{\eta}}}\lessdot& \frac{\mathcal{N}^{(j_0+5)^2}}{\acute{\eta}M_{\beta-\gamma,\mu-\nu}(\varrho^\prime)^{\tau+n}} e^{C\mathcal{M}\mathcal{N}\varepsilon s_{j_0}/\eta}\Vert R_{j_0}^{\alpha\beta\gamma\mu\nu}\Vert_{D_{j_0}\times \Pi_{\acute{\eta}}}\\
\le &\frac{\mathcal{N}^{(j_0+5)^2+20(\tau+n)}}{\acute{\eta} M_{\beta-\gamma,\mu-\nu}}\Vert R_{j_0}^{\alpha\beta\gamma\mu\nu}\Vert_{D_{j_0}\times \Pi_{\acute{\eta}}},
\end{align*}
since $C\mathcal{M}\mathcal{N}\varepsilon s_{j_0}/\eta\lessdot1$.
Similar to subcase 1, one can obtain the estimate for $ \partial_\xi G_{j_0}^{\alpha\beta\gamma\mu\nu}$ and prove \eqref{hou1}.

\medskip

Combining {\bf Cases 1, 2, 3} and {\bf 4} together, $F_{j_0}^{\alpha\beta\gamma\mu\nu}$ has the estimate
\begin{equation*}
\Vert F_{j_0}^{\alpha\beta\gamma\mu\nu}\Vert_{D^{\prime\prime}_{j_0}\times \Pi_{\acute{\eta}}}\lessdot \frac{\mathcal{N}^{2(j_0+5)^2+40(\tau+n)}}{\acute{\eta}^2M_{\beta-\gamma,\mu-\nu}}\Vert R_{j_0}^{\alpha\beta\gamma\mu\nu}\Vert_{D_{j_0}\times \Pi_{\acute{\eta}}}
\end{equation*}
for all $2|\alpha|+|\beta|+|\gamma|+|\mu|+|\nu|=j_0+1\le \mathcal{M}+2, |\mu|+|\nu|\le2$.

Therefore, from \eqref{xr}, we have
\begin{align}
|||F_{j_0}|||^*_{D^{\prime\prime}_{j_0}\times \Pi_{\acute{\eta}}}\lessdot& \frac{\mathcal{N}^{2(j_0+5)^2+40(\tau+n)}}{\acute{\eta}^2}|||R_{j_0}|||_{D_{j_0}\times\Pi_{\acute{\eta}}}\nonumber\\
\le &\frac{\mathcal{N}^{2(j_0+5)^2+40(\tau+n)}}{\acute{\eta}^2}\cdot \rho^3\left(\frac{\mathcal{N}^{2C_0(j_0+5)^2}}{\acute{\eta}^2}\rho\right)^{j_0-2}\nonumber\\
\le &\rho^2\left(\frac{\mathcal{N}^{2C_0(j_0+5)^2}}{\acute{\eta}^2}\rho\right)^{j_0-1},\label{fujb}
\end{align}
and
\begin{align}\label{pxf}
\Vert X_{F_{j_0}}\Vert_{p,D^{\prime\prime}_{j_0}\times \Pi_{\acute{\eta}}}\lessdot& \frac{\mathcal{N}^{2(j_0+5)^2+40(\tau+n)}}{\acute{\eta}^2}\Vert X_{R_{j_0}}\Vert_{p-1, D_{j_0}\times\Pi_{\acute{\eta}}}\nonumber\\
\le &\frac{\mathcal{N}^{2(j_0+5)^2+40(\tau+n)}}{\acute{\eta}^2}\cdot \rho\left(\frac{\mathcal{N}^{2C_0(j_0+5)^2}}{\acute{\eta}^2}\rho\right)^{j_0-2}\nonumber\\
\le &\left(\frac{\mathcal{N}^{2C_0(j_0+5)^2}}{\acute{\eta}^2}\rho\right)^{j_0-1},
\end{align}
which implies the estimates \eqref{psi1} and \eqref{psi2} in the iterative lemma.
\subsection{Estimate $Z_{j_0+1}, R_{j_0+1}, Q_{j_0+1}$ and $T_{j_0+1}$}
We rewrite $Z_{j_0+1}$ as
$$Z_{j_0+1}(x,y,z,\bar{z};\xi)=\sum_{3\le j\le j_0+1}Z_{(j_0+1)j}(x,y,z,\bar{z};\xi),$$
where
$$Z_{(j_0+1)j}(x,y,z,\bar{z};\xi)=\sum_{ |\alpha|+2|\beta|+2|\mu|=j, |\mu|\le1 }Z_{j_0}^{\alpha\beta\beta\mu\mu}(x;\xi)y^{\alpha}\acute{z}^{\beta}\bar{\acute{z}}^{\beta}\hat{z}^{\mu}\bar{\hat{z}}^{\mu}.$$
Then from \eqref{zj0},\eqref{hatz}, \eqref{hatz1} and \eqref{zj0+}, one gets
$$Z_{(j_0+1)j}(x,y,z,\bar{z};\xi)=Z_{j_0j}(x,y,z,\bar{z};\xi),\quad 3\le j\le j_0,$$
and
$$Z_{(j_0+1)(j_0+1)}(x,y,z,\bar{z};\xi)=\widehat{Z_{j_0}}(x,y,z,\bar{z};\xi).$$
Thus, for $3\le j\le j_0$, by iterative assumption \eqref{xz}, we have,
\begin{align*}
|||{Z_{(j_0+1)j}}|||_{D_{j_0+1}\times \Pi_{\acute{\eta}}}\lessdot \rho^3\left(\frac{\mathcal{N}^{2C_0(j_0+4)^2}}{\acute{\eta}^2}\rho\right)^{j-3}\le \rho^3\left(\frac{\mathcal{N}^{2C_0(j_0+5)^2)}}{\acute{\eta}^2}\rho\right)^{j-3}.
\end{align*}
Similarly, by \eqref{xz}, one can obtain
$$\Vert X_{Z_{(j_0+1)j}}\Vert_{p-1,D_{j_0+1}\times \Pi_{\acute{\eta}}}\lessdot \rho\left(\frac{\mathcal{N}^{2C_0(j_0+5)^2}}{\acute{\eta}^2}\rho\right)^{j-3}, \ \ 3\le j\le j_0.$$

For $j= j_0+1$, from \eqref{hatz} and \eqref{hatz1}, we have
\begin{align*}
&|||{Z_{(j_0+1)(j_0+1)}}|||_{D_{j_0+1}\times \Pi_{\acute{\eta}}}
=|||{\widehat{Z_{j_0}}}|||_{D_{j_0+1}\times \Pi_{\acute{\eta}}}
\\
\le& |||{R_{j_0(j_0+1)}}|||_{D_{j_0+1}\times \Pi_{\acute{\eta}}}
\lessdot \rho^3\left(\frac{\mathcal{N}^{2C_0(j_0+5)^2}}{\acute{\eta}^2}\rho\right)^{j_0-2}.
\end{align*}
Similarly, one can obtain
$$\Vert X_{Z_{(j_0+1)(j_0+1)}}\Vert_{p-1,D_{j_0+1}\times \Pi_{\acute{\eta}}}\le \Vert X_{R_{j_0(j_0+1)}}\Vert_{p-1,D_{j_0+1}\times \Pi_{\acute{\eta}}}\lessdot \rho\left(\frac{\mathcal{N}^{2C_0(j_0+5)^2}}{\acute{\eta}^2}\rho\right)^{j_0-2}.$$
The estimate for $Z_{j_0+1}$ has been established until now.

\medskip

To obtain the remaining estimate, we consider the term $\int_0^1(1-t)\{\{\breve{N},F_{j_0}\},F_{j_0}\}\circ X_{F_{j_0}}^tdt$. From Taylor's formula, we have
$$\int_0^1(1-t)\{\{\breve{N},F_{j_0}\},F_{j_0}\}\circ X_{F_{j_0}}^tdt=\sum_{j\geq 2}\frac{1}{j!}\breve{N}^{(j)},$$
where $\breve{N}^{(0)}=\breve{N}$, $\breve{N}^{(j)}=\{\breve{N}^{(j-1)},F_{j_0}\}$ for $j\geq 1.$ In addition, $\breve{N}^{(j)}$ has the following form
$$\breve{N}^{(j)}=\sum_{2|\alpha|+|\beta|+|\gamma|+|\mu|+|\nu|=j(j_0+1)+2-2j}\breve{N}^{(j)\alpha\beta\gamma\mu\nu}(x;\xi)y^{\alpha}\acute{z}^{\beta}\bar{\acute{z}}^{\gamma}\hat{z}^{\mu}\bar{\hat{z}}^{\nu}.$$
By homological equation \eqref{hh}, one gets
$$\breve{N}^{(1)}=-R_{j_0(j_0+1)}+\widehat{Z_{j_0}}+\widehat{T_{j_0}}.$$
In what follows, we are going to give an estimate for
$$\widehat{T_{j_0}}=\{\breve{N}, (1-\Gamma_K)G_{j_0}\}.$$
Since
$$\{\breve{N},G_{j_0}\}=-R_{j_0(j_0+1)}+\widehat{Z_{j_0}},$$
we have
$$\{[\breve{N}],G_{j_0}\}=-R_{j_0(j_0+1)}+\widehat{Z_{j_0}}-\{\breve{N}-[\breve{N}],G_{j_0}\},$$
where $[\breve{N}]$ denotes the mean value of the function $\breve{N}$ with respect to the variable $x$ over $\mathbb{T}^n$. Therefore,
$$\{[\breve{N}], (1-\Gamma_K)G_{j_0}\}=(1-\Gamma_K)(-R_{j_0(j_0+1)}+\widehat{Z_{j_0}}-\{\breve{N}-[\breve{N}],G_{j_0}\}).$$
Let $K=\frac{\mathcal{M}}{\sigma^\prime}|\log \frac{1}{\rho}|=\mathcal{M}\mathcal{N}^2|\log \frac{1}{\rho}|$. It is easy to obtain
\begin{equation}\label{hatp11}
|||\{[\breve{N}], (1-\Gamma_K)G_{j_0}\}|||_{D^{\prime\prime}_{j_0}\times \Pi_{\acute{\eta}}}\lessdot\rho^3\left(\frac{\mathcal{N}^{2C_0(j_0+5)^2}}{\acute{\eta}^2}\rho\right)^{\mathcal{M}},
\end{equation}
and
\begin{equation}\label{hatp22}
|||\{\breve{N}-[\breve{N}], (1-\Gamma_K)G_{j_0}\}|||_{D^{\prime\prime}_{j_0}\times \Pi_{\acute{\eta}}} \lessdot\rho^3\left(\frac{\mathcal{N}^{2C_0(j_0+5)^2}}{\acute{\eta}^2}\rho\right)^{\mathcal{M}}.
\end{equation}
Combining \eqref{hatp11} and \eqref{hatp22}, one obtains
\begin{equation}\label{hung}
|||\widehat{T_{j_0}}|||_{D^{\prime\prime}_{j_0}\times \Pi_{\acute{\eta}}} =|||\{\breve{N}, (1-\Gamma_K)G_{j_0}\}|||_{D^{\prime\prime}_{j_0}\times \Pi_{\acute{\eta}}} \lessdot\rho^3\left(\frac{\mathcal{N}^{2C_0(j_0+5)^2}}{\acute{\eta}^2}\rho\right)^{\mathcal{M}}.
\end{equation}
Accordingly, one can obtain
\begin{equation}\label{hub1}
\Vert X_{\widehat{T_{j_0}}}\Vert_{p-1,D^{\prime\prime}_{j_0}\times \Pi_{\acute{\eta}}}\lessdot\rho\left(\frac{\mathcal{N}^{2C_0(j_0+5)^2}}{\acute{\eta}^2}\rho\right)^{\mathcal{M}}.
\end{equation}
From \eqref{xr}, \eqref{hatz}, \eqref{hatz1} and \eqref{hung}, we get
\begin{equation*}
|||{\breve{N}^{(1)}}|||_{D^{\prime\prime}_{j_0}\times \Pi_{\acute{\eta}}}\le |||{R_{j_0(j_0+1)}}|||_{D^{\prime\prime}_{j_0}\times \Pi_{\acute{\eta}}}+ |||\widehat{Z_{j_0}}|||_{D^{\prime\prime}_{j_0}\times \Pi_{\acute{\eta}}}+|||\widehat{T_{j_0}}|||_{D_{j_0+1}\times \Pi_{\acute{\eta}}}.
\end{equation*}
\begin{equation}\label{4.27}
\lessdot \rho^3\left(\frac{\mathcal{N}^{2C_0(j_0+5)^2)}}{\acute{\eta}^2}\rho\right)^{j_0-2}.
\end{equation}
Similarly, one can get
\begin{equation}\label{huzuui}
\Vert X_{\breve{N}^{(1)}}\Vert_{p-1,D^{\prime\prime}_{j_0}\times \Pi_{\acute{\eta}}}\lessdot \rho\left(\frac{\mathcal{N}^{2C_0(j_0+5)^2)}}{\acute{\eta}^2}\rho\right)^{j_0-2}.
\end{equation}

For $2\le j\le \mathcal{M}+2$, let
$$\varrho_j=\frac{1}{12j\mathcal{N}^2},\quad\rho_j=\frac{\rho}{12j\mathcal{M}}.
$$
Thus, by \eqref{4.27}, \eqref{fujb} and Lemma \ref{brak111}, we obtain
\begin{align}
&\frac{1}{j!}|||{\breve{N}^{(j)}}|||_{D_{j_0+1}\times\Pi_{\acute{\eta}}}\nonumber\\
\lessdot& \frac{1}{j!}\left(\frac{K^{n+1} (2\mathcal{N}+n)^{j_0+1-2}}{\varrho_j\rho_j\rho}\right)^{j-1}\cdot |||\breve{N}^{(1)}|||_{D_{j_0}\times\Pi_{\acute{\eta}}}\cdot (|||{F_{j_0}}|||_{D^{\prime\prime}_{j_0}\times\Pi_{\acute{\eta}}}^*)^{j-1}\nonumber\\
\lessdot&\frac{j^{2(j-1)}}{j!}\left(\frac{K^{n+1}(2\mathcal{N}+n)^{j_0-1}\mathcal{N}^{2}}{\rho^2}\right)^{j-1}\cdot \rho^3 \left(\frac{\mathcal{N}^{2C_0(j_0+5)^2}}{\acute{\eta}^2}\rho\right)^{j_0-2}\nonumber\\
&\cdot \left(\rho^2\Big(\frac{\mathcal{N}^{2C_0(j_0+5)^2}}{\acute{\eta}^2}\rho\Big)^{j_0-1}\right)^{j-1}\nonumber\\
\lessdot &\rho^3\left(\frac{\mathcal{N}^{2C_0(j_0+5)^2}}{\acute{\eta}^2}\rho\right)^{j_0-2}\cdot \left(K^{n+1}\mathcal{N}^2\Big(\frac{(2\mathcal{N}+n)\mathcal{N}^{2C_0(j_0+5)^2}}{\acute{\eta}^2}\rho\Big)^{j_0-1}\right)^{j-1}\nonumber\\
\lessdot & \rho^3\left(\frac{\mathcal{N}^{2C_0(j_0+5)^2}}{\acute{\eta}^2}\rho\right)^{j_0-2}\cdot \left(\frac{\mathcal{N}^{2C_0(j_0+6)^2}}{\acute{\eta}^2}\rho\right)^{(j_0-1)(j-1)}\nonumber\\
\lessdot &\rho^3\Big(\frac{\mathcal{N}^{2C_0(j_0+6)^2}}{\acute{\eta}^2}\rho\Big)^{j(j_0+1)-1-2j}
=\rho^3\Big(\frac{\mathcal{N}^{2C_0(j_0+6)^2}}{\acute{\eta}^2}\rho\Big)^{[j(j_0+1)+2-2j]-3},\label{nee1}
\end{align}
where we use $j\le \mathcal{M}+2$. Similarly,  by using Cauchy estimates and \eqref{huzuui}, \eqref{pxf}, we have
\begin{equation}\label{nee}
\frac{1}{j!}\Vert X_{\breve{N}^{(j)}}\Vert_{p-1,D_{j_0+1}\times \Pi_{\acute{\eta}}}\lessdot
\rho\Big(\frac{\mathcal{N}^{2C_0(j_0+6)^2}}{\acute{\eta}^2}\rho\Big)^{[j(j_0+1)+2-2j]-3}.
\end{equation}

Denote
$$j_*(j_0+1)+2-2j_*= \mathcal{M}+3.$$
Then, we can rewrite
$$\int_0^1(1-t)\{\{\breve{N},F_{j_0}\},F_{j_0}\}\circ X_{F_{j_0}}^tdt=\sum_{2\le j\le j_*-1}\frac{1}{j!}\breve{N}^{(j)}+\int_0^1(1-t)^{j_*}\breve{N}^{(j_*)}\circ\Psi_{j_0}dt$$
Thus,  the term $\int_0^1(1-t)\{\{\breve{N},F_{j_0}\},F_{j_0}\}\circ X_{F_{j_0}}^tdt$ is divided into two parts, the first part $\sum_{2\le j\le j_*-1}\frac{1}{j!}\breve{N}^{(j)}$ can be divided into two functions
$$R^{(1)}=\sum_{2\le j\le j_*-1}\frac{1}{j!}\sum_{\substack{2|\alpha|+|\beta|+|\gamma|+|\mu|+|\nu|=j(j_0+1)+2-2j\\
 |\mu|+|\nu|\le2} }\breve{N}^{(j)\alpha\beta\gamma\mu\nu}(x;\xi)y^{\alpha}\acute{z}^{\beta}\bar{\acute{z}}^{\gamma}\hat{z}^{\mu}\bar{\hat{z}}^{\nu},$$
 and
$$Q^{(1)}=\sum_{2\le j\le j_*-1}\frac{1}{j!}\sum_{\substack{2|\alpha|+|\beta|+|\gamma|+|\mu|+|\nu|=j(j_0+1)+2-2j\\
 |\mu|+|\nu|\geq 3} }\breve{N}^{(j)\alpha\beta\gamma\mu\nu}(x;\xi)y^{\alpha}\acute{z}^{\beta}\bar{\acute{z}}^{\gamma}\hat{z}^{\mu}\bar{\hat{z}}^{\nu},$$
where $R^{(1)}$ contains the terms has variables $(z_j,\bar{z}_j)_{|j|\geq \mathcal{N}+1}$ at most order $2$ and $Q^{(1)}$ contains the terms has variables $(z_j,\bar{z}_j)_{|j|\geq \mathcal{N}+1}$ at least order $3$. Therefore, $R^{(1)}$ and $Q^{(1)}$ can be put into $R_{j_0+1}$ and $Q_{j_0+1}$, respectively. The second part
\begin{equation}\label{p111}
P_1:=\int_0^1(1-t)^{j_*}\breve{N}^{(j_*)}\circ\Psi_{j_0}dt
\end{equation}
can be put into $T_{j_0+1}$. Moreover, from \eqref{nee}, we have the following estimate for $P_1$
\begin{equation}\label{hub2}
\Vert X_{P_1}\Vert_{p-1,D_{j_0+1}\times \Pi_{\acute{\eta}}}=\Vert X_{\breve{N}^{(j_*)}\circ \Psi_{j_0}}\Vert_{p-1,D_{j_0+1}\times \Pi_{\acute{\eta}}}\lessdot \rho\Big(\frac{\mathcal{N}^{2C_0(j_0+6)^2}}{\acute{\eta}^2}\rho\Big)^{\mathcal{M}}.
\end{equation}

Next, we will give the estimate for the term $\int_0^1\{R_{j_0(j_0+1)}+Z_{j_0},F_{j_0}\}\circ X_{F_{j_0}}^tdt$. Let
$$W_i=Z_{j_0i},\quad 3\le i\le j_0$$
and
$$W_{j_0+1}=R_{j_0(j_0+1)}.$$
Then
$$\sum_{3\le i\le j_0+1}W_i=R_{j_0(j_0+1)}+Z_{j_0}.$$
Denote
$$W_i^{(0)}=W_i\circ X_{F_{j_0}}^t|_{t=0}=W_i,$$
and
$$W_i^{(j)}=\{W_i^{(j-1)},F_{j_0}\},\quad j\geq 1.$$
By Taylor's formula, we have
$$\int_0^1\{R_{j_0(j_0+1)}+Z_{j_0},F_{j_0}\}\circ X_{F_{j_0}}^tdt=\sum_{3\le i\le j_0+1}(\sum_{j\geq 1}\frac{1}{j!}W_i^{(j)}).$$
Note that $W_i^{(j)}(x,y,z,\bar{z};\xi)$ has the following form
$$W_i^{(j)}=\sum_{2|\alpha|+|\beta|+|\gamma|+|\mu|+|\nu|=j(j_0+1)+i-2j}W_i^{(j)\alpha\beta\gamma\mu\nu}(x;\xi)y^{\alpha}\acute{z}^{\beta}\bar{\acute{z}}^{\gamma}\hat{z}^{\mu}\bar{\hat{z}}^{\nu}.$$
Using the discussion similar to ${\breve{N}^{(j)}}$ (see \eqref{nee1}, \eqref{nee}) and the inequalities \eqref{xz}, \eqref{xr}, one obtains that, for $1\le j\le \mathcal{M}+2$,
$$\frac{1}{j!}|||{W_i^{(j)}}|||_{D_{j_0+1}\times\Pi_{\acute{\eta}}}\lessdot \rho^3\left(\frac{\mathcal{N}^{2C_0(j_0+6)^2}}{\acute{\eta}^2}\rho\right)^{(j(j_0+1)+i-2j)-3},$$
and
\begin{equation}\label{rubi}
\frac{1}{j!}\Vert X_{W_i^{(j)}}\Vert_{p-1,D_{j_0+1}\times \Pi_{\acute{\eta}}} \lessdot \rho\left(\frac{\mathcal{N}^{2C_0(j_0+6)^2}}{\acute{\eta}^2}\rho\right)^{(j(j_0+1)+i-2j)-3}.
\end{equation}
In addition, we have
$$2|\alpha|+|\beta|+|\gamma|+|\mu|+|\nu|=j(j_0+1)+i-2j\geq j_0+2,$$
since $j\geq1, j_0\geq2, i\geq 3$.

For each $3\le i\le j_0+1$, we denote
$$j^i_*(j_0+1)+i-2j_*^i=\mathcal{M}+3.$$
Then, we rewrite the term
$$\int_0^1\{R_{j_0(j_0+1)}+Z_{j_0},F_{j_0}\}\circ X_{F_{j_0}}^tdt$$
$$=\sum_{3\le i\le j_0+1}\left(\sum_{1\le j\le j_*^i}\frac{1}{j!}W_i^{(j)}+\int_0^1(1-t)^{j^i_*}W_i^{(j_*^i)}\circ\Psi_{j_0}dt\right)$$
 As splitting the term $\sum_{2\le j\le j_*-1}\frac{1}{j!}\breve{N}^{(j)}$, we can write $\sum_{3\le i\le j_0+1}\sum_{1\le j\le j_*^i}\frac{1}{j!}W_i^{(j)}$ into two functions $R^{(2)}$ and $Q^{(2)}$ with
$$R^{(2)}=\sum_{3\le i\le j_0+1}\sum_{1\le j\le j_*^i}\frac{1}{j!}\sum_{\substack{2|\alpha|+|\beta|+|\gamma|+|\mu|+|\nu|=j(j_0+1)+i-2j\\
 |\mu|+|\nu|\le2} }W_i^{(j)\alpha\beta\gamma\mu\nu}(x;\xi)y^{\alpha}\acute{z}^{\beta}\bar{\acute{z}}^{\gamma}\hat{z}^{\mu}\bar{\hat{z}}^{\nu},$$
$$Q^{(2)}=\sum_{3\le i\le j_0+1}\sum_{1\le j\le j_*^i}\frac{1}{j!}\sum_{\substack{2|\alpha|+|\beta|+|\gamma|+|\mu|+|\nu|=j(j_0+1)+i-2j\\
 |\mu|+|\nu|\geq 3} }W_i^{(j)\alpha\beta\gamma\mu\nu}(x;\xi)y^{\alpha}\acute{z}^{\beta}\bar{\acute{z}}^{\gamma}\hat{z}^{\mu}\bar{\hat{z}}^{\nu},$$
where $R^{(2)}$ contains the terms has variables $(z_j,\bar{z}_j)_{|j|\geq \mathcal{N}+1}$ at most order $2$ and $Q^{(2)}$ contains the terms has variables $(z_j,\bar{z}_j)_{|j|\geq \mathcal{N}+1}$ at least order $3$. Therefore,  we can put $R^{(2)}$ and $Q^{(2)}$ into $R_{j_0+1}$ and $Q_{j_0+1}$, respectively. We denote the second term by
\begin{equation}\label{p222}
P_2:=\sum_{3\le i\le j_0+1}\int_0^1(1-t)^{j^i_*}W_i^{(j_*^i)}\circ\Psi_{j_0}dt.
\end{equation}
 Moreover, by \eqref{rubi}% and $\Psi_{j_0}:D_{j_0+1}\rightarrow D_{j_0}$
 , we have the estimate
\begin{equation}\label{hub3}
\Vert X_{P_2}\Vert_{p-1,D_{j_0+1}\times \Pi_{\acute{\eta}}}\lessdot \rho \Big(\frac{\mathcal{N}^{2C_0(j_0+6)^2}}{\acute{\eta}^2}\rho\Big)^{\mathcal{M}}.
\end{equation}

\medskip

 Finally, we give the estimate for the term
$(R_{j_0}-R_{j_0(j_0+1)}+Q_{j_0})\circ X_{F_{j_0}}^t|_{t=1}$. Let
$$U_i=R_{j_0i},\quad i\geq j_0+2,$$
and
$$V_i=Q_{j_0i},\quad i\geq 3.$$
Then
$$\sum_{i\geq j_0+2}U_i+\sum_{i\geq 3}V_i=R_{j_0}-R_{j_0(j_0+1)}+Q_{j_0}.$$
For simplicity, denote $A_i=U_i$ or $V_i$. Let
$$A_i^{(0)}=A_i\circ X_{F_{j_0}}^t|_{t=0}=A_i,$$
and
$$A_i^{(j)}=\{A_i^{(j-1)},F_{j_0}\},\quad j\geq 1.$$
Using Taylor's formula, we have
\begin{align*}
&(R_{j_0}-R_{j_0(j_0+1)}+Q_{j_0})\circ X_{F_{j_0}}^t|_{t=1}\\
=&R_{j_0}-R_{j_0(j_0+1)}+Q_{j_0}+\sum_{i\geq j_0+2}(\sum_{j\geq 1}\frac{1}{j!}U_i^{(j)})+\sum_{i\geq 3}(\sum_{j\geq 1}V_i^{(j)}).
\end{align*}

%$$\sum_{i\geq j_0+2}(\sum_{1\le j\le \mathcal{M}+2}\frac{1}{j!}U_i^{(j)}+\int_0^1\frac{1}{(\mathcal{M}+3)!}(1-t)^{\mathcal{M}+3}U_i^{(\mathcal{M}+3)}\circ X_{F_{j_0}}^tdt)$$

For $j\geq 1$, $A_i^{(j)}(x,y,z,\bar{z};\xi)$ with $A_i=U_i $ or $V_i$ has the following form
$$A_i^{(j)}=\sum_{2|\alpha|+|\beta|+|\gamma|+|\mu|+|\nu|=j(j_0+1)+i-2j}A_i^{(j)\alpha\beta\gamma\mu\nu}(x;\xi)y^{\alpha}\acute{z}^{\beta}\bar{\acute{z}}^{\gamma}\hat{z}^{\mu}\bar{\hat{z}}^{\nu}.$$
Using the discussion similar to ${\breve{N}^{(j)}}$ (see \eqref{nee1}, \eqref{nee}) and the inequalities \eqref{xr}, \eqref{xy}, one obtains, for $j\le \mathcal{M}+2$,
$$\frac{1}{j!}|||{A_i^{(j)}}|||_{D_{j_0+1}\times\Pi_{\acute{\eta}}}\lessdot \rho^3\left(\frac{\mathcal{N}^{2C_0(j_0+6)^2}}{\acute{\eta}^2}\rho\right)^{(j(j_0+1)+i-2j)-3},$$
and
\begin{equation}\label{wohen}
\frac{1}{j!}\Vert X_{A_i^{(j)}}\Vert_{p-1,D_{j_0+1}\times \Pi_{\acute{\eta}}}\lessdot \rho\left(\frac{\mathcal{N}^{2C_0(j_0+6)^2}}{\acute{\eta}^2}\rho\right)^{(j(j_0+1)+i-2j)-3}.
\end{equation}
In addition, we have
$$2|\alpha|+|\beta|+|\gamma|+|\mu|+|\nu|=j(j_0+1)+i-2j\geq j_0+2,$$
since $j\geq1, j_0\geq2, i\geq 3$.

Let
$$j_*^i(j_0+1)+i-2j_*^i=\mathcal{M}+3.$$
Then we rewrite
\begin{equation*}
  \begin{split}
     &(R_{j_0}-R_{j_0(j_0+1)}+Q_{j_0})\circ X_{F_{j_0}}^t|_{t=1}\\
     =&R_{j_0}-R_{j_0(j_0+1)}+Q_{j_0} \\
     &+\sum_{i\geq j_0+2}\left(\sum_{1\le j\le j_*^i}\frac{1}{j!}U_i^{(j)}+\int_0^1(1-t)^{j^i_*}U_i^{(j_*^i)}\circ\Psi_{j_0}dt\right)\\
     &+\sum_{i\geq 3}\left(\sum_{1\le j\le j_*^i}\frac{1}{j!}V_i^{(j)}+\int_0^1(1-t)^{j^i_*}V_i^{(j_*^i)}\circ\Psi_{j_0}dt\right).
   \end{split}
\end{equation*}

We rewrite
$$\sum_{i\geq j_0+2}\sum_{1\le j\le j_*^i}\frac{1}{j!}U_i^{(j)}=R^{(3)}+Q^{(3)},$$
$$\sum_{i\geq 3}\sum_{1\le j\le j_*^i}\frac{1}{j!}V_i^{(j)}=R^{(4)}+Q^{(4)}$$
with
$$R^{(3)}=\sum_{i\geq j_0+2}\sum_{1\le j\le j_*^i}\frac{1}{j!}\sum_{\substack{2|\alpha|+|\beta|+|\gamma|+|\mu|+|\nu|=j(j_0+1)+i-2j\\
 |\mu|+|\nu|\le2} }U_i^{(j)\alpha\beta\gamma\mu\nu}(x;\xi)y^{\alpha}\acute{z}^{\beta}\bar{\acute{z}}^{\gamma}\hat{z}^{\mu}\bar{\hat{z}}^{\nu},$$
$$Q^{(3)}=\sum_{i\geq j_0+2}\sum_{1\le j\le j_*^i}\frac{1}{j!}\sum_{\substack{2|\alpha|+|\beta|+|\gamma|+|\mu|+|\nu|=j(j_0+1)+i-2j\\
 |\mu|+|\nu|\geq 3} }U_i^{(j)\alpha\beta\gamma\mu\nu}(x;\xi)y^{\alpha}\acute{z}^{\beta}\bar{\acute{z}}^{\gamma}\hat{z}^{\mu}\bar{\hat{z}}^{\nu},$$
 $$R^{(4)}=\sum_{i\geq 3}\sum_{1\le j\le j_*^i}\frac{1}{j!}\sum_{\substack{2|\alpha|+|\beta|+|\gamma|+|\mu|+|\nu|=j(j_0+1)+i-2j\\
 |\mu|+|\nu|\le2} }V_i^{(j)\alpha\beta\gamma\mu\nu}(x;\xi)y^{\alpha}\acute{z}^{\beta}\bar{\acute{z}}^{\gamma}\hat{z}^{\mu}\bar{\hat{z}}^{\nu},$$
$$Q^{(4)}=\sum_{i\geq 3}\sum_{1\le j\le j_*^i}\frac{1}{j!}\sum_{\substack{2|\alpha|+|\beta|+|\gamma|+|\mu|+|\nu|=j(j_0+1)+i-2j\\
 |\mu|+|\nu|\geq 3} }V_i^{(j)\alpha\beta\gamma\mu\nu}(x;\xi)y^{\alpha}\acute{z}^{\beta}\bar{\acute{z}}^{\gamma}\hat{z}^{\mu}\bar{\hat{z}}^{\nu}.$$
Therefore, we can put the terms $R_{j_0}-R_{j_0(j_0+1)}+R^{(3)}+R^{(4)}$ and $Q_{j_0}+Q^{(3)}+Q^{(4)}$ into $R_{j_0+1}$ and $Q_{j_0+1}$, respectively. Denote
\begin{equation}\label{p333}
P_3=\sum_{i\geq j_0+2}\int_0^1(1-t)^{j^i_*}U_i^{(j_*^i)}\circ\Psi_{j_0}dt+\sum_{i\geq 3}\int_0^1(1-t)^{j^i_*}V_i^{(j_*^i)}\circ\Psi_{j_0}dt.
\end{equation}
By \eqref{wohen}, we have the estimate
\begin{equation}\label{hub4}
\Vert X_{P_3}\Vert_{p-1,D_{j_0+1}\times \Pi_{\acute{\eta}}}\lessdot \rho \Big(\frac{\mathcal{N}^{2C_0(j_0+6)^2}}{\acute{\eta}^2}\rho\Big)^{\mathcal{M}}.
\end{equation}

As a consequence, we obtain that
$$R_{j_0+1}=R_{j_0}-R_{j_0(j_0+1)}+R^{(1)}+R^{(2)}+R^{(3)}+R^{(4)},$$
$$Q_{j_0+1}=Q_{j_0}+Q^{(1)}+Q^{(2)}+Q^{(3)}+Q^{(4)}.$$
In addition, from \eqref{nee}, \eqref{rubi} and \eqref{wohen}, we have for $j\le \mathcal{M}+2$,
\begin{align*}
&\Vert X_{R_{(j_0+1)j}}\Vert_{p-1,D_{j_0+1}\times \Pi_{\acute{\eta}}},,\ \frac{1}{\rho^2}|||{R_{(j_0+1)j}}|||_{D_{j_0+1}\times \Pi}\lessdot \rho\left(\frac{\mathcal{N}^{2C_0(j_0+6)^2}}{\acute{\eta}^2}\rho\right)^{j-3},\label{xr}\\
 &\Vert X_{Q_{(j_0+1)j}}\Vert_{p-1,D_{j_0+1}\times \Pi_{\acute{\eta}}}, \frac{1}{\rho^2} ||||{Q_{(j_0+1)j}}|||_{D_{j_0+1}\times \Pi_{\acute{\eta}}}\lessdot \rho\left(\frac{\mathcal{N}^{2C_0(j_0+6)^2}}{\acute{\eta}^2}\rho\right)^{j-3}.
 \end{align*}

 We rewrite $T_{j_0+1}$ as
$$T_{j_0+1}(x,y,z,\bar{z};\xi)=\widehat{T_{j_0}}+T_{j_0} \circ \Psi_{j_0} +P_1+P_2+P_3,$$
where $P_{1},P_{2},P_{3}$ are defined by \eqref{p111}, \eqref{p222} and \eqref{p333}.
Thus, noting \eqref{hub1}, \eqref{hub2}, \eqref{hub3} and \eqref{hub4}, we finally obtain the estimates for $T_{j_0+1}$
$$\Vert X_{T_{j_0+1}}\Vert_{p-1,D_{j_0+1}\times \Pi_{\acute{\eta}}}\lessdot  \rho \Big(\frac{\mathcal{N}^{2C_0(j_0+6)^2}}{\acute{\eta}^2}\rho\Big)^{\mathcal{M}}.$$

 Thus, we finish the proof of Lemma \ref{iterm}.

\subsection{Proof of Theorem \ref{parth}}
By Theorem \ref{normaltheorem}, we obtain a normal form of order $2$ around KAM tori, which reads
$$\breve{H}(x,y,z,\bar{z};\xi)=\breve{N}(x,y,z,\bar{z};\xi)+\breve{P}(x,y,z,\bar{z};\xi),$$
where
$$\breve{N}(x,y,z,\bar{z};\xi)=\sum_{j=1}^n\breve{\omega}_j(\xi)y_j+\sum_{j\geq 1}\frac{\breve{\Omega}_j(x;\xi)}{j}z_j\bar{z}_j,$$
%$$=\sum_{j=1}^n\breve{\omega}_j(\xi)y_j+\sum_{|j|\le \mathcal{N}}\breve{\Omega}_j(\xi)z_j\bar{z}_j+\sum_{|j|\geq\mathcal{N}+ 1}\breve{\Omega}_j(x;\xi)z_j\bar{z}_j$$
and
$$\breve{P}(x,y,z,\bar{z};\xi)=\sum_{\substack{\alpha\in\mathbb{N}^n,\beta,\gamma\in \mathbb{N}^{\mathbb{N}},\\ 2|\alpha|+|\beta+\gamma|\geq 3}}\breve{P}^{\alpha\beta\gamma}(x;\xi)y^{\alpha}z^{\beta}\bar{z}^{\gamma}.$$

In the following, we use the notations in Lemma \ref{iterm}. Let
$$H_2(x,y,z,\bar{z};\xi)=\breve{H}(x,y,z,\bar{z};\xi),$$
and write $\breve{P}$
\begin{equation}\label{huaddp}
\breve{P}(x,y,z,\bar{z};\xi)=Z_{2}(x,y,z,\bar{z};\xi)+R_{2}(x,y,z,\bar{z};\xi)+Q_{2}(x,y,z,\bar{z};\xi)+T_{2}(x,y,z,\bar{z};\xi),
\end{equation}
where
\begin{align}
&Z_{2}(x,y,z,\bar{z};\xi)=\sum_{3\le j\le 2}Z_{2j}(x,y,z,\bar{z};\xi)=0,\nonumber\\
%&A_{2}(x,y,z,\bar{z};\xi)=\sum_{3\le j\le 2}A_{2j}(x,y,z,\bar{z};\xi)=0,\\
&R_{2}(x,y,z,\bar{z};\xi)=\sum_{3\le j\le \mathcal{M}+2}R_{2j}(x,y,z,\bar{z};\xi),\nonumber\\
&Q_{2}(x,y,z,\bar{z};\xi)=\sum_{3\le j\le \mathcal{M}+2}Q_{2j}(x,y,z,\bar{z};\xi),\nonumber\\
&T_{2}(x,y,z,\bar{z};\xi)=\sum_{j> \mathcal{M}+2}T_{2j}(x,y,z,\bar{z};\xi),\label{addt}
\end{align}
with
$$Z_{2j}(x,y,z,\bar{z};\xi)=\sum_{ 2|\alpha|+2|\beta|+2|\mu|=j, |\mu|\le1 }Z_{2}^{\alpha\beta\beta\mu\mu}(x;\xi)y^{\alpha}\acute{z}^{\beta}\bar{\acute{z}}^{\beta}\hat{z}^{\mu}\bar{\hat{z}}^{\mu},$$
$$R_{2j}(x,y,z,\bar{z};\xi)=\sum_{2|\alpha|+|\beta|+|\gamma|+|\mu|+|\nu|=j, |\mu|+|\nu|\le2 }R_{2}^{\alpha\beta\gamma\mu\nu}(x;\xi)y^{\alpha}\acute{z}^{\beta}\bar{\acute{z}}^{\gamma}\hat{z}^{\mu}\bar{\hat{z}}^{\nu},$$
$$Q_{2j}(x,y,z,\bar{z};\xi)=\sum_{2|\alpha|+|\beta|+|\gamma|+|\mu|+|\nu|=j,|\mu|+|\nu|\geq3 }Q_2^{\alpha\beta\gamma\mu\nu}(x;\xi)y^{\alpha}\acute{z}^{\beta}\bar{\acute{z}}^{\gamma}\hat{z}^{\mu}\bar{\hat{z}}^{\nu},$$
and
$$T_{2j}(x,y,z,\bar{z};\xi)=\sum_{2|\alpha|+|\beta|+|\gamma|+|\mu|+|\nu|=j }T_{2}^{\alpha\beta\gamma\mu\nu}(x;\xi)y^{\alpha}\acute{z}^{\beta}\bar{\acute{z}}^{\gamma}\hat{z}^{\mu}\bar{\hat{z}}^{\nu}.$$
Denote $W_j(x,y,z,\bar{z}:\xi)=R_{2j}(x,y,z,\bar{z}:\xi)$, $Q_{2j}(x,y,z,\bar{z}:\xi)$ or  $T_{2j}(x,y,z,\bar{z}:\xi)$ for $j\geq 3$.
Then, by the definition of the norm $|||\cdot |||$ in definition \ref{defnd} and \eqref{brevep}, one gets
$$|||{W_j}|||_{D(s_0/2, 5\rho,5\rho)\times\Pi_\eta}\le |||{\breve{P}}|||_{D(s_0/2, r_0/2,r_0/2)\times\Pi_\eta}\left(\frac{10\rho}{r_0}\right)^{j}\le \epsilon(1+c\eta^6 \varepsilon)\left(\frac{10\rho}{r_0}\right)^{j}.$$
By the choice of $\mathcal{N}$ (noting that larg $\mathcal{N}$ is depending on $s_0,r_0,n,\mathcal{M},\tau$ and $p$) and $0<\breve{\eta}<1$, we have
$$|||{W_j}|||_{D(s_0/2, 5\rho,5\rho)\times\Pi_\eta}\le \rho^3\left(\frac{\mathcal{N}}{\acute{\eta}^2}\rho\right)^{j-3},\quad j\geq 3.$$
From Lemma \ref{guji}, it implies that
\begin{align}
\Vert X_{W_j}\Vert_{p-1, D(5\varrho, 5\rho,5\rho)\times\Pi_\eta}\le&  c^{j-2} \frac{j^{p+2}}{(\frac{s_0}{2}-5\varrho)^p} |||W_j|||_{D(s_0/2, 5\rho,5\rho)\times \Pi_\eta}\rho^{-2}\nonumber\\
\le &c\rho\left(\frac{\mathcal{N}}{\acute{\eta}^2}\rho\right)^{j-3},\quad j\geq 3.\label{addwe}
\end{align}
By \eqref{huaddp}, \eqref{addt} and \eqref{addwe}, we have
\begin{align*}
\Vert X_{T_2}\Vert_{p-1,D(5\varrho, 5\rho,5\rho)\times\Pi_\eta}\le& \sum_{j\geq \mathcal{M}+3}  \Vert X_{T_{2j}}\Vert_{p-1,D(5\varrho, 5\rho,5\rho)\times\Pi_\eta}\\
\le&c\sum_{j\geq \mathcal{M}+3}  \rho \left(\frac{\mathcal{N}}{\acute{\eta}^2}\rho\right)^{j-3}\le\rho \left(\frac{\mathcal{N}^2}{\acute{\eta}^2}\rho\right)^{\mathcal{M}} .
\end{align*}

Let $\Psi=\Psi_2\circ\cdots\circ\Psi_{\mathcal{M}+1}$.
Then by Lemma \ref{iterm}, one gets
$$\breve{\breve{H}}=\breve{H}\circ\Phi=\breve{N}(x,y,z,\bar{z};\xi)+Z_{\mathcal{M}+2}(x,y,z,\bar{z};\xi)$$
$$+R_{\mathcal{M}+2}(x,y,z,\bar{z};\xi)+Q_{\mathcal{M}+2}(x,y,z,\bar{z};\xi)+T_{\mathcal{M}+2}(x,y,z,\bar{z};\xi),$$
where
\begin{align*}
&Z_{\mathcal{M}+2}(y,z,\bar{z};\xi)=\sum_{3\le j\le \mathcal{M}+2}Z_{(\mathcal{M}+2)j}(x,y,z,\bar{z};\xi),\\
%&A_{\mathcal{M}+2}(x,y,z,\bar{z};\xi)=\sum_{3\le j\le \mathcal{M}+2}A_{(\mathcal{M}+2)j}(x,y,z,\bar{z};\xi),\\
%&R_{\mathcal{M}+2}(x,y,z,\bar{z};\xi)=\sum_{j\geq \mathcal{M}+3}R_{(\mathcal{M}+2)j}(x,y,z,\bar{z};\xi),\\
&Q_{\mathcal{M}+2}(x,y,z,\bar{z};\xi)=\sum_{3\le j\le \mathcal{M}+2}Q_{(\mathcal{M}+2)j}(x,y,z,\bar{z};\xi),\\
&T_{\mathcal{M}+2}(x,y,z,\bar{z};\xi)=\sum_{j\geq 3}T_{(\mathcal{M}+2)j}(x,y,z,\bar{z};\xi),\\
&R_{\mathcal{M}+2}(x,y,z,\bar{z};\xi)=0,
\end{align*}
with
\begin{align*}
Z_{(\mathcal{M}+2)j}(y,z,\bar{z};\xi)=&\sum_{ |\alpha|+2|\beta|+2|\mu|=j, |\mu|=0 }Z_{\mathcal{M}+2}^{\alpha\beta\beta\mu\mu}(\xi)y^{\alpha}\acute{z}^{\beta}\bar{\acute{z}}^{\beta}\hat{z}^{\mu}\bar{\hat{z}}^{\mu}\\
&+\sum_{ |\alpha|+2|\beta|+2|\mu|=j, |\mu|=1 }Z_{\mathcal{M}+2}^{\alpha\beta\beta\mu\mu}(x;\xi)y^{\alpha}\acute{z}^{\beta}\bar{\acute{z}}^{\beta}\hat{z}^{\mu}\bar{\hat{z}}^{\mu},
\end{align*}
%$$R_{(\mathcal{M}+2)j}(x,y,z,\bar{z};\xi)=\sum_{|\alpha|+|\beta|+|\gamma|+|\mu|+|\nu|=j, |\mu|+|\nu|\le2 }R_{\mathcal{M}+2}^{\alpha\beta\gamma\mu\nu}(x;\xi)y^{\alpha}\acute{z}^{\beta}\bar{\acute{z}}^{\gamma}\hat{z}^{\mu}\bar{\hat{z}}^{\nu},$$
$$Q_{(\mathcal{M}+2)j}(x,y,z,\bar{z};\xi)=\sum_{|\alpha|+|\beta|+|\gamma|+|\mu|+|\nu|=j,|\mu|+|\nu|\geq3 }Q_{\mathcal{M}+2}^{\alpha\beta\gamma\mu\nu}(x;\xi)y^{\alpha}\acute{z}^{\beta}\bar{\acute{z}}^{\gamma}\hat{z}^{\mu}\bar{\hat{z}}^{\nu}.$$
$$T_{(\mathcal{M}+2)j}(x,y,z,\bar{z};\xi)=\sum_{2|\alpha|+|\beta|+|\gamma|+|\mu|+|\nu|=j }T_{j_0}^{\alpha\beta\gamma\mu\nu}(x;\xi)y^{\alpha}\acute{z}^{\beta}\bar{\acute{z}}^{\gamma}\hat{z}^{\mu}\bar{\hat{z}}^{\nu},\quad j\geq \mathcal{M}+3,$$
and
$$T_{(\mathcal{M}+2)j}(x,y,z,\bar{z};\xi)=\sum_{\substack{2|\alpha|+|\beta|+|\gamma|+|\mu|+|\nu|=j \\\beta\neq \gamma, \mu\neq\nu, |\mu|+|\nu|\le 2}}T_{j_0}^{\alpha\beta\gamma\mu\nu}(x;\xi)y^{\alpha}\acute{z}^{\beta}\bar{\acute{z}}^{\gamma}\hat{z}^{\mu}\bar{\hat{z}}^{\nu},\quad 3\le j\le \mathcal{M}+2.$$
Moreover, the functions $Z_{(\mathcal{M}+2)j}(x,y,z,\bar{z};\xi)$, $Q_{(\mathcal{M}+2)j}(x,y,z,\bar{z};\xi)$ and $T_{(\mathcal{M}+2)}(x,y,z,\bar{z};\xi)$ satisfy the following estimates,
\begin{align*}
&\Vert X_{{Z_{(\mathcal{M}+2)j}}}\Vert_{p-1, D(4\varrho,4\rho,4\rho)\times \Pi_{\acute{\eta}}},\ \frac{1}{\rho^2}|||{Z_{(\mathcal{M}+2)j}}|||_{D(4\varrho,4\rho,4\rho)\times \Pi_{\acute{\eta}}}\lessdot \rho\left(\frac{\mathcal{N}^{2C_0(\mathcal{M}+6)^2}}{\acute{\eta}^2}\rho\right)^{j-3},\\
 &\Vert X_{Q_{(\mathcal{M}+2)j}}\Vert_{p-1,D(4\varrho,4\rho,4\rho)\times \Pi_{\acute{\eta}}},\ \frac{1}{\rho^2}|||{Q_{(\mathcal{M}+2)j}}|||_{D(4\varrho,4\rho,4\rho)\times \Pi_{\acute{\eta}}}\lessdot \rho\left(\frac{\mathcal{N}^{2C_0(\mathcal{M}+7)^2}}{\acute{\eta}^2}\rho\right)^{j-3},\\
 &\Vert X_{T_{(\mathcal{M}+2)}}\Vert_{p-1,D(4\varrho,4\rho,4\rho)\times \Pi_{\acute{\eta}}}\lessdot \rho\left(\frac{\mathcal{N}^{2C_0(\mathcal{M}+7)^2}}{\acute{\eta}^2}\rho\right)^{\mathcal{M}}.
 \end{align*}
By the choice of $\mathcal{N}$, we have
$$\frac{\mathcal{N}^{2C_0(\mathcal{M}+7)^2}}{\acute{\eta}^2}\rho<\frac{1}{2}.$$
Then by the definition of $\Psi$, we have
$$\Vert \Psi-id\Vert_{p,D(4\varrho,4\rho,4\rho)}\le c\frac{\rho \mathcal{N}^C}{\acute{\eta}^2},$$
$$|| D\Psi-Id||_{p,D(4\varrho,4\rho,4\rho)}\le c\frac{\mathcal{N}^C}{\acute{\eta}^2},$$
where $C$ is some constant depending on $n$ and $\tau$.
%Finally, we compute
%\begin{align*}
%&|||{Z_{\mathcal{M}+2}}|||\le \sup_{3\le j\le \mathcal{M}+2}|||{Z_{(\mathcal{M}+2)j}}|||\\
%=&\sup_{4\le j\le \mathcal{M}+2}|||{Z_{(\mathcal{M}+2)j}}|||
%&\le \sum_{4\le j\le \mathcal{M}+2}\rho(\frac{\mathcal{N}^{2(\mathcal{M}+6)^2}\rho}{\acute{\eta}^2})^{j-3}
%\le\rho\left(\frac{\mathcal{N}^{2(\mathcal{M}+6)^2}\rho}{\acute{\eta}^2}\right),
%\end{align*}
%\begin{align*}
%|||{P_{\mathcal{M}+2}}|||%\le &\sum_{3\le j\le \mathcal{M}+2}|||{P_{(\mathcal{M}+2)j}}|||
%\le \sup_{ j\geq \mathcal{M}+3}|||{P_{(\mathcal{M}+2)j}}|||
%&\le \sum_{j\geq \mathcal{M}+3}\rho(\frac{\mathcal{N}^{2(\mathcal{M}+6)^2}\rho}{\acute{\eta}^2})^{j-3}
%\le\rho\left(\frac{\mathcal{N}^{2(\mathcal{M}+6)^2}\rho}{\acute{\eta}^2}\right)^{\mathcal{M}},
%\end{align*}
%and
%\begin{align*}
%|||X_{Q_{\mathcal{M}+2}}|||\le &\sum_{3\le j\le \mathcal{M}+2}|||X_{Q_{(\mathcal{M}+2)j}}|||\le \sum_{j\geq 3}|||X_{Q_{(\mathcal{M}+2)j}}|||\\
%&\le \sum_{j\geq 3}\rho(\frac{\mathcal{N}^{2(\mathcal{M}+6)^2}\rho}{\acute{\eta}^2})^{j-3}\le\rho.
%\end{align*}
\subsection{Measure estimate}\label{Mease}
In this section, we will show
$${\rm Meas}\ \Pi_{\acute{\eta}}\geq ({\rm Meas}\ \Pi_{\eta})(1-c\acute{\eta}),$$
where $c$ is a constant depending on $n$.

Firstly, we want to compute the number of the non-empty non-resonant sets $\mathcal{R}_{k\acute{l}\hat{l}}$ defined by \eqref{rlll}. For this purpose, we consider the term $\langle \hat{l},[\hat{\Omega}]\rangle$ which is classified into the following three cases:

{\bf Case 1}: $|\hat{l}|=1$ and $|\langle \hat{l},[\hat{\Omega}]\rangle|=|[\breve{\Omega}_i]|$, for some $i\geq \mathcal{N}+1.$
Then, it follows that
$$|\langle \hat{l},[\hat{\Omega}]\rangle|\geq \frac{c_1}{2}i^2\geq \frac{c_1}{2}i.$$

{\bf Case 2}: $|\hat{l}|=2$ and $|\langle \hat{l},[\hat{\Omega}]\rangle|=|[\breve{\Omega}_i]+[\breve{\Omega}_j]|$ or $|\langle \hat{l},[\hat{\Omega}]\rangle|=|[\breve{\Omega}_i]-[\breve{\Omega}_j]|$, ($ i\neq \pm j$).
Then
$$|\langle \hat{l},[\hat{\Omega}]\rangle|\geq c_1(i^2\pm j^2)-(i+j)c\eta^6 \epsilon\geq \frac{c_1}{2}|i+j||i-j|\geq \frac{c_1}{2}\max\{|i|,|j|\}.$$
In case 1 and case 2, assume $|\Omega_j(\xi)|\le c_2|j|^2$. Thus, if
$$\max\{|i|,|j|\}\geq \frac{8}{c_1}\left(|k||\breve{\omega}(\xi)|+c_2(\mathcal{M}+2)\mathcal{N}^2\right),$$
 then
\begin{align*}
&|\langle k,\breve{\omega}\rangle+\langle \acute{l},[\acute{\Omega}]\rangle+\langle \hat{l},[\hat{\Omega}]\rangle|\geq |\langle \hat{l},[\hat{\Omega}]\rangle|-|\langle k,\breve{\omega}\rangle+\langle \acute{l},[\acute{\Omega}]\rangle|\\
%\geq &\frac{c_1}{4}\max\{|i|,|j|\}+\frac{c_1}{4}\max\{|i|,|j|\}\\
\geq &\frac{c_1}{4}\max\{|i|,|j|\}+ 2(|k||\breve{\omega}(\xi)|+c_2(\mathcal{M}+2)\mathcal{N}^2+1)\\
&-(1+c\eta^6\epsilon)(|k||\breve{\omega}(\xi)|+c_2|\acute{l}|\mathcal{N}^2)\geq \frac{c_1}{4}\max\{|i|,|j|\},
\end{align*}
which implies $\mathcal{R}_{k\acute{l}\hat{l}}=\emptyset$ for $\max\{|i|,|j|\}\geq \frac{8}{c_1}\left(|k||\breve{\omega}(\xi)|+c_2(\mathcal{M}+2)\mathcal{N}^2\right)$. Therefore, for the above two cases, we only need to consider
$$|i|,|j| \le j_*:=\frac{8}{c_1}\left(|k||\omega(\xi)|+c_2(\mathcal{M}+2)\mathcal{N}^2\right).$$

 It means that, for the above two cases, the number of non empty sets $\mathcal{R}_{k\acute{l}\hat{l}}$ is less than
$$A_1:=(2\mathcal{N}+1)^{|\acute{l}|}\left(\frac{8}{c_1}(|k||\omega(\xi)|+c_2(\mathcal{M}+2)\mathcal{N}^2)\right)^2.$$

{\bf Case 3}: $|\hat{l}|=2$, $|\langle \hat{l},[\hat{\Omega}]\rangle|=|[\breve{\Omega}_i]-[\breve{\Omega}_j]|$ for some $|i|,|j|\geq \mathcal{N}+1$ and $i= -j.$ In this case, by momentum conservation, for the set $\mathcal{R}_{k\acute{l}\hat{l}}$ we have the relation
\begin{equation}\label{rmom}
-\sum_{b=1}^nj_bk_b-2j+\sum_{|i|\le \mathcal{N}} i(\beta_i-\gamma_i)=0,
\end{equation}
where $\acute{l}=(\beta,\gamma)$. It implies that $j$ is bounded by
$$|j|\le \frac{1}{2} |k|\max_{b}\{j_b\}+\frac{1}{2}(\mathcal{M}+2)\mathcal{N}\le j_{**}:=C|k|+\frac{1}{2}(\mathcal{M}+2)\mathcal{N}.$$

 It means that for cases 3 the number of non empty sets $\mathcal{R}_{k\acute{l}\hat{l}}$ is less than
$$A_2:=(2\mathcal{N}+1)^{|\acute{l}|}(C|k|+\frac{1}{2}(\mathcal{M}+2)\mathcal{N})^2.$$

\begin{Rem}From equation \eqref{rmom}, it is easy to observe that if $|k|=|\acute{l}|=0$, we must have $j=0$. That is, there is no such set $\mathcal{R}_{k\acute{l}\hat{l}}$ with $|k|=|\acute{l}|=0$ and $i=-j$.
\end{Rem}

In the following, we are going to estimate the measure of $\mathcal{R}_{k\acute{l}\hat{l}}$.

 For $k\neq 0$, $|i|,|j|<\max\{j_*,j_{**}\}$, denote
$$j_{max}=\max\{j: \acute{l}_j\, {\rm or}\, \hat{l}_j \neq 0\},\quad j_{min}=\min\{j: \acute{l}_j\, {\rm or}\, \hat{l}_j \neq 0\}
$$
\begin{equation}\label{JJ}
J=\left\{
\begin{array}{l}
j_{max},\quad {\rm if}\ j_{max}>-j_{min},\\
j_{min},\quad {\rm if}\ j_{max}\le-j_{min}.\\
\end{array}
\right.
\end{equation}
Without loss of generality, we assume $\acute{l}_J\neq 0$.

For $|k|<|J|$, in view of \eqref{pgyy}, \eqref{iii} and $2\mathcal{M}\epsilon \eta^6 <\frac{1}{2}$, we obtain
\begin{align*}
&|\partial_{\xi_J}(\langle k,\breve{\omega}\rangle+\langle \acute{l},[\acute{\Omega}]\rangle+\langle \hat{l},[\hat{\Omega}]\rangle)|\\
\geq &|\acute {l}_J||\partial_{\xi_{J}}[\acute{\Omega}]_J|-|\partial_{\xi_{J}}(\langle k,\breve{\omega}\rangle+\langle \acute{l},[\acute{\Omega}]\rangle+\langle \hat{l},[\hat{\Omega}]\rangle-\acute{l}_J[\acute{\Omega}]_J)|\\
\geq&|\acute {l}_J|(1-c\eta^6\epsilon)-|k|\frac{c\eta\epsilon}{|J|}-\sum_{j\neq J}|\acute{l}_j|\frac{|j|c\eta^6\epsilon}{|J|}-\sum_{j\geq \mathcal{N}, |\hat{l}|\le 2}|\hat{l}_j|\frac{|j|c\eta^6\epsilon}{|J|}\\
\geq& |\acute {l}_J|-|k|\frac{c\eta^6\epsilon}{|J|}-(\mathcal{M}+2)c\eta^6\epsilon \\
\geq & |\acute {l}_J|-[c\eta^6\epsilon+(\mathcal{M}+2)c\eta^6\epsilon]\\
\geq &\frac{1}{4}.
\end{align*}
%Since $|j|\xi_j=\xi^\prime_j$, for some function $f(\xi)$, we have $\partial_{\xi^\prime_j}f(\xi)=\frac{1}{|j|}\partial_{\xi_j}f(\xi)$.
%Therefore, one has
%\begin{align*}
%|\partial_{\xi^\prime_{p}}(\langle k,\breve{\omega}\rangle+\langle \acute{l},[\acute{\Omega}]\rangle+\langle \hat{l},[\hat{\Omega}]\rangle)|\geq &\frac{1}{4 |p|}\geq \frac{1}{4 M_{\acute{l},\hat{l}}}.
%\end{align*}
Hence
$${\rm Meas}\, \mathcal{R}_{k\acute{l}\hat{l}}\le \frac{4\acute{\eta}M_{\acute{l},\hat{l}}}{4^{\mathcal{M}}(|k|+1)^{\tau}C(\mathcal{N},\acute{l})}{\rm Meas}\,\Pi_\eta.$$

For $|k|\geq |J|$, let
$$|k_1|=\max_{1\le i\le n}\{|k_1|,\cdots,|k_n|\}.$$
Then, by $2n\mathcal{M}\epsilon \eta^6 <\frac{1}{2}$, we have
\begin{align*}
&|\partial_{\xi_1}(\langle k,\breve{\omega}\rangle+\langle \acute{l},[\acute{\Omega}]\rangle+\langle \hat{l},[\hat{\Omega}]\rangle)|\\
\geq &|k_1||\partial_{\xi_1}\breve{\omega}_1|-|\partial_{\xi_1}(\sum_{i=2}^nk_i,\breve{\omega}_i+\langle \acute{l},[\acute{\Omega}]\rangle+\langle \hat{l},[\hat{\Omega}]\rangle)|\\
\geq& |k_1|(1-c\eta^6\epsilon)-\sum_{i=2}^n|k_i|c\eta^6\epsilon-(|\acute{l}|+|\hat{l}|)c\eta^6\epsilon\cdot |J|\\
\geq& |k_1|-|k|c\eta^6\epsilon-(\mathcal{M}+2)c\eta^6\epsilon |J|\\
\geq & |k_1|-n|k_1|c\eta^6\epsilon-(\mathcal{M}+2)c\eta^6\epsilon |k|\quad (|J|\le |k|,\quad  |k|\le n|k_1|)\\
\geq &\frac{1}{4}|k_1|\\
\geq& \frac{1}{4}.
\end{align*}
%Similarly, we have
%\begin{align*}
%&|\partial_{\xi^\prime_1}(\langle k,\breve{\omega}\rangle+\langle \acute{l},[\acute{\Omega}]\rangle+\langle \hat{l},[\hat{\Omega}]\rangle)|\\
%\geq& \frac{1}{4|j_1|}\geq \frac{1}{4C_J},
%\end{align*}
%where $C_{J_n}$ is a constant depending on the index set $J$.
Hence
$${\rm Meas}\, \mathcal{R}_{k\acute{l}\hat{l}}\le \frac{4\ \acute{\eta}M_{\acute{l},\hat{l}}}{4^{\mathcal{M}}(|k|+1)^{\tau}C(\mathcal{N},\acute{l})}{\rm Meas}\,\Pi_\eta.$$

If $|k|=0$ and $|i|,|j|<\max\{j_*,j_{**}\}$, we let $J\in\mathbb{Z}\setminus\{0\}$ be a number defined in \eqref{JJ}. Without loss of generality, we assume $\acute{l}_J\neq 0$. Then
\begin{align*}
|\partial_{\xi_{J}}(\langle k,\breve{\omega}\rangle+\langle \acute{l},[\acute{\Omega}]\rangle+\langle \hat{l},[\hat{\Omega}]\rangle)|\geq &|\acute {l}_J||\partial_{\xi_{J}}[\acute{\Omega}]_J|-|\partial_{\xi_{J}}(\langle \acute{l},[\acute{\Omega}]\rangle+\langle \hat{l},[\hat{\Omega}]\rangle-\acute{l}_J[\acute{\Omega}]_J)|\\
\geq& |\breve{l}_J|(1-c\eta\epsilon)-(\sum_{j\neq J}|\acute{l}_j|\frac{|j|c\eta^6\epsilon}{|J|}+\sum_{j\geq \mathcal{N}, |\hat{l}|\le 2}|\hat{l}|_j\frac{|j|c\eta^6\epsilon}{|J|}) \\
\geq& |\breve{l}_J|-(\mathcal{M}+2)c\eta^6\epsilon\\
\geq &\frac{1}{4}, \quad (\mathcal{M}\le (2c\eta^6\epsilon)^{-1}).
\end{align*}
%Similarly, one has
%\begin{align*}
%|\partial_{\xi^\prime_{t}}(\langle k,\breve{\omega}\rangle+\langle \acute{l},[\acute{\Omega}]\rangle+\langle \hat{l},[\hat{\Omega}]\rangle)|\geq &\frac{1}{4 |t|}\geq \frac{1}{4 M_{\acute{l},\hat{l}}}.
%\end{align*}
Hence,
$${\rm Meas}\, \mathcal{R}_{k\acute{l}\hat{l}}\le \frac{4\acute{\eta}M_{\acute{l},\hat{l}}}{4^{\mathcal{M}}(|k|+1)^{\tau}C(\mathcal{N},\acute{l})}{\rm Meas}\,\Pi_\eta.$$

%If $|k|=0, |\acute{l}|=0$ and $1\le |\hat{l}|\le 2$, then it is easy to see that $|\langle \hat{l},[\hat{\Omega}]\rangle|$ is not small, i.e.,
%$$\mathcal{R}_{k\acute{l}\hat{l}}=\emptyset,\quad \forall |k|=0, |\acute{l}|=0, 1\le |\hat{l}|\le 2.$$
%In the following, we will count the number of the resonant set $\mathcal{R}_{k\acute{l}\hat{l}}$.

%Since $M_{\acute{l},\hat{l}}\le \max\{j_*,j_{**},\mathcal{N}\}$,  it means that the number of non empty set $\mathcal{R}_{k\acute{l}\hat{l}}$ is less than
%$$A:=(2\mathcal{N}+1)^{|\acute{l}|}(\frac{4}{c_1}((|k|+1)(|\omega(\xi)|+1)+c_2(\mathcal{M}+2)\mathcal{N}^2+1))^2.$$

Therefore, we have
\begin{align*}
{\rm Meas}\,\mathcal{R}\le &\sum_{\substack{|k|+|\acute{l}|+|\hat{l}|\neq 0,\\|\acute{l}|+|\hat{l}|\le \mathcal{M}+2, |\hat{l}|\le 2}} {\rm Meas}\,\mathcal{R}_{k\acute{l}\hat{l}}\\
\le & \sum_{k\in \mathbb{Z}^n, |\acute{l}|\le \mathcal{M}+2}  \frac{4\acute{\eta}M_{\acute{l},\hat{l}}(A_1+A_2)}{4^{\mathcal{M}}(|k|+1)^{\tau}C(\mathcal{N},\acute{l})}\ {\rm Meas}\,\Pi_\eta\\
\le & \sum_{k\in \mathbb{Z}^n, |\acute{l}|\le \mathcal{M}+2}  \frac{4\acute{\eta}\max\{j_*,j_{**},\mathcal{N}\}(A_1+A_2)}{4^{\mathcal{M}}(|k|+1)^{\tau}C(\mathcal{N},\acute{l})}\ {\rm Meas}\,\Pi_\eta\\
\le & c\acute{\eta}\, {\rm Meas}\,\Pi_\eta,
\end{align*}
where $c$ is a constant depending on $c_1, c_2$ and $n$. Thus, we finish the proof of measure estimate.

\subsection{Proof of Theorem \ref{lthe}} \label{lthea}
Based on Theorem \ref{parth}, for a given positive integer $0\le \mathcal{M}\le (c\eta^6\epsilon)^{-1/2}$ and $0<\acute{\eta}<1$, there exists a small $\delta_0$ depending on $s_0,r_0,n,\acute{\eta}$ and $\mathcal{M}$ such that, for each $0<\delta<\delta_0$, $\xi\in\Pi_{\acute{\eta}}$ and the positive $\mathcal{N}$ satisfying
\begin{equation}\label{mathn}
\delta^{-\frac{\mathcal{M}+1}{2(p-1)}}\le \mathcal{N}+1<\delta^{-\frac{\mathcal{M}+1}{2(p-1)}}+1,
\end{equation}
 and there is a map
$$\Phi:D(4\varrho,4\delta,4\delta)\rightarrow D(5\varrho, 5\delta,5\delta),$$
such that
\begin{equation*}
  \begin{split}
     \breve{\breve{H}}(x,y,z,\bar{z};\xi) &= \breve{H}\circ\Phi=\breve{N}(x,y,z,\bar{z};\xi)+Z(x,y,z,\bar{z};\xi) \\
       & +Q(x,y,z,\bar{z};\xi)+T(x,y,z,\bar{z};\xi),
  \end{split}
\end{equation*}

with
\begin{equation}\label{nyou}
\breve{N}(x,y,z,\bar{z};\xi)=\sum_{j=1}^n\breve{\omega}_j(\xi)y_j+\sum_{j\in\mathbb{Z}_*}\frac{\breve{\Omega}_j(x;\xi)}{j}z_j\bar{z}_j,
\end{equation}
%$$=\sum_{j=1}^n\breve{\omega}_j(\xi)y_j+\sum_{|j|\le \mathcal{N}}\breve{\Omega}_j(\xi)z_j\bar{z}_j+\sum_{|j|\geq \mathcal{N}+1}\breve{\Omega}_j(x;\xi)z_j\bar{z}_j,$$
\begin{align}
Z(x,y,z,\bar{z};\xi)=&\sum_{4\le 2|\alpha|+2|\beta|+2|\mu|\le\mathcal{M}+2, |\mu|=0 }Z^{\alpha\beta\beta\mu\mu}(\xi)y^{\alpha}\acute{z}^{\beta}\bar{\acute{z}}^{\beta}\hat{z}^{\mu}\bar{\hat{z}}^{\mu}\nonumber\\
&+\sum_{4\le 2|\alpha|+2|\beta|+2|\mu|\le\mathcal{M}+2, |\mu|=1 }Z^{\alpha\beta\beta\mu\mu}(x;\xi)y^{\alpha}\acute{z}^{\beta}\bar{\acute{z}}^{\beta}\hat{z}^{\mu}\bar{\hat{z}}^{\mu},\label{zyou}
\end{align}
%$$R(x,y,z,\bar{z};\xi)=\sum_{2|\alpha|+|\beta|+|\gamma|+|\mu|+|\nu|\geq\mathcal{M}+3, |\mu|+|\nu|\le2 }R^{\alpha\beta\gamma\mu\nu}(x;\xi)y^{\alpha}\acute{z}^{\beta}\bar{\acute{z}}^{\gamma}\hat{z}^{\mu}\bar{\hat{z}}^{\nu},$$
$$Q(x,y,z,\bar{z};\xi)=\sum_{3\le |\alpha|+|\beta|+|\gamma|+|\mu|+|\nu|\le\mathcal{M}+2, |\mu|+|\nu|\geq3 }Q^{\alpha\beta\gamma\mu\nu}(x;\xi)y^{\alpha}\acute{z}^{\beta}\bar{\acute{z}}^{\gamma}\hat{z}^{\mu}\bar{\hat{z}}^{\nu}.$$
and
%$$T(x,y,z,\bar{z};\xi)=\sum_{3\le2|\alpha|+|\beta|+|\gamma|+|\mu|+|\nu|\le\mathcal{M}+2, |\mu|+|\nu|\le2 }A^{\alpha\beta\gamma\mu\nu}(x;\xi)y^{\alpha}\acute{z}^{\beta}\bar{\acute{z}}^{\gamma}\hat{z}^{\mu}\bar{\hat{z}}^{\nu}.$$
\begin{align*}
T(x,y,z,\bar{z};\xi)=&\sum_{2|\alpha|+|\beta|+|\gamma|+|\mu|+|\nu|\geq \mathcal{M}+3}T^{\alpha\beta\gamma\mu\nu}(x;\xi)y^{\alpha}\acute{z}^{\beta}\bar{\acute{z}}^{\gamma}\hat{z}^{\mu}\bar{\hat{z}}^{\nu},\\
&+\sum_{ \substack{3\le 2|\alpha|+|\beta|+|\gamma|+|\mu|+|\nu|\le \mathcal{M}+2 \\\beta\neq \gamma, \mu\neq\nu, |\mu|+|\nu|\le 2}}T^{\alpha\beta\gamma\mu\nu}(x;\xi)y^{\alpha}\acute{z}^{\beta}\bar{\acute{z}}^{\gamma}\hat{z}^{\mu}\bar{\hat{z}}^{\nu}.
\end{align*}
Moreover, the following estimates hold:

(1) for each $\xi\in\Pi_{\acute{\eta}}$, the symplectic map $\Phi$ satisfies
$$\Vert \Phi-id\Vert_{p,D(4\varrho,4\delta,4\delta)\times\Pi_{\acute{\eta}}}\le \frac{c\mathcal{N}^{C}\delta}{\acute{\eta}^2},$$
$$|| D\Phi-Id||_{p,D(4\varrho,4\delta,4\delta)\times\Pi_{\acute{\eta}}}\le \frac{c\mathcal{N}^{C}}{\acute{\eta}^2}.$$

(2)The functions ${Z}, Q$ and $T$ satisfy
\begin{equation}\label{y1}
\Vert X_{Z_j}\Vert_{p-1,D(4\varrho,4\delta,4\delta)\times\Pi_{\acute{\eta}}},\, \frac{1}{\delta^2}|||Z_j|||_{D(4\varrho,4\delta,4\delta)\times\Pi_{\acute{\eta}}}\le c\delta\left(\frac{\mathcal{N}^{2C_0(\mathcal{M}+6)^2}}{\acute{\eta}^2}\delta\right)^{j-3},
\end{equation}
%$$|||R_j|||_{p,D(s_0/4,4\delta,4\delta)\times\Pi_{\acute{\eta}}}\le c\delta^3\left(\frac{\mathcal{N}^{2(\mathcal{M}+6)^2+2C(\tau+n)}}{\acute{\eta}^2}\delta\right)^{j-3},$$
\begin{equation}\label{y2}
\Vert X_{Q_j}\Vert_{p-1,D(4\varrho,4\delta,4\delta)\times\Pi_{\acute{\eta}}},\,\frac{1}{\delta^2}|||Q_j|||_{D(4\varrho,4\delta,4\delta)\times\Pi_{\acute{\eta}}}\le c\delta\left(\frac{\mathcal{N}^{2C_0(\mathcal{M}+7)^2}}{\acute{\eta}^2}\delta\right)^{j-3},
\end{equation}
and
\begin{equation}\label{y3}
\Vert X_T\Vert_{p-1,D(4\varrho,4\delta,4\delta)\times\Pi_{\acute{\eta}}}\le c\delta\left(\frac{\mathcal{N}^{2C_0(\mathcal{M}+7)^2}}{\acute{\eta}^2}\delta\right)^{\mathcal{M}},
\end{equation}
where $c$ is a constant depending on $s_0,r_0,n,\mathcal{M}$.

By \eqref{y1},\eqref{y2} and \eqref{y3}, we have
\begin{align}
\Vert X_T\Vert_{p-1,D(4\varrho,4\delta,4\delta)\times\Pi_{\acute{\eta}}}\le&  c\delta \left(\frac{\mathcal{N}^{2C_0(\mathcal{M}+7)^2}}{\acute{\eta}^2}\delta\right)^{\mathcal{M}}\nonumber\\
\le&\delta^{\mathcal{M}+1}\left(\frac{c\mathcal{N}^{2C_0(\mathcal{M}+7)^2\mathcal{M}}}{\acute{\eta}^{2\mathcal{M}}}\right)\nonumber\\
\overset{\eqref{mathn}}{\le} &\delta^{\mathcal{M}+1}\left(\frac{c}{\acute{\eta}^{2\mathcal{M}}}\delta^{-\frac{2C_0(\mathcal{M}+7)^2\mathcal{M}(\mathcal{M}+1)}{p-1}}\right)\nonumber\\
\le & \delta^{\mathcal{M}+1}\left(\frac{c}{\acute{\eta}^{2\mathcal{M}}}\delta^{\frac{1}{4}}\right)\nonumber\\
&({\rm by}\ \ p\geq 8C_0(\mathcal{M}+7)^4+1)\nonumber\\
\le &\delta^{\mathcal{M}+1},\label{tti}
\end{align}

\begin{align*}
\Vert X_Q\Vert_{p-1,D(4\varrho,4\delta,4\delta)\times\Pi_{\acute{\eta}}}\le &\sum_{j\geq 3} c\delta \left(\frac{\mathcal{N}^{2C_0(\mathcal{M}+7)^2}}{\acute{\eta}^2}\delta\right)^{j-3}\\
\le &\sum_{j\geq 3}c\delta\left(\frac{\mathcal{N}^{2C_0(\mathcal{M}+7)^2}}{\acute{\eta}^2}\delta\right)^{j-3}\\
\le & \delta,
\end{align*}
and
\begin{align*}
\Vert X_Z\Vert_{p-1,D(4\varrho,4\delta,4\delta)\times\Pi_{\acute{\eta}}}\le &\sum_{j\geq 4} c\delta \left(\frac{\mathcal{N}^{2C_0(\mathcal{M}+6)^2}}{\acute{\eta}^2}\delta\right)^{j-3}\\
\le &\sum_{j\geq 4}c\delta\left(\frac{\mathcal{N}^{2C_0(\mathcal{M}+6)^2}}{\acute{\eta}^2}\delta\right)^{j-3}\\
\le & \delta^2.
\end{align*}

In addition, one can obtain the following inequality from \eqref{mathn}
\begin{align}
\Vert \hat{z}\Vert_1=&\sqrt{\sum_{j\geq \mathcal{N}+1}|z_j|^2j^2}\nonumber\\
=&\sqrt{\sum_{j\geq \mathcal{N}+1}|j|^{2p}|z_j|^{2p}/j^{2(p-1)}}\nonumber\\
\le &\frac{ \Vert \hat{z}\Vert_p}{(\mathcal{N}+1)^{p-1}}\nonumber\\
\le &\delta^{\mathcal{M}+1} \Vert \hat{z}\Vert_p.\label{2p2}
\end{align}
By Lemma \ref{guji}, we get on the domain $D(3\varrho,4\delta,4\delta)\times\Pi_{\acute{\eta}}$
\begin{align*}
 \frac{1}{4\delta} \Vert (-{\bf i} m \partial_{z_m}Q)_{m\in\mathbb{Z}_*}\Vert_{p/2} \le & c\sum_{3\le j\le \mathcal{M}+2} |||Q_j|||_{D(s,r,r)\times \Pi}(4\delta)^{-(j+1)}\nonumber, \\
    & \cdot\sum_{\alpha}(\frac{1}{(\varrho)^{p/2}}|y|^{|\alpha|}\Vert z\Vert_1^{j-2|\alpha|-1}+\Vert z\Vert_p \Vert z\Vert_1^{j-2|\alpha|-2}),\\
      \le & c\sum_{3\le j\le \mathcal{M}+2} \delta^3 \left(\frac{\mathcal{N}^{2C_0(\mathcal{M}+7)^2}}{\acute{\eta} ^2}\delta\right)^{j-3}(4\delta)^{-(j+1)}\\
& \cdot \sum_{\alpha}\left(\left(\frac{\mathcal N^2}{\mathcal{M}}\right)^{p/2}(4\delta)^{j-1}\delta^{2(\mathcal{M}+1)}+(4\delta)^{j-1}\delta^{\mathcal{M}+1}\right)\\
&({\rm since}\ Q_j=O(|\hat{z}|^3) \ {\rm and}\ \eqref{2p2} )\\
\le & c\delta^{\mathcal{M}+1} \sum_{3\le j\le \mathcal{M}+2} c\delta \left(\frac{\mathcal{N}^{2C_0(\mathcal{M}+7)^2}}{\acute{\eta} ^2}\delta\right)^{j-3}\\
\le &  \delta^{\mathcal{M}+1},
\end{align*}
where $c$ is a constant depending on $s_0,n,r_0$ and $\mathcal{M}$.
Similarly, the following estimates hold
$$\frac{1}{(4\delta)^2}\Vert \partial_x Q\Vert, \quad \Vert \partial_y Q\Vert,\quad \frac{1}{4\delta} \Vert ({\bf i} m \partial_{\bar{z}_m}Q)_{m\in\mathbb{Z}_*}\Vert_{p/2}\le   \delta^{\mathcal{M}+1}.$$
In conclusion, we have
\begin{align}\label{452}
\Vert X_Q\Vert_{p-1,D(3\varrho,4\delta,4\delta)\times\Pi_{\acute{\eta}}}
\le & \delta^{\mathcal{M}+1}.
\end{align}

Observe that the function $Z$ in the Hamiltonian  $\breve{\breve{H}}(x,y,z,\bar{z};\xi)$ depends on the angular variable $x$. Therefore, the above normal form is not enough to give the long time stability of KAM tori. In order to prove the long time stability of KAM tori, we will do some further transformations. By the definition of $\breve{N}$ and $Z$ in \eqref{nyou} and \eqref{zyou}, we can rewrite the Hamiltonian $\breve{\breve{H}}$ as
$$\breve{\breve{H}}=\breve{N}(x,y,z,\bar{z};\xi)+Z(x,y,z,\bar{z};\xi)+Q(x,y,z,\bar{z};\xi)+T(x,y,z,\bar{z};\xi),$$
with
\begin{align*}
\breve{N}(x,y,z,\bar{z};\xi)%=&\sum_{j=1}^n\breve{\omega}_j(\xi)y_j+\sum_{j\geq 1}\frac{[\breve{\Omega}_j](\xi)}{j}|z_j|^2+\sum_{j\geq 1}\frac{(\breve{\Omega}_j-[\breve{\Omega}_j])(x;\xi)}{j}|z_j|^2,\\
=&\sum_{j=1}^n\breve{\omega}_j(\xi)y_j+\sum_{j\geq 1}\frac{[\breve{\Omega}_j](\xi)}{j}|z_j|^2+\sum_{j\geq 1}\frac{(\breve{\Omega}_j-[\breve{\Omega}_j])(x;\xi)}{j}(|z_j|^2-|z_j(0)|^2)\\
&+\sum_{j\geq 1}\frac{(\breve{\Omega}_j-[\breve{\Omega}_j])(x;\xi)}{j}|z_j(0)|^2,
\end{align*}
and
\begin{align*}
Z(x,y,z,\bar{z};\xi)=&\sum_{4\le 2|\alpha|+2|\beta|\le\mathcal{M}+2, |\mu|\le1 }Z^{\alpha\beta}(x;\xi)y^{\alpha}|z|^{2\beta}\\
=&\sum_{4\le 2|\alpha|+2|\beta|\le\mathcal{M}+2, |\mu|\le1 }Z^{\alpha\beta}(x)y^{\alpha}(|z|^{2}-|z(0)|^2+|z(0)|^2)^{\beta}\\
=&\sum_{|\alpha|\geq 2}Z^{\alpha}(\xi)y^{\alpha}+\sum_{|\alpha|\geq 1,|\beta|\geq 1}Z^{\alpha\beta}(x;\xi)y^{\alpha}|z(0)|^{2\beta}\\%+\sum_{|\alpha|=1,|\beta|\geq 1}Z^{\alpha\beta}y^\alpha |z(0)|^{2\beta}\\
&+\sum_{|\alpha|,|\gamma|\geq 0,|\beta|\geq 1}\widetilde{Z}^{\alpha\beta\gamma}(x;\xi)y^{\alpha}|z(0)|^{2\gamma}(|z|^{2}-|z(0)|^2)^{\beta}\\
=&\sum_{|\alpha|\geq 2}\widetilde{Z}^{\alpha}(x;\xi)y^{\alpha}%+\sum_{|\alpha|=1,\beta\geq 1}Z^{\alpha\beta}y^\alpha |z(0)|^{2\beta}
+\sum_{|\alpha|,|\gamma|\geq 0, |\beta|\geq 1}\widetilde{Z}^{\alpha\beta\gamma}(x;\xi)y^{\alpha}|z(0)|^{2\gamma}(|z|^{2}-|z(0)|^2)^{\beta},
\end{align*}
where
$$\sum_{|\alpha|\geq 1}\widetilde{Z}^{\alpha}(x;\xi)y^{\alpha}=\sum_{|\alpha|\geq 2}Z^{\alpha}(\xi)y^{\alpha}+\sum_{|\alpha|\geq 1,|\beta|\geq 1}Z^{\alpha\beta}(x;\xi)y^{\alpha}|z(0)|^{2\beta}.$$
%Let $\omega^0_j=\breve{\omega}_j+\sum_{\beta\geq 1}[Z^{\alpha^j\beta}] |z(0)|^{2\beta}$ where $\alpha^j=e_j$ and $[Z^{\alpha^j\beta}]$ is the mean value of the function $Z^{\alpha^j\beta}$ over $\mathbb{T}^n$.
Therefore, the Hamiltonian becomes
\begin{align*}
\breve{\breve{H}}=&\sum_{j=1}^n\breve{\omega}_j(\xi)y_j+\sum_{j\geq 1}\frac{[\breve{\Omega}_j](\xi)}{j}|z_j|^2+\sum_{j\geq 1}\frac{(\breve{\Omega}_j-[\breve{\Omega}_j])(x;\xi)}{j}(|z_j|^2-|z_j(0)|^2)\\
&+\sum_{j\geq 1}\frac{(\breve{\Omega}_j-[\breve{\Omega}_j])(x;\xi)}{j}|z_j(0)|^2+\sum_{|\alpha|\geq 1}\widetilde{Z}^\alpha(x;\xi) y^{\alpha}\\
&+\sum_{|\alpha|,|\gamma|\geq 0,|\beta|\geq 1}\widetilde{Z}^{\alpha\beta\gamma}(x;\xi)y^{\alpha}|z|^{2\gamma}(|z|^{2}-|z(0)|^2)^{\beta}+Q+T\\
=&N+A+\widetilde{Z}+B+Q+T,
\end{align*}
where
\begin{align*}
N=&\sum_{j=1}^n\breve{\omega}_j(\xi)y_j+\sum_{j\geq 1}\frac{[\breve{\Omega}_j](\xi)}{j}|z_j|^2+\sum_{j\geq 1}\frac{(\breve{\Omega}_j-[\breve{\Omega}_j])(x;\xi)}{j}(|z_j|^2-|z_j(0)|^2)\\
A=&\sum_{j\geq 1}\frac{(\breve{\Omega}_j-[\breve{\Omega}_j])(x;\xi)}{j}|z_j(0)|^2+\widetilde{Z}^1(x;\xi)  y\\
\widetilde{Z}=&\sum_{|\alpha|\geq 2}\widetilde{Z}^\alpha(x;\xi) y^{\alpha},\quad  B=\sum_{|\alpha|,|\gamma|\geq 0, |\beta|\geq 1}\widetilde{Z}^{\alpha\beta\gamma}(x;\xi)y^{\alpha}|z(0)|^{2\gamma}(|z|^{2}-|z(0)|^2)^{\beta},\\
Q=&\sum_{\substack{3\le |\alpha|+|\beta|+|\gamma|+|\mu|+|\nu|\le\mathcal{M}+2\\ |\mu|+|\nu|\geq3} }Q^{\alpha\beta\gamma\mu\nu}(x;\xi)y^{\alpha}\acute{z}^{\beta}\bar{\acute{z}}^{\gamma}\hat{z}^{\mu}\bar{\hat{z}}^{\nu},\\
T=&\sum_{2|\alpha|+|\beta|+|\gamma|+|\mu|+|\nu|\geq \mathcal{M}+3}T^{\alpha\beta\gamma\mu\nu}(x;\xi)y^{\alpha}\acute{z}^{\beta}\bar{\acute{z}}^{\gamma}\hat{z}^{\mu}\bar{\hat{z}}^{\nu},\\
&+\sum_{3\le \substack{2|\alpha|+|\beta|+|\gamma|+|\mu|+|\nu|\le \mathcal{M}+2 \\\beta\neq \gamma, \mu\neq\nu, |\mu|+|\nu|\le 2}}T^{\alpha\beta\gamma\mu\nu}(x;\xi)y^{\alpha}\acute{z}^{\beta}\bar{\acute{z}}^{\gamma}\hat{z}^{\mu}\bar{\hat{z}}^{\nu},
\end{align*}
where $Q$ and $T$ satisfying the inequalities \eqref{452} and \eqref{y3}, respectively.

Since
$$\sup_{j} \left\Vert\frac{(\breve{\Omega}_j-[\breve{\Omega}_j])(x;\xi)}{j}\right\Vert_{D(4\varrho)\times \Pi}\le c\eta^6\epsilon,$$
and $\Vert z(0)\Vert_{p}\le \delta,$
we have
$$\left\Vert \sum_{j\geq 1}\frac{(\breve{\Omega}_j-[\breve{\Omega}_j])(x;\xi)}{j}|z_j(0)|^2\right\Vert_{D(4\varrho)\times \Pi}\le \delta^2.$$
From the estimate of $Z$ in \eqref{y1}, one can easily obtain
$$\Vert\widetilde{Z}^1(x;\xi)\Vert_{D(4\varrho)\times \Pi_{\acute{\eta}}}\le \delta^{2-},$$
and
$$\Vert\widetilde{Z}^\alpha(x;\xi)-[\widetilde{Z}^\alpha](\xi)\Vert_{D(4\varrho)\times \Pi_{\acute{\eta}}}\le \delta^{2-},\quad {\rm for}\ \ |\alpha|\geq 2.$$
Thus, on the domain $D(3\varrho,4\delta)\times \Pi_{\acute{\eta}}$
\begin{equation}\label{hhha}
\Vert X_{A(x,y;\xi)}\Vert_{D(3\varrho,4\delta)\times \Pi_{\acute{\eta}}}\le \delta^{2-}.
\end{equation}
\begin{equation}\label{hhhz}
\Vert X_{\widetilde{Z}(x,y;\xi)}\Vert_{D(3\varrho,4\delta)\times \Pi_{\acute{\eta}}}\le \delta^{2-},
\end{equation}
where the domain $D(s,r):=\{(x,y): |{\rm Im}\ x|\le s, |y|\le r^2\}$. In addition, for the function $B$, if we let
$$B=\sum_{|\beta|\geq 1} B^\beta(x,y;\xi)(|z|^{2}-|z(0)|^2)^{\beta}, \quad B^\beta(x,y;\xi)=\sum_{|\alpha|,|\gamma|\geq 0}\widetilde{Z}^{\alpha\beta\gamma}(x;\xi)y^{\alpha}|z(0)|^\gamma.$$
From \eqref{y1} and the definition of $B$, we can easily to compute
\begin{equation}\label{Bhu}
\Vert B^\beta(x,y;\xi)\Vert_{D(3\varrho,4\delta)\times \Pi_{\acute{\eta}}} \le \delta^{2-}.
\end{equation}

We rewrite the function $A$ by $A(x,y)=\Gamma_KA(x,y)+(1-\Gamma_K)A(x,y)$ with
\begin{equation}\label{newk}
K=\left|\frac{1}{\varrho}\log \delta^{-\mathcal{M}}\right|
\end{equation}
and here $\Gamma_KA$ denotes the Fourier series of $A=A(x,y)$ truncated at $K$-th order. Then, on the domain $D(2\varrho,3\delta)$, we have $\Vert X_{(1-\Gamma_K)A(x,y)}\Vert\le \delta^{\mathcal{M}+2}$.

 In the following, we will perform some transformations such that $A(x,y)$ become small enough for long time stability estimate. For this purpose, we will perform KAM iteration in $(x,y)$-variables. Let $D_j=\{(x,y): |{\rm Im}\ x|\le \varrho_j, |y|\le \delta_j^2\}$ where $\varrho_j=3\varrho-j \varrho/\log_2\mathcal{M}$ and $\delta_j=4\delta-j \delta/\log_2\mathcal{M}$.

Let $\Phi=X_F^t|_{t=1}$ with
$$F=F^0(x)+F^1(x)y.$$
Then, we obtain
\begin{align}
\breve{\breve{H}}\circ \Phi=&(N+A+\widetilde{Z}+B+T+Q
)\circ \Phi\nonumber\\
=&N+ \{N,F\}+\int_0^1(1-t)\{\{ N,F\},F\}\circ X_F^tdt\nonumber\\
&+\Gamma_KA+(1-\Gamma_K)A+\int_0^t t\{A, F\}\circ X_F^tdt\nonumber\\
&+(\widetilde{Z}+B+Q+T)\circ \Phi.\label{remain term}
\end{align}
Then, we want to solve the following homological equation.
$$ \{N,F\}=\Gamma_KA-[A],$$
where $[A]=[\widetilde{Z}^1]y$. Here, the notation $[\cdot]$ represents the average of a function with respect to the variable $x$ over $\mathbb{T}^n$. By Fourier expansion, in both sides of the homological equation, we have
$$\sum_{0\neq|k|\le K}{\bf i}\langle k,\breve{\omega}\rangle \widehat{F^0_k}e^{{\bf i}\langle k,x\rangle}=\sum_{0\neq|k|\le K} \widehat{A^0_k}e^{{\bf i}\langle k,x\rangle},$$
$$\sum_{0\neq|k|\le K}{\bf i}\langle k,\breve{\omega}\rangle \widehat{F^1_k}e^{{\bf i}\langle k,x\rangle}=\sum_{0\neq|k|\le K} \widehat{A^1_k}e^{{\bf i}\langle k,x\rangle}.$$
Thus, we have
$${\bf i}\langle k,\breve{\omega}\rangle \widehat{F^0_k}=\widehat{A^0_k},\quad \ {\bf i}\langle k,\breve{\omega}\rangle \widehat{F^1_k}=\widehat{A^1_k}.$$
%Since $\breve{\omega}$ satisfies the diophantine condition, we obtain
%\begin{align*}
%|\langle k,\omega^*\rangle|\geq& |\langle k,\breve{\omega}\rangle|-|\langle k,\omega^*-\breve{\omega}\rangle|\\
%\geq&\frac{\eta}{|k|^\tau}-K\delta^2\geq \frac{\eta}{2|k|^\tau}.
%\end{align*}
For each $\xi \in \Pi_{\acute{\eta}}$,
if $\breve{\omega}$ satisfies the Diophantine condition, i.e.,
\begin{equation}\label{dioo}
|\langle k, \breve{\omega}\rangle|\geq \frac{\acute{\eta}}{1+|k|^\tau},
\end{equation}
then one can easily obtain
\begin{align}\label{y10}
\Vert X_F\Vert_{D_1} \le \acute{\eta}^{-1} K^\tau \delta^2.
\end{align}
It implies the transformation $\Phi: D_1\rightarrow D_0$.
Then, we are going to estimate the remaining terms in \eqref{remain term}. By \eqref{hhha} and \eqref{y10}, we have
$$\Big\Vert X_{\int_0^1(1-t)\{\{N,F\},F\}\circ X_F^tdt}\Big\Vert_{D_1}\le \delta^{4-},$$
$$\Big\Vert X_{\int_0^1t\{A,F\}\circ X_F^tdt}\Big\Vert_{D_1}\le \delta^{4-}.$$

Under the transformation $\Phi$ and by \eqref{hhhz}, the function $\tilde{Z}=\sum_{|\alpha|\geq 2}\tilde{Z}^\alpha y^\alpha$ has the estimate
\begin{align*}
\Vert X_{\tilde{Z}\circ \Phi-\tilde{Z}}\Vert_{D(4\varrho,4\delta)}=\Big\Vert \int_0^1t\{\tilde{Z},F\}\circ X_F^tdt\Big\Vert_{D(4\varrho)}
\le& \delta^{4-} .
\end{align*}
%Moreover, if we denote $\tilde{Z}\circ \Phi=\tilde{Z}+\sum_{|\alpha|\geq 0}Z_+^\alpha(x) y^\alpha$, then we immediately have $Z_+^\alpha(x)$ will depend on the power of $\partial_xF$ and $\Vert Z_+^0(x)\Vert\le (\eta^{-1} K^\tau \delta^2)^2\le \delta^{4-}$.
For the term $B, Q$, we still have the form
$$B\circ \Phi=\sum_{|\beta|\geq 1}B_+^{\beta}(x,y;\xi)(|z|^{2}-|z(0)|^2)^{\beta},$$
$$Q\circ\Phi=\sum_{|\mu|+|\nu|\geq3 }Q_+^{\alpha\beta\gamma\mu\nu}(x;\xi)y^{\alpha}\acute{z}^{\beta}\bar{\acute{z}}^{\gamma}\hat{z}^{\mu}\bar{\hat{z}}^{\nu},$$
since we only do change variables for the variables $x, y$. For the term $B, Q$ and $T$,  since transformation $\Phi:  D(\varrho_1,\delta_1,4\delta)\rightarrow D(3\varrho,4\delta,4\delta)$ and the estimate \eqref{tti}, \eqref{452} and \eqref{Bhu}, we have
$$\Vert B_+^\beta(x,y;\xi)\Vert_{D(\varrho_1,\delta_1)\times \Pi_{\acute{\eta}}} \le \delta^{2-}.
$$
$$\Vert X_{Q\circ\Phi}\Vert_{p/2, D(\varrho_1,\delta_1,4\delta)},\quad \Vert X_{T\circ\Phi}\Vert_{p-1, D(\varrho_1,\delta_1,4\delta)}\le \delta^{\mathcal{M}+1}.$$
Thus, under the transformation $\Phi$, the Hamiltonian becomes
$$H=N_++A_++Z_++B_++Q_++T_+,$$
with
$$N_+=\sum_{j=1}^n\omega_{+j}(\xi)y_j+\sum_{j\geq 1}\frac{[\breve{\Omega}_j](\xi)}{j}|z_j|^2+\sum_{j\geq 1}\frac{\breve{(\Omega}_j-[\breve{\Omega}_j])(x;\xi)}{j}(|z_j|^2-|z_j(0)|^2),$$
$$A_+=A_+^0(x)+A_+^1(x)y,\quad Z_+=\sum_{|\alpha|\geq 2}Z_+^\alpha(x) y^\alpha,$$
$$B_+=\sum_{|\beta|\geq 1}B_+^{\beta}(x,y;\xi)(|z|^{2}-|z(0)|^2)^{\beta},$$
$$Q_+=\sum_{|\mu|+|\nu|\geq3 }Q_+^{\alpha\beta\gamma\mu\nu}(x;\xi)y^{\alpha}\acute{z}^{\beta}\bar{\acute{z}}^{\gamma}\hat{z}^{\mu}\bar{\hat{z}}^{\nu},$$
where the functions $A_+,\ Z_+, B_+, Q_+$ and $T_+$ satisfies $\Vert X_{A_+}\Vert_{D_1}\le \delta^{4-}, \ \Vert X_{Z_+}\Vert_{D_1}\le \delta^{2-}$, $\Vert B_+^\beta(x,y;\xi)\Vert_{D(\varrho_1,\delta_1)} \le \delta^{2-}$ and $\Vert X_{Q_+}\Vert_{p/2,D(\varrho_1,\delta_1,4\delta)}, \Vert X_{T_+}\Vert_{p-1,D(\varrho_1,\delta_1,4\delta)}\le \delta^{\mathcal{M}+1}$, respectively. Here,
$$\omega_+=\breve{\omega}+[\tilde{Z}^1].$$
In what follows, we will verify the new frequency $\omega_+$ satisfy the condition
$$|\langle k,\omega_+\rangle|\geq \frac{\acute{\eta}}{2(1+|k|^\tau)},\quad 0\neq k\in\mathbb{Z}^n,$$
if $\breve{\omega}$ satisfies the Diophantine condition \eqref{dioo}.  Actually, by \eqref{dioo}, one derives
\begin{align*}
|\langle k,\omega_+\rangle|\geq& |\langle k,\breve{\omega}\rangle|-|\langle k,\omega_+-\breve{\omega}\rangle|\\
\geq&\frac{\acute{\eta}}{1+|k|^\tau}-K\delta^{2-}\geq \frac{\acute{\eta}}{2(1+|k|^\tau)},
\end{align*}
since we have $K(1+K|^\tau)\delta^{2-}\le \frac{\acute{\eta}}{2}$ by the choice of $K$ in \eqref{newk}.

After we perform at most $\log_2\mathcal{M}$ times iteration, there exists a symplectic transformation such that $\Phi^*$ the Hamiltonian $H$ becomes
$$H^*=H\circ\Phi^{*}=N+Z+B+Q+T,$$
where %$ \Gamma(x,y)$ has the estimate $\Vert\Gamma\Vert\le \delta^{\mathcal{M}}$ with can be put into T.
$$N=\sum_{j=1}^n\omega^*_j(\xi)y_j+\sum_{j\geq 1}\frac{[\breve{\Omega}_j](\xi)}{j}|z_j|^2+\sum_{j\geq 1}\frac{\breve{(\Omega}_j-[\breve{\Omega}_j])(x;\xi)}{j}(|z_j|^2-|z_j(0)|^2),$$
$$Z=\sum_{|\alpha|\geq 2}Z^\alpha(x;\xi) y^\alpha,\quad B=\sum_{ |\beta|\geq 1}B^{\beta}(x,y;\xi)(|z|^{2}-|z(0)|^2)^{\beta},$$
$$Q=\sum_{|\mu|+|\nu|\geq3 }Q^{\alpha\beta\gamma\mu\nu}(x;\xi)y^{\alpha}\acute{z}^{\beta}\bar{\acute{z}}^{\gamma}\hat{z}^{\mu}\bar{\hat{z}}^{\nu},$$
with the coefficients of $Z$, $Q$, $B$ and $T$ satisfy the following estimates
\begin{align}
&\Vert Z^\alpha(x;\xi)-[Z^\alpha](\xi)\Vert_{D(2\varrho)}\le \delta^{2-},\quad \Vert B^\beta(x,y;\xi)\Vert_{D(\varrho_1,\delta_1)} \le \delta^{2-},\label{huuz}\\
&\Vert X_Q\Vert_{p/2,D(2\varrho,3\delta,4\delta)},\quad \Vert X_T\Vert_{p-1,D(2\varrho,3\delta,4\delta)}\le \delta^{\mathcal{M}+1}.
\end{align}
 In addition, the frequency satisfies the small divisor condition
\begin{equation}\label{o*}
|\langle k,\omega^*\rangle|\geq \frac{\acute{\eta}}{2(1+|k|^\tau)},\quad 0\neq k\in\mathbb{Z}^n.
\end{equation}

\medskip

In what follows, we are going to eliminate $Z^2(x;\xi)y^2$. For this purpose, we write $Z^2(x;\xi)y^{2}=\Gamma_KZ^{2}(x;\xi)y^2+(1-\Gamma_K)Z^2(x;\xi)y^2$.  Let $\tilde D_j=\{(x,y): |{\rm Im}\ x|\le \tilde\varrho_j, |y|\le \tilde\delta_j^2\}$ where $\tilde\varrho_j=2\varrho-j \varrho/\mathcal{M}$ and $\tilde\delta_j=3\delta-j \delta/\mathcal{M}$.

Let $\Phi=X_F^t|_{t=1}$ be the time-1 map of the vector field $X_{F}$ of Hamiltonion
\begin{equation}\label{yaof}
F(x,y;\xi)=F^2(x;\xi)y^2.
\end{equation}
Then, we obtain
\begin{equation}\label{y20}
\begin{split}
  H^*\circ \Phi=&(N+Z+B+Q+T
)\circ \Phi\\
=&N+ \{N,F\}+\int_0^1\{(1-t)\{ N,F\}+t Z^2(x)y^2,F\}\circ X_F^tdt\\
&+\Gamma_KZ^{2}(x)y^2+(1-\Gamma_K)Z^2(x)y^2+((Z-Z^2(x)y^2)+B+Q+T)\circ \Phi.
\end{split}
\end{equation}
From the above, the homological equation reads,
$$\{N,F\}=\Gamma_KZ^{2}y^2-[Z^{2}]y^2,$$
where $F$ is defined by \eqref{yaof}. By Fourier expansion, in both sides of the homological equation, we have
$$\sum_{0\neq|k|\le K}{\bf i}\langle k,\omega^*\rangle \widehat{F^2_k}e^{{\bf i}\langle k,x\rangle}=\sum_{0\neq|k|\le K} \widehat{Z^2_k}e^{{\bf i}\langle k,x\rangle}.$$
Thus, we have
$${\bf i}\langle k,\omega^*\rangle \widehat{F^2_k}=\widehat{Z^2_k}.$$
According to the small divisor condition \eqref{o*},
we have the solution $F(x,y;\xi)=F^2(x;\xi)y^2$ obeys the estimate
\begin{equation}\label{haoshuo}
\Vert F^2(x;\xi)\Vert_{\tilde D_1}\le \delta^{2-}.\end{equation}

Now let us deal with the remaining terms in \eqref{y20}. We have
\begin{align*}
\int_0^1\{(1-t)\{ N,F\}+t Z^2(x;\xi)y^2,F\}\circ X_F^tdt=\sum_{|\alpha|\geq 3}Z^\alpha_{+1}(x;\xi)y^\alpha,
\end{align*}
with $\Vert Z^\alpha_{+1}\Vert_{\tilde D_1}\le \delta^{2-}$ by \eqref{haoshuo} and \eqref{huuz}.

For the term $(1-\Gamma_K)Z^2y^2$, we have $\Vert (1-\Gamma_K)Z^2\Vert\le \delta^{\mathcal{M}+2}$. %For $\Gamma\circ \Phi=\Gamma_+(x,y)$, we have $|||\Gamma_+(x,y)|||\le \delta^{\mathcal{M}}$.

For the other terms, we have
$$(Z- Z^2(x;\xi)y^2)\circ \Phi=\sum_{|\alpha|\geq 3}Z_{+2}^\alpha(x;\xi) y^\alpha,$$
$$B_+=B\circ \Phi=\sum_{|\alpha|\geq 0, |\beta|\geq 1}B_+^{\beta}(x,y;\xi)(|z|^{2}-|z(0)|^2)^{\beta},$$
$$Q\circ \Phi=\sum_{|\mu|+|\nu|\geq3 }Q^{\alpha\beta\gamma\mu\nu}(x;\xi)y^{\alpha}\acute{z}^{\beta}\bar{\acute{z}}^{\gamma}\hat{z}^{\mu}\bar{\hat{z}}^{\nu},$$
%and
%$$R\circ\Phi=\sum_{\alpha,\beta,\gamma }R^{\alpha\beta\gamma}(x;\xi)y^{\alpha}z^{\beta}\bar{z}^{\gamma}$$
with the coefficients of $Z$, $B$, $Q$ and $T$ satisfy
$$\Vert B_+^{\beta}(x,y;\xi)\Vert_{\tilde D_1},\quad \Vert Z_{+2}^\alpha(x;\xi)-[Z_{+2}^\alpha](\xi)\Vert_{\tilde D_1}\le \delta^{2-},$$
and
$$\Vert X_{Q\circ \Phi}\Vert_{p/2,D(\tilde\varrho_1,\tilde\delta_1,4\delta)},\quad \Vert X_{T\circ \Phi}\Vert_{p-1,D(\tilde\varrho_1,\tilde\delta_1,4\delta)}\le \delta^{\mathcal{M}+1},$$
respectively.

Performing this process for $\mathcal{M}-2$ times, under some transformation $\Phi^{**}$  we finally get a Hamiltonian
$$H^{**}=H^{*}\circ\Phi^{**}=N+Z+B+Q+T,$$
with
%with $ \Gamma=\sum_{|\alpha|\le \mathcal{M}+2}\Gamma_\alpha(x)y^\alpha$ and $\Vert\Gamma\Vert\le \delta^{\mathcal{M}}$.
$$N=\sum_{j=1}^n\omega^*_j(\xi)y_j+\sum_{j\geq 1}\frac{[\breve{\Omega}_j](\xi)}{j}|z_j|^2+\sum_{j\geq 1}\frac{\breve{(\Omega}_j-[\breve{\Omega}_j])(x;\xi)}{j}(|z_j|^2-|z_j(0)|^2),$$
\begin{equation}\label{zhuuum}
Z=\sum_{|\alpha|\le \mathcal{M}+2}Z^\alpha(\xi) y^\alpha+\sum_{|\alpha|\geq \mathcal{M}+3}Z^\alpha(x;\xi) y^\alpha,
\end{equation}
$$B=\sum_{|\beta|\geq 1}B^{\beta}(x,y;\xi)(|z|^{2}-|z(0)|^2)^{\beta},$$
$$Q=\sum_{|\mu|+|\nu|\geq3 }Q^{\alpha\beta\gamma\mu\nu}(x;\xi)y^{\alpha}\acute{z}^{\beta}\bar{\acute{z}}^{\gamma}\hat{z}^{\mu}\bar{\hat{z}}^{\nu},$$
%and
%$$R=\sum_{\alpha,\beta,\gamma }R^{\alpha\beta\gamma}(x;\xi)y^{\alpha}z^{\beta}\bar{z}^{\gamma}$$
with the coefficients of $Z$, $B,\ Q$ and $T$ satisfy
\begin{equation}\label{yongz}
\Vert Z^\alpha(x;\xi)\Vert_{D(\varrho)}\le \delta^{2-},\quad |\alpha|\geq\mathcal{M}+3,
\end{equation}
\begin{equation}\label{yongb}
\Vert B^\beta(x,y;\xi)\Vert_{D(\varrho,2\delta)}\le \delta^{2-},\quad |\beta|\geq 1,
\end{equation} and
\begin{equation}\label{zuihout}
\Vert X_Q\Vert_{p/2,D(\varrho,2\delta,4\delta)},\quad \Vert X_T\Vert_{p-1,D(\varrho,2\delta,4\delta)}\le \delta^{\mathcal{M}+1},\end{equation}
respectively.

\medskip

Having the previous preparation, we are now in a position to do the long time stability estimate for the KAM tori obtained in Remark \ref{rem2.7}.

Assume $w(t)=(x(t),y(t),z(t),\overline{z}(t))$ is the solution of Hamiltonian vector field $X_{H^{**}}(x,y,z,\overline{z})$ with the initial datum
$$w(0)=(x(0),y(0),z(0),\bar{z}(0))\in D(\varrho, 2\delta,4\delta)$$
satisfying
$$d_p(w(0),\mathcal{T}_0)\le \delta.$$
Noting $\mathcal{T}_0=\mathbb{T}^n\times\{y=0\}\times\{z=0\}\times\{\bar{z}=0\}$, we have
$$\frac{1}{4\delta}||y(0)||+||z(0)||_{p}+||{\bar{z}}(0)||_{p}\le \delta.$$
The flow on the torus $\mathcal{T}_0$ is $x_0(t)=x_0+\breve{\omega} t$. %Without loss of generality, we assume $x_0=0$.
 Now, we will estimate $||z(t)||_{p/2}$ and $||{\bar{z}}(t)||_{p/2}$. For this purpose, we let
$$\tilde{N}(z,\bar{z}):=||z||_{p/2}^2=\sum_{j\geq 1}|z_j|^2|j|^{2p/2}.$$
By initial data, we have
$$\tilde{N}(w(0))\le \delta^2.$$
Define
$$T^*:=\inf\{|t|: \tilde{N}(w(t))\geq 4\delta^2\}.$$
Then $\tilde{N}(w(t))\le 4\delta^2$ for all $|t|\le T^*$. In particular, we have
 $$\tilde{N}(w(t))= 4\delta^2,\quad{\rm for}\,\, t= T^*\,{\rm or}-T^*.$$
Without loss of generality, we assume
\begin{equation}\label{1w}
\tilde{N}(w(T^*))= 4\delta^2.
\end{equation}
Now we are in a position to show $T^*>\delta^{-\mathcal{M}/2}$. For each $(x,y,z,\bar{z})\in D(\varrho, 2\delta,2\delta)$, from the definition of Poisson bracket under the symplectic structure \eqref{symplectic structure} and \eqref{zuihout}, we have
\begin{align}\label{*}
|\{T, \tilde{N}\}(z,\bar{z})|\le&( \Vert {\bf i}\kappa \partial_zT\Vert_{p-1}+\Vert  -{\bf i}\kappa \partial_{\bar{z}}T\Vert_{p-1})\Vert z\Vert_p\nonumber\\
\le &2\delta ||X_T||_{p-1,D(\varrho, 2\delta,2\delta)}\Vert z\Vert_p\nonumber\\
\le &16\delta^{\mathcal{M}+\frac{5}{2}}.
\end{align}
By \eqref{zuihout}, one can obtain
\begin{align}\label{**}
|\{Q, \tilde{N}\}(x,y,z,\bar{z}))|\le&( \Vert {\bf i}\kappa \partial_zT\Vert_{p/2}+\Vert  -{\bf i}\kappa \partial_{\bar{z}}T\Vert_{p/2})\Vert z\Vert_p\nonumber\\
\le &2\delta ||X_T||_{p/2,D(\varrho, 2\delta,2\delta)}\Vert z\Vert_p\nonumber\\
\le &16\delta^{\mathcal{M}+\frac{5}{2}}.
\end{align}
Thus, one gets
\begin{equation}\label{ntie}
\begin{split}
|\tilde{N}(w(t))-\tilde{N}(w(0))|=&|\int_0^t\{H^{**},\tilde{N}\}(w(s))ds|\\
=&|\int_0^t\{N+Z+B+Q+T,\tilde{N}\}(w(s))ds|\\
\le &\int_0^t|\{Q+T,\tilde{N}\}(w(s))|ds\\
&(N,\ Z\ {\rm and}\ B\ {\rm depend}\ {\rm on}\ |z|^2)\\
\le &32\delta^{\mathcal{M}+\frac{5}{2}}|t|\\
&({\rm by}\ \eqref{*}\ {\rm and}\ \eqref{**})\\
\le &32\delta^{\mathcal{M}+\frac{5}{2}}T^*.
\end{split}
\end{equation}

Proof by contradiction: If $T^*<\delta^{-\mathcal{M}/2}$, then
\begin{align*}
4\delta^2=&|\tilde{N}(w(T))|\le |\tilde{N}(w(0))|+|\tilde{N}(w(t))-\tilde{N}(w(0))|\\
\le &\delta^2+ 32\delta^{(\mathcal{M}+{5})/{2}}\delta^{-\mathcal{M}/2}<2\delta^2,
\end{align*}
which is impossible.

In the following, we will estimate $||y(t)||$. For $1\le j\le n$, we let
$$\tilde{Y}_j(y):=y_j.$$
By the initial data, we have
$$\tilde{Y}(w(0))\le 4\delta^2.$$
Define
$$T_j:=\inf\{|t|: \tilde{Y}_j(w(t))\geq 8\delta^2\}.$$
Then $\tilde{Y}_j(w(t))\le 8\delta^2$ for all $|t|\le T_j$. In particular, we have
 $$\tilde{Y}_j(w(t))= 8\delta^2,\quad{\rm for}\,\,|t|= T_j\,{\rm or}-T_j.$$
Without loss of generality, we assume
\begin{equation}\label{2w}
\tilde{Y}_j(w(T_{j}))= 8\delta^2.
\end{equation}
Now we are in a position to show $T_{j}>\delta^{-\mathcal{M}/4}$.

 %For each $(x,y,z,\bar{z})\in D(4\sigma, 4\delta,4\delta)$, since
% \begin{align*}
%&\sum_{j}|j|^{2(p-1)}||z_j(t)|^2-|z_j(0)|^2|\\
%\le& \sum_{j}|j|^{2(p-1)}\Big|\int_0^t {\bf i} j(Z_{\bar{z}_j}+R_{\bar{z}_j}+Q_{\bar{z}_j})\bar{z}_j-{\bf i} j(Z_{z_j}+R_{z_j}+Q_{z_j})z_jdt\Big|^2\\
%&\le |t| \Vert X_{Z+R+Q}\Vert^2_{p-1} \delta^2\\
%&\le( |t|\delta^{\mathcal{M}+\frac{5}{2}})^2.
%\end{align*}
From \eqref{ntie},  we have $|\sum_{j}|j|^{2p/2}(|z_j(t)|^2-|z_j(0)|^2)|\le  |t|\delta^{(\mathcal{M}+{5})/{2}}$. It follows that
\begin{equation}\label{1q}
 \left|\left\{\sum_{j\geq 1}\frac{(\breve{\Omega}_j-[\breve{\Omega}_j])(x;\xi)}{j}(|z_j|^2-|z_j(0)|^2), \tilde{Y}_j\right\}\right|\le |t|\delta^{(\mathcal{M}+{5})/{2}},
 \end{equation}
 and by \eqref{yongb} and Cauchy's estimate
\begin{align}
&|\{B, \tilde{Y}_j\}(x,y,z,\bar{z})|\le ||\partial_xB||\nonumber\\
\le& || \sum_{|\beta|\geq 1}\partial_x B^\beta(x,y;\xi)(|z|^2-|z(0)|^2)^\beta ||\nonumber\\
\le& \frac{1}{\varrho}\delta^{2-} \sum_{|\beta|\geq 1}  (|t|\delta^{(\mathcal{M}+{5})/{2}})^\beta\nonumber\\
\le & 16\delta^{(\mathcal{M}+{5})/{2}}|t|,\label{2q}
\end{align}
for some $|t|\le T^*$.
 In addition, from \eqref{zhuuum} and \eqref{yongz}, we have
 \begin{align}
&|\{Z, \tilde{Y}_j\}(x,y,z,\bar{z})|\le ||\partial_xZ||\nonumber\\
\le & || \sum_{|\alpha|\geq \mathcal{M}+3}\partial_xZ^\alpha(x;\xi) y^\alpha||\nonumber\\
\le&  \sum_{|\alpha|\geq \mathcal{M}+3}  \frac{1}{\varrho}\delta^{2-}\cdot\delta^{2\mathcal{M}}\nonumber\\
\le& 16\delta^{\mathcal{M}+\frac{5}{2}}.\label{3q}
\end{align}
By \eqref{zuihout}, we obtain
\begin{align}
&|\{T, \tilde{Y}_j\}(x,y,z,\bar{z})|\le ||\partial_xT||\le 16\delta^{\mathcal{M}+\frac{5}{2}},\label{3q}
\end{align}
and
\begin{align}\label{4q}
|\{Q, \tilde{Y}_j\}(x,y,z,\bar{z})|\le ||\partial_xQ||\le 16\delta^{\mathcal{M}+\frac{5}{2}}.
\end{align}
Thus, by \eqref{1q}, \eqref{2q}, \eqref{3q} and \eqref{4q}, one gets
\begin{align*}
|\tilde{Y}_j(w(t))-\tilde{Y}_j(w(0))|=&\left|\int_0^t\{H^{**},\tilde{Y}_j\}(w(s))ds\right|\\
=&\left|\int_0^t\{N+Z+B+Q+T,\tilde{Y}_j\}(w(s))ds\right|\\
%\le &32\delta^{\mathcal{M}+\frac{5}{2}}||t|\\
\le &32\delta^{(\mathcal{M}+{5})/{2}}T_j^2.
\end{align*}
Proof by contradiction: If $T_j<\delta^{-\mathcal{M}/4}$, then
\begin{align*}
8\delta^2=&|\tilde{N}(w(T))|\le |\tilde{N}(w(0))|+|\tilde{N}(w(t))-\tilde{N}(w(0))|\\
\le &4\delta^2+ 64\delta^{(\mathcal{M}+{5})/{2}}\delta^{-\mathcal{M}/2}<5\delta^2,
\end{align*}
which is impossible.

As a consequence, for all $|t|\le\min\{T_*, T_j (1\le j\le n)\}=\delta^{-\mathcal{M}/4}$, from \eqref{1w} and \eqref{2w} we have
$$\frac{1}{4\delta}||w_y(t)||+||w_z(t)||_{p/2}+||w_{\bar{z}}(t)||_{p/2}\le 2\delta.$$
That is,
\begin{equation}\label{etw}
d_{p/2}(w(t),\mathcal{T}_0)\le 2\delta.
\end{equation}
%\begin{Rem}
%Actually, in view of \eqref{} and \eqref{}, the better estimates about $\tilde{N}(w(t))$ and $\tilde{Y}_j(w(t))$ can be obtained:
%$$|\tilde{N}(w(t))-\tilde{N}(w(0))|,\quad |\tilde{Y}_j(w(t))-\tilde{Y}_j(w(0))|\le 32 \delta^{\frac{5}{2}}$$
%\end{Rem}
Based on the partial normal form \eqref{pnf} constructed in Theorem \ref{parth}, for each $\xi\in\mathcal{T}$ of the original Hamiltonian $H$ can be defined by $\mathcal{T}=(\Psi\circ\Phi\circ\Phi^*\circ\Phi^{**})^{-1}\mathcal{T}_0$. Assume $\tilde{w}(t)$ is a solution to the original Hamiltonian vector field $X_H$ with the initial datum $\tilde{w}(0)=(\tilde{w}_x(0),\tilde{w}_y(0),\tilde{w}_z(0),\tilde{w}_{\bar{z}}(0))$ satisfying
$$d_p(\tilde{w}(0),\mathcal{T})\le \delta.$$
Then there exists $\tilde{w}^*\in\mathcal{T}$ such that
$$d_p(\tilde{w}(0),\tilde{w}^*)\le \frac{9}{8}\delta.$$
Hence,
\begin{align*}
&d_{p}(\Psi\circ\Phi\circ\Phi^*\circ\Phi^{**}\circ \tilde{w}(0),\mathcal{T}_0)\\
\le & d_{p}(\Psi\circ\Phi\circ\Phi^*\circ\Phi^{**}\circ \tilde{w}(0),\Psi\circ\Phi\circ\Phi^*\circ\Phi^{**}\circ \tilde{w}^*)\\
\le & d_{p}(\Psi\circ\Phi\circ\Phi^*\circ\Phi^{**}\circ \tilde{w}(0),\tilde{w}(0))+d_{p}(\tilde{w}(0), \tilde{w}^*)+d_{p}(\tilde{w}^*,\Psi\circ\Phi\circ\Phi^*\circ\Phi^{**}\circ \tilde{w}^*)\\
\le &4\delta ||\Psi\circ\Phi\circ\Phi^*\circ\Phi^{**}\circ \tilde{w}(0)-\tilde{w}(0)||_{p, D(\varrho, 2\rho,4\rho)}+\frac{9}{8}\delta\\
&+4\delta||\tilde{w}^*,\Psi\circ\Phi\circ\Phi^*\circ\Phi^{**}\circ \tilde{w}^*||_{p, D(\varrho, 2\rho,4\rho)}\\
\le & \frac{4}{3}\delta,
\end{align*}
where the last inequality comes from \eqref{lho1} and \eqref{lho2}.
Therefore, according to the estimate \eqref{etw}, we have that, for $w(t)=\Psi\circ\Phi\circ\Phi^*\circ\Phi^{**}\circ \tilde{w}(t)$,
$$d_{p/2}(\Psi\circ\Phi\circ\Phi^*\circ\Phi^{**}\circ \tilde{w}(t),\mathcal{T}_0)\le  \frac{5}{3}\delta,\quad {\rm for}\ \ {\rm all}\ \ |t|<\delta^{-\mathcal{M}/4}.$$
Moreover, we have
$$d_{p/2}(\tilde{w}(t),(\Psi\circ\Phi\circ\Phi^*\circ\Phi^{**})^{-1}\mathcal{T}_0)\le 2\delta,\quad {\rm for}\ \ {\rm all}\ \ |t|<\delta^{-\mathcal{M}/4}.$$

\subsection{Proof of main Theorem \ref{theorem1.1}}
By subsection \ref{2.1.1}, we obtain a Hamiltonian $H(x,y, z,\bar{z};\xi)$ of the following form
$$H(x,y, z,\bar{z};\xi)=N(y, z,\bar{z};\xi)+ P(x,y, z,\bar{z};\xi),$$
where
$$
N(y, z,\bar{z};\xi)=\sum_{i=1}^n\omega_i(\xi)y_i+\sum_{j\in\mathbb{Z}_*}\frac{\Omega_j(\xi)}{j}z_j\bar{z}_j,
$$
with the tangent frequency
\begin{equation}\label{oi}
\omega(\xi)=(\omega_{1}(\xi),\cdots,\omega_{n}(\xi))=\left(j_1^2+\xi_{j_1},\cdots,j_n^2+\xi_{j_n}\right),
\end{equation}
and the normal frequency $\left(\frac{\Omega_{j}(\xi)}{j}\right)_{j\in\mathbb{Z}_*}$ where
\begin{equation}\label{Oi}
(\Omega_{j}(\xi))_{j\in\mathbb{Z}_*}=\left(j^2+\xi_{j}\right)_{j\in\mathbb{Z}_*}.
\end{equation}
The perturbation becomes $ P(x,y,z,\bar{z})$ which is independent of $\xi$.

In view of \eqref{oi} and \eqref{Oi}, it is easy to verify the {\bf Assumption A, B} and {\bf C} in Theorem \ref{normaltheorem} hold. Moreover, the perturbation $ P(x,y,z,\bar{z})$ satisfies $|||P|||_{D(s_0,r_0,r_0)\times \Pi}<\epsilon$. Hence, all assumptions in Theorem \ref{normaltheorem} hold. Therefore, there exists a subset $\Pi_\eta\subset \Pi$ with the estimate
$${\rm Meas}\,\Pi_\eta \geq {\rm Meas}\,\Pi (1-O(\eta)).$$
For each $\xi \in \Pi_\eta$, there is symplectic transformation
$$\Psi: D(s_0/2,r_0/2,r_0/2)\rightarrow D(s_0,r_0,r_0),$$
such that
$$\breve{H}(x,y,z,\bar{z};\xi)=\breve{N}(x,y,z,\bar{z};\xi)+\breve{P}(x,y,z,\bar{z};\xi),$$
where
$$\breve{N}(y,z,\bar{z};\xi)=\sum_{j=1}^n\breve{\omega}_j(\xi)y_j+\sum_{j\geq 1}\breve{\Omega}_j(x;\xi)z_j\bar{z}_j$$
$$=\sum_{j=1}^n\breve{\omega}_j(\xi)y_j+\sum_{|j|\le \mathcal{N}}\breve{\Omega}_j(\xi)z_j\bar{z}_j+\sum_{|j|\geq\mathcal{N}+ 1}\breve{\Omega}_j(x;\xi)z_j\bar{z}_j,$$
and
$$\breve{P}(x,y,z,\bar{z};\xi)=\sum_{\substack{\alpha\in\mathbb{N}^n,\beta,\gamma\in \mathbb{N}^{\mathbb{N}},\\ 2|\alpha|+|\beta+\gamma|\geq 3}}\breve{P}^{\alpha\beta\gamma}(x;\xi)y^{\alpha}z^{\beta}\bar{z}^{\gamma}.$$
Based on Theorem \ref{lthe}, there is a subset ${\Pi}_{\acute{\eta}}$ such that for each $\xi \in{\Pi}_{\acute{\eta}}$ equation \eqref{main equation} possess a KAM tori $\mathcal{T}_\xi$ in $H_0^{p}([0,\pi])$. In addition, for any solution $u(t,x)$ of equation \eqref{main equation} with the initial datum satisfying
$$d_{H_0^{p}([0,\pi])}(u(0,x), \mathcal{T}_\xi)\le \delta,$$
then
$$d_{H_0^{p/2}([0,\pi])}(u(t,x), \mathcal{T}_\xi)\le2 \delta,\quad {\rm for}\,\, {\rm all}\,\, |t|\le \delta^{-\mathcal{M}/4}.$$

\section{Some technical lemmata}

\begin{Lem}\label{lewxg}
Consider a function $W(x;\xi)=\sum_{k\in\mathbb{Z}^n}W(k;\xi)e^{{\bf i}\langle k,x\rangle}$ defined on the domain $D(s)\times \Pi$, which is analytic in $(x)\in D(s)$ and $C^1$ in $\xi\in\Pi$, then the following generalized Cauchy estimate holds
$$\Vert \partial_xW\Vert_{D(s-\sigma,r)\times\Pi}\le\frac{1}{e\sigma}\Vert W\Vert_{D(s,r)\times\Pi},$$
%$$\Vert W_y\Vert_{D(s,r-\sigma^\prime)\times\Pi}\le\frac{1}{r\sigma^\prime}\Vert W\Vert_{D(s,r)\times\Pi},$$
where $0<\sigma<s$. %and $0<\sigma^\prime<r/2$.
\end{Lem}

This proof of this Lemma can be found in Lemma 7.3 in \cite{CLY}.

\begin{Lem}\label{brak}
Consider two functions
\begin{equation*}
  \begin{split}
    U(x,y,z,\bar{z};\xi)&=\sum_{\kappa\in \mathbb{N}^n,\alpha,\beta\in\mathbb{N}^{\mathbb{Z}_*}} U^{\kappa\alpha\beta}(x;\xi)y^\kappa z^\alpha\bar{z}^\beta\in\mathcal{H}\\
      & =\sum_{k\in\mathbb{Z}^n,\kappa\in \mathbb{N}^n,\alpha,\beta\in\mathbb{N}^{\mathbb{Z}_*}} U^{\kappa\alpha\beta}(k;\xi)e^{{\bf i}\langle k,x\rangle}y^\kappa z^\alpha\bar{z}^\beta,
  \end{split}
\end{equation*}

and
$$ V(x,y,z,\bar{z};\xi)=\sum_{\substack{\tau\in \mathbb{N}^n,\mu,\gamma\in\mathbb{N}^{\mathbb{Z}_*}\\2|\tau|+|\mu|+|\gamma|\le 2}} V^{\tau\mu\gamma}(x;\xi)y^\tau z^\mu\bar{z}^\gamma\in\mathcal{H}_*$$
$$=\sum_{\substack{|k|\le K, \tau\in \mathbb{N}^n,\mu,\gamma\in\mathbb{N}^{\mathbb{Z}_*}\\2|\tau|+|\mu|+|\gamma|\le 2}} V^{\tau\mu\gamma}(k;\xi)e^{{\bf i}\langle k,x\rangle}y^\tau z^\mu\bar{z}^\gamma,$$
where the functions $U,  V$ defined on the domain $D(s,r,r)\times \Pi$ are analytic in $(x,y,z,\bar{z})\in D(s,r,r)$ and $C^1$ in $\xi \in \Pi$. %In addition, $U$ and $V$ satisfy for real variables $x$
%\begin{equation}\label{realim}
%U^{\kappa\alpha\beta}(x;\xi)=U^{\kappa\beta\alpha}(x;\xi)\in\mathbb{R},\quad V^{\tau\mu\gamma}(x;\xi)=-V^{\tau\gamma\mu}(x;\xi)\in {\bf i}\mathbb{R}.
%\end{equation}
Then we have that $\{U,V\}$ belongs to $\mathcal{H}$ with the following estimate
%$$h\geq 0, 2> g\geq 0:\quad ||| \{U,V\}|||_{D(s-\sigma)\times \Pi}\le \frac{(2K)^n+n}{e\sigma} ||| U|||_{D(s)\times \Pi}||| V|||^*_{D(s)\times \Pi}.$$
$$ ||| \{U,V\}|||_{D(s-\sigma,r-\sigma^\prime,r-\sigma^\prime)\times \Pi}\le C\max\Big\{\frac{K^n}{r\sigma^\prime},\frac{1}{r^2\sigma}\Big\} ||| U|||_{D(s,r,r)\times \Pi}||| V|||^*_{D(s,r,r)\times \Pi}.$$
%Moreover, for real variables $x$, the coefficients $\{U,V\}^{q^\prime l^\prime k^\prime}$ satisfy
%$$\{U,V\}^{q^\prime l^\prime k^\prime}=\{U,V\}^{q^\prime k^\prime l^\prime }\in \mathbb{R}.$$
\end{Lem}

\begin{proof}
%We will consider the case $h>3, g\geq 2$, the other cases are easier to estimate, here we omit the proof.
We consider the functions
$$ U(x,y,z,\bar{z};\xi)=\sum_{2|\kappa|+|\alpha+\beta|=h} U^{\kappa\alpha\beta}(x;\xi)y^\kappa z^\alpha\bar{z}^\beta\in\mathcal{H}$$
$$=\sum_{k\in\mathbb{Z}^n,2|\kappa|+|\alpha+\beta|=h} U^{\kappa\alpha\beta}(k;\xi)e^{{\bf i}\langle k,x\rangle}y^\kappa z^\alpha\bar{z}^\beta,$$
and
$$ V(x,y,z,\bar{z};\xi)=\sum_{2|\tau|+|\mu+\gamma|=2} V^{\tau\mu\gamma}(x;\xi)y^\tau z^\mu\bar{z}^\gamma\in\mathcal{H}_*,$$
$$=\sum_{|k|\le K, 2|\tau|+|\mu+\gamma|=2} V^{\tau\mu\gamma}(k;\xi)e^{{\bf i}\langle k,x\rangle}y^\tau z^\mu\bar{z}^\gamma.$$
By the definition of Possion bracket, we have
$$\{U,V\}=\sum_{2|q^\prime|+|l^\prime+k^\prime|=h}\{U,V\}^{q^\prime l^\prime k^\prime}y^{q^\prime}z^{l^\prime}\bar{z}^{k^\prime},$$
where, for $2|q^\prime|+|l^\prime+k^\prime|=h$,
$$\{U,V\}^{q^\prime l^\prime k^\prime}=\{U,V\}^{q^\prime l^\prime k^\prime}_1+\{U,V\}^{q^\prime l^\prime k^\prime}_2,$$
with
\begin{align}
\{U,V\}^{q^\prime l^\prime k^\prime}_1=&\sum_{\substack{\kappa+\tau-e_j=q^{\prime}, \alpha+\mu=l^\prime,\beta+\gamma=k^\prime\\2|\kappa|+|\alpha+\beta|=h,2|\tau|+|\mu+\gamma|=2}}\tau_j\partial_{x_j}U^{\kappa\alpha\beta} V^{\tau\mu\gamma}\label{lpkp2}\\
&+\sum_{\substack{\kappa-e_j+\tau=q^{\prime}, \alpha+\mu=l^\prime,\beta+\gamma=k^\prime\\2|\kappa|+|\alpha+\beta|=h,2|\tau|+|\mu+\gamma|=2}}\kappa_jU^{\kappa\alpha\beta} \partial_{x_j}V^{\tau\mu\gamma},\label{lpkp}
\end{align}
and
\begin{align}
\{U,V\}^{q^\prime l^\prime k^\prime}_2=&{\bf i}\sum_{j\in\mathbb{Z}^*}\Bigg(\sum_{\substack{\tau+\kappa=q^\prime,(\alpha-e_j)+\mu=l^\prime,\beta+(\gamma-e_j)=k^\prime\\2|\kappa|+|\alpha+\beta|=h,2|\tau|+|\mu+\gamma|=2}} j\alpha_j\gamma_jU^{\kappa\alpha\beta}V^{\tau\mu\gamma}\nonumber\\
&-
\sum_{\substack{\tau+\kappa=q^\prime,\alpha+(\mu-e_j)=l^\prime,(\beta-e_j)+\gamma=k^\prime\\2|\kappa|+|\alpha+\beta|=h,2|\tau|+|\mu+\gamma|=2}} j\beta_j\mu_jU^{\kappa\alpha\beta}V^{\tau\mu\gamma}\Bigg).\label{lpkp1}
\end{align}
Consider the term $\{U,V\}^{q^\prime l^\prime k^\prime}_2$. In what follows, we will show there are at most
\begin{equation}\label{shushu}
2(2K)^nh
\end{equation}
terms at the right hand of \eqref{lpkp1}.

Consider the nonzero components of vector $q^\prime, k^\prime$ and $l^\prime$. Taking multiplicity into account, we have that $(q^\prime, l^\prime, k^\prime)$ has $|q^\prime|+| l^\prime|+| k^\prime|$ non zero components whose values are $1$. We denote them by
\begin{align*}
q_{j_{1}},\cdots, q_{j_{|q^\prime|}}, l_{j_{|q^\prime|+1}},\cdots,l_{j_{|q^\prime|+|l^\prime|}},k_{j_{|q^\prime|+|l^\prime|+1}},\cdots, k_{j_{|q^\prime|+|l^\prime|+|k^\prime|}}.
\end{align*}
If $(\tau,\mu,\gamma-e_j)$ is fixed, then $(\kappa,\alpha-e_j,\beta)$ can be determined by
\begin{equation}\label{ss11}
\kappa+\tau=q^\prime,\quad (\alpha-e_j)+\mu=l^\prime,\quad \beta+(\gamma-e_j)=k^\prime,
\end{equation}
or
\begin{equation}\label{ss22}
\kappa+\tau=q^\prime,\quad \alpha+(\mu-e_j)=l^\prime,\quad(\beta-e_j)+\gamma=k^\prime.
\end{equation}
Therefore we choose $|\tau|+|\mu|+|\gamma|-1$ positions from $|q^\prime|+| l^\prime|+| k^\prime|$ positions. Actually, there are $\left(\begin{array}{l}
q^\prime|+| l^\prime|+| k^\prime|\\
|\tau|+|\mu|+|\gamma|-1
\end{array}
\right)$ choices.
If  $(\tau,\mu,\gamma-e_j)$ is fixed, the last position $j$ i.e., $e_j$ will have $(2K)^n$ choices, according to the momentum conservation:
$$-\sum_{i=1}^nk_ij_i+\sum_{j\in \mathbb{Z}_*}j(\mu_j-\gamma_j)=0.$$
In conclusion, the number of pairs $(\tau,\mu,\gamma)$ and $(\kappa,\alpha,\beta)$ satisfying \eqref{ss11} and \eqref{ss22} equals
\begin{align}
2(2K)^n\left(\begin{array}{l}
|q^\prime|+| l^\prime|+| k^\prime|\\
|\tau|+|\mu|+|\gamma|-1
\end{array}
\right)\le& 2(2K)^n\left(\begin{array}{l}
|q^\prime|+| l^\prime|+| k^\prime|+|\tau|\\
2|\tau|+|\mu|+|\gamma|-1
\end{array}
\right)\nonumber\\
=& 2(2K)^n\left(\begin{array}{l}
|q^\prime|+| l^\prime|+| k^\prime|+|\tau|\\
1
\end{array}
\right)\nonumber\\
\le & 2(2K)^n\left(\begin{array}{l}
h\\
1
\end{array}
\right)=2(2K)^nh,\label{zzza}
\end{align}
where we use $\left(\begin{array}{l}
a\\
b
\end{array}
\right)\le \left(\begin{array}{l}
a+k\\
b+k
\end{array}
\right)$ and  $\left(\begin{array}{l}
a\\
b
\end{array}
\right)\le \left(\begin{array}{l}
a+l\\
b
\end{array}
\right)$.

Similarly, one can prove there is at most
\begin{equation}\label{shushu}
nh
\end{equation}
terms at the right hand of \eqref{lpkp}. In addition, by the special form of  $V$, it is easy to find that there are $n$ terms in \eqref{lpkp2}. Therefore, we have
$$\Vert \eqref{lpkp2}\Vert \le n \sup_{\kappa,\alpha,\beta}\Vert \partial_{x}U^{\kappa\alpha\beta}\Vert \sup_{\tau,\mu,\gamma} \Vert V^{\tau\mu\gamma}\Vert$$
$$\le \frac{n}{e\sigma} \sup_{\kappa,\alpha,\beta}\Vert U^{\kappa\alpha\beta}\Vert \sup_{\tau,\mu,\gamma}\Vert V^{\tau\mu\gamma}\Vert$$
and
$$\Vert \eqref{lpkp}\Vert \le n h\sup_{\kappa,\alpha,\beta}\Vert U^{\kappa\alpha\beta}\Vert \sup_{\tau,\mu,\gamma} \Vert \partial_{x}V^{\tau\mu\gamma}\Vert$$
$$\le n hK\sup_{\kappa,\alpha,\beta}\Vert U^{\kappa\alpha\beta}\Vert \sup_{\tau,\mu,\gamma}\Vert V^{\tau\mu\gamma}\Vert,$$
since $V$ is a truncated Fourier series of $K$-th order.

From \eqref{zzza} and \eqref{shushu}, we obtain
\begin{equation*}
  \begin{split}
    \Vert\{U,V\}^{q^\prime l^\prime k^\prime}\Vert_{D(s-\sigma)\times \Pi} &\le \max\Big\{\big(2(2K)^n+nK\big)h,\frac{1}{e\sigma}n\Big\} \\
      & \cdot \sup_{\kappa,\alpha,\beta}\Vert U^{\kappa\alpha\beta}\Vert_{D(s)\times \Pi}M_{\mu\gamma}\sup_{\tau,\mu,\gamma}\Vert V^{\tau\mu\gamma}\Vert_{D(s)\times \Pi}.
  \end{split}
\end{equation*}
%Following the proof of Proposition 2.5 in \cite{MR4420616}, one can obtain for the term $\{U,V\}^{q^\prime l^\prime k^\prime}_2$
%$$\Vert\{U,V\}^{l^\prime k^\prime}\Vert_{D(s-\sigma)\times \Pi}\le 2^{h+g} \sup_{\kappa,\alpha,\beta}\Vert U^{\kappa\alpha\beta}\Vert_{D(s)\times \Pi} M_{\mu\gamma}\sup_{\tau,\mu,\gamma}\Vert V^{\tau\mu\gamma}\Vert_{D(s)\times \Pi}.$$
%That is
%$$\Vert\{U,V\}^{q^\prime l^\prime k^\prime}\Vert_{D(s-\sigma)\times \Pi}$$
%$$\le C\max\Big\{K^nh,\frac{1}{\sigma}\Big\}\sup_{\kappa,\alpha,\beta}\Vert U^{\kappa\alpha\beta}\Vert_{D(s)\times \Pi}\cdot M_{\mu\gamma}\sup_{\tau,\mu,\gamma}\Vert V^{\tau\mu\gamma}\Vert_{D(s)\times \Pi},$$
%where the constant $C$ only depends on $n$.
Then, we have
\begin{align*}
&(r-\sigma^\prime)^{h}\Vert\{U,V\}^{q^\prime l^\prime k^\prime}\Vert_{D(s-\sigma)\times \Pi}\\
\le& C\max\Big\{K^nh,\frac{1}{\sigma}\Big\}\sup_{\kappa,\alpha,\beta}\Vert U^{\kappa\alpha\beta}\Vert_{D(s)\times \Pi}\cdot M_{\mu\gamma}\sup_{\tau,\mu,\gamma}\Vert V^{\tau\mu\gamma}\Vert_{D(s)\times \Pi} (r-\sigma^\prime)^{h}\\
\le & C\max\Big\{\frac{K^n}{\sigma^\prime},\frac{1}{r\sigma}\Big\}\sup_{\kappa,\alpha,\beta}\Vert U^{\kappa\alpha\beta}\Vert_{D(s)\times \Pi}\cdot M_{\mu\gamma}\sup_{\tau,\mu,\gamma}\Vert V^{\tau\mu\gamma}\Vert_{D(s)\times \Pi} r^{h+1}\\
\le &   C\max\Big\{\frac{K^n}{r\sigma^\prime},\frac{1}{r^2\sigma}\Big\}\sup_{\kappa,\alpha,\beta}\Vert U^{\kappa\alpha\beta}\Vert_{D(s)\times \Pi} r^{h}\cdot M_{\mu\gamma}\sup_{\tau,\mu,\gamma}\Vert V^{\tau\mu\gamma}\Vert_{D(s)\times \Pi} r^{2},
\end{align*}
where the constant $C$ only depends on $n$.
The second inequality comes from $h(r-\sigma^\prime)^{h-1}\sigma^\prime\le r^{h}$.
It follows that:
\begin{equation}\label{simuv}
||| \{U,V\}|||_{D(s-\sigma,r-\sigma^\prime,r-\sigma^\prime)\times \Pi}\le C\max\Big\{\frac{K^n}{r\sigma^\prime},\frac{1}{r^2\sigma}\Big\} ||| U|||_{D(s,r,r)\times \Pi}||| V|||^*_{D(s,r,r)\times \Pi}.
\end{equation}
Similarly, we can prove that \eqref{simuv} still holds after replacing $U, V$ by
$$ U(x,y,z,\bar{z};\xi)=\sum_{2|\kappa|+|\alpha+\beta|=h} U^{\kappa\alpha\beta}(x;\xi)y^\kappa z^\alpha\bar{z}^\beta\in\mathcal{H}, $$
and
$$V(x,y,z,\bar{z};\xi)=\sum_{2|\tau|+|\mu+\gamma|=g} V^{\tau\mu\gamma}(x;\xi)y^\tau z^\mu\bar{z}^\gamma\in\mathcal{H}_*,\quad g=0,1.$$

Finally, we consider more general functions
$$U(x,y,z,\bar{z};\xi)=\sum_{h\geq 0}U_h(x,y,z,\bar{z};\xi),$$
 $$V(x,y,z,\bar{z};\xi)=\sum_{g\geq 0}U_g(x,y,z,\bar{z};\xi),$$
where
$$ U_h(x,y,z,\bar{z};\xi)=\sum_{\substack{\alpha\in \mathbb{N}^n,\mu,\gamma\in\mathbb{N}^{\mathbb{Z}_*}\\2|\alpha|+|\mu+\gamma|=h}} U^{\alpha\mu\gamma}(x;\xi)y^\alpha z^\mu\bar{z}^\gamma,$$
$$ V_g(x,y,z,\bar{z};\xi)=\sum_{\substack{\alpha\in \mathbb{N}^n,\mu,\gamma\in\mathbb{N}^{\mathbb{Z}_*}\\2|\alpha|+|\mu+\gamma|=g}} V^{\alpha\mu\gamma}(x;\xi)y^\alpha z^\mu\bar{z}^\gamma.$$
By a direct calculation, we have
\begin{align*}
&|||\{U,V\}|||_{D(s-\sigma,r-\sigma^\prime,r-\sigma^\prime)\times \Pi}\\
=&|||\{\sum_{h\geq 0}U_h,\sum_{g\geq 0}V_g\}|||_{D(s-\sigma,r-\sigma^\prime,r-\sigma^\prime)\times \Pi}\\
\le &\sum_{h,g\geq 0}|||\{U_h,V_g\}|||_{D(s-\sigma,r-\sigma^\prime,r-\sigma^\prime)\times \Pi}\\
\le & \sum_{h,g\geq 0}CK^n\frac{1}{r\sigma\sigma^\prime} ||| U_h|||_{D(s,r,r)\times \Pi}||| V_g|||^*_{D(s,r,r)\times \Pi}\\
\le& C\max\Big\{\frac{K^n}{r\sigma^\prime},\frac{1}{r^2\sigma}\Big\}\left( \sum_{h\geq 0}||| U_h|||_{D(s,r,r)\times \Pi}\right)\left(\sum_{g\geq 0}||| V_g|||_{D(s,r,r)\times \Pi}\right)\\
\le&C\max\Big\{\frac{K^n}{r\sigma^\prime},\frac{1}{r^2\sigma}\Big\}||| U|||_{D(s,r,r)\times \Pi}||| V|||^*_{D(s,r,r)\times \Pi}.
\end{align*}
Thus, we finish the proof of this Lemma.
\end{proof}

\begin{Lem}\label{lemma5.3}
Consider two Hamiltonians $U(x,y,z,\bar{z};\xi)\in \mathcal{H}$ and $V(x,y,z,\bar{z};\xi)\in\mathcal{H}^*$ defined on the domain $D(s,r,r)\times \Pi$ for some $0<s,r\le 1$. In addition, we assume $U,V$ satisfies the condition in Lemma \ref{brak}. Given $0<\sigma<s, 0<\sigma^\prime<r$. Suppose
$$|||V|||^*_{D(s,r,r)\times \Pi}\le \frac{1}{2A},$$
where
$$A=C_0e\max\Big\{\frac{K^n}{r\sigma^\prime},\frac{1}{r^2\sigma}\Big\},$$
and $C_0$ is the constant given in Lemma \ref{brak}. Then for each $|t|\le 1$, we have
$$|||U\circ X_V^t|||_{D(s-\sigma,r-\sigma^\prime,r-\sigma^\prime)\times \Pi}\le 2 ||| U|||_{{D(s,r,r)\times \Pi}}.$$
\end{Lem}

\begin{proof}
Let $$W^{(j)}=\{\cdots\{U,V\},V\},\underbrace{V\}\cdots\},V}_{j-fold}\}.$$

For $j\geq 1$, let $\sigma_j=\frac{\sigma}{j}$ and $\sigma_j^\prime=\frac{\sigma_j^\prime}{j}$. Hence, we obtain
\begin{align*}
&|||W^{(j)}|||_{D(s-\sigma,r-\sigma^\prime,r-\sigma^\prime)\times \Pi}\\
=&|||W^{(j)}|||_{D(s-j\sigma_j,r-j\sigma_j^\prime,r-j\sigma_j^\prime)\times \Pi}\\
\le & \left(C_0\max\Big\{\frac{K^n}{r\sigma_j^\prime},\frac{1}{r^2\sigma_j}\Big\}\right)^j\left(||| U|||_{_{D(s,r,r)\times \Pi}}\right)\left(||| V|||_{_{D(s,r,r)\times \Pi}}\right)^j\\
\le & j^{j}\left(C_0\max\Big\{\frac{K^n}{r\sigma^\prime},\frac{1}{r^2\sigma}\Big\}\right)^j\left(||| U|||_{_{D(s,r,r)\times \Pi}}\right)\left(||| V|||_{_{D(s,r,r)\times \Pi}}\right)^j.
\end{align*}
Using the inequality
$$j^j\le j!e^j,$$
one obtains
\begin{align*}
&\frac{1}{j!}|||W^{(j)}|||_{D(s-\sigma,r-\sigma^\prime,r-\sigma^\prime)\times \Pi}\\
\le & e^{j}\left(C_0\max\Big\{\frac{K^n}{r\sigma^\prime},\frac{1}{r^2\sigma}\Big\}\right)^j\left(||| U|||_{_{D(s,r,r)\times \Pi}}\right)\left(||| V|||_{_{D(s,r,r)\times \Pi}}\right)^j\\
=&\left(||| U|||_{_{D(s,r,r)\times \Pi}}\right)\left(A||| V|||_{_{D(s,r,r)\times \Pi}}\right)^j\\
\le &\frac{1}{2^j}||| U|||_{_{D(s,r,r)\times \Pi}}.
\end{align*}
Expanding the Hamiltonian $U\circ X_V^t$ into Taylor series about $t$ at $t=0$, we get
\begin{align*}
&|||U\circ X_V^t|||_{D(s-\sigma,r-\sigma^\prime,r-\sigma^\prime)\times \Pi}\\
=&|||\sum_{j\geq 0}\frac{t^j}{j!}W^{(j)}|||_{D(s-\sigma,r-\sigma^\prime,r-\sigma^\prime)\times \Pi}\\
\le &\sum_{j\geq 0} \frac{1}{j!}|||W^{(j)}|||_{D(s-\sigma,r-\sigma^\prime,r-\sigma^\prime)\times \Pi}\\
\le & \sum_{j\geq 0} \frac{1}{2^j}||| U|||_{_{D(s,r,r)\times \Pi}}\\
\le & 2||| U|||_{_{D(s,r,r)\times \Pi}}.
\end{align*}

\end{proof}

\begin{Lem}\label{brak111}
Consider two functions
$$ U(x,y,z,\bar{z};\xi)=\sum_{\substack{|k|\in \mathbb{Z}^n,\alpha\in \mathbb{N}^n,\mu,\gamma\in\mathbb{N}^{\mathbb{Z}_*}\\2|\kappa|+|\alpha+\beta|=h}} U^{\kappa\alpha\beta}(k;\xi)e^{{\bf i}\langle k,x\rangle}y^\kappa z^\alpha\bar{z}^\beta,$$
and
$$ V(x,y,z,\bar{z};\xi)=\sum_{\substack{|k|\le K, \tau\in \mathbb{N}^n,\mu,\gamma\in\mathbb{N}^{\mathbb{Z}_*}\\2|\tau|+|\mu|+|\gamma|=g}} V^{\tau\mu\gamma}(k;\xi)e^{{\bf i}\langle k,x\rangle}y^\tau z^\mu\bar{z}^\gamma,$$
where the functions $U,  V$ defined on the domain $D(s,r,r)\times \Pi$ are analytic in $(x,y,z,\bar{z})\in D(s,r,r)$ and $C^1$ in $\xi \in \Pi$. %In addition, $U$ and $V$ satisfy for real variables $x$
%\begin{equation}\label{realim}
%U^{\kappa\alpha\beta}(x;\xi)=U^{\kappa\beta\alpha}(x;\xi)\in\mathbb{R},\quad V^{\tau\mu\gamma}(x;\xi)=-V^{\tau\gamma\mu}(x;\xi)\in {\bf i}\mathbb{R}.
%\end{equation}
We assume $V$ has at most degree $2$ in $(z_j,\bar{z}_j)_{|j|>N}$.
Then we have that $\{U,V\}$ belongs to $\mathcal{H}$ with the following estimate
%$$h\geq 0, 2> g\geq 0:\quad ||| \{U,V\}|||_{D(s-\sigma)\times \Pi}\le \frac{(2K)^n+n}{e\sigma} ||| U|||_{D(s)\times \Pi}||| V|||^*_{D(s)\times \Pi}.$$
$$ ||| \{U,V\}|||_{D(s-\sigma,r-\sigma^\prime,r-\sigma^\prime)\times \Pi}\le C\frac{K^{n+1} (2N+n)^{g-2}}{\sigma\sigma^\prime r}||| U|||_{D(s-\sigma,r,r)\times \Pi} ||| V|||^*_{D(s,r,r)\times \Pi}.$$
%In addition, one has
%$$ ||| \{U,V\}|||_{D(s-\sigma,r-\sigma^\prime,r-\sigma^\prime)\times \Pi}\le C\frac{K^{n+1} (2N+h)^{g-2}}{(s_0-s+\sigma)\sigma^\prime r} ||| U|||_{D(s_0,r,r)\times \Pi}||| V|||^*_{D(s,r,r)\times \Pi}.$$

%Moreover, for real variables $x$, the coefficients $\{U,V\}^{q^\prime l^\prime k^\prime}$ satisfy
%$$\{U,V\}^{q^\prime l^\prime k^\prime}=\{U,V\}^{q^\prime k^\prime l^\prime }\in \mathbb{R}.$$
\end{Lem}
\begin{Rem}
The proof of this lemma follows from the proof of Lemma 4.6 of \cite{yuan2014long}.
\end{Rem}
\begin{proof}
By the definition of Possion bracket, we have
$$\{U,V\}=\sum_{2|q^\prime|+|l^\prime+k^\prime|=h+g-2}\{U,V\}^{q^\prime l^\prime k^\prime}y^{q^\prime}z^{l^\prime}\bar{z}^{k^\prime},$$
where, for $2|q^\prime|+|l^\prime+k^\prime|=h+g-2$,
$$\{U,V\}^{q^\prime l^\prime k^\prime}=\{U,V\}^{q^\prime l^\prime k^\prime}_1+\{U,V\}^{q^\prime l^\prime k^\prime}_2,$$
with
\begin{align}
\{U,V\}^{q^\prime l^\prime k^\prime}_1=&\sum_{\substack{\kappa+\tau-e_j=q^{\prime}, \alpha+\mu=l^\prime,\beta+\gamma=k^\prime\\2|\kappa|+|\alpha+\beta|=h,2|\tau|+|\mu+\gamma|=g}}\tau_j\partial_{x_j}U^{\kappa\alpha\beta} V^{\tau\mu\gamma}\label{qu1}\\
&+\sum_{\substack{\kappa-e_j+\tau=q^{\prime}, \alpha+\mu=l^\prime,\beta+\gamma=k^\prime\\2|\kappa|+|\alpha+\beta|=h,2|\tau|+|\mu+\gamma|=g}}\kappa_jU^{\kappa\alpha\beta} \partial_{x_j}V^{\tau\mu\gamma}\label{qu2},
\end{align}
and
\begin{align}
\{U,V\}^{q^\prime l^\prime k^\prime}_2=&{\bf i}\sum_{j\in\mathbb{Z}^*}\Bigg(\sum_{\substack{\tau+\kappa=q^\prime,(\alpha-e_j)+\mu=l^\prime,\beta+(\gamma-e_j)=k^\prime\\2|\kappa|+|\alpha+\beta|=h,2|\tau|+|\mu+\gamma|=g}} j\alpha_j\gamma_jU^{\kappa\alpha\beta}V^{\tau\mu\gamma}\nonumber\\
&-
\sum_{\substack{\tau+\kappa=q^\prime,\alpha+(\mu-e_j)=l^\prime,(\beta-e_j)+\gamma=k^\prime\\2|\kappa|+|\alpha+\beta|=h,2|\tau|+|\mu+\gamma|=g}} j\beta_j\mu_jU^{\kappa\alpha\beta}V^{\tau\mu\gamma}\Bigg).
\end{align}
In addition, we write
$$\{U,V\}^{q^\prime l^\prime k^\prime}_2:=\{U,V\}^{q^\prime l^\prime k^\prime}_{21}+\{U,V\}^{q^\prime l^\prime k^\prime}_{22},$$
with
\begin{align*}
\{U,V\}^{q^\prime l^\prime k^\prime}_{21}=&{\bf i}\sum_{0<|j|\le N}\Bigg(\sum_{\substack{\tau+\kappa=q^\prime,(\alpha-e_j)+\mu=l^\prime,\beta+(\gamma-e_j)=k^\prime\\2|\kappa|+|\alpha+\beta|=h,2|\tau|+|\mu+\gamma|=g}} j\alpha_j\gamma_jU^{\kappa\alpha\beta}V^{\tau\mu\gamma}\nonumber\\
&-
\sum_{\substack{\tau+\kappa=q^\prime,\alpha+(\mu-e_j)=l^\prime,(\beta-e_j)+\gamma=k^\prime\\2|\kappa|+|\alpha+\beta|=h,2|\tau|+|\mu+\gamma|=g}} j\beta_j\mu_jU^{\kappa\alpha\beta}V^{\tau\mu\gamma}\Bigg),\\
 \{U,V\}^{q^\prime l^\prime k^\prime}_{22}=&{\bf i}\sum_{|j|> N}\Bigg(\sum_{\substack{\tau+\kappa=q^\prime,(\alpha-e_j)+\mu=l^\prime,\beta+(\gamma-e_j)=k^\prime\\2|\kappa|+|\alpha+\beta|=h,2|\tau|+|\mu+\gamma|=g}} j\alpha_j\gamma_jU^{\kappa\alpha\beta}V^{\tau\mu\gamma}\nonumber\\
&-
\sum_{\substack{\tau+\kappa=q^\prime,\alpha+(\mu-e_j)=l^\prime,(\beta-e_j)+\gamma=k^\prime\\2|\kappa|+|\alpha+\beta|=h,2|\tau|+|\mu+\gamma|=g}} j\beta_j\mu_jU^{\kappa\alpha\beta}V^{\tau\mu\gamma}\Bigg).
\end{align*}
Consider the nonzero components of vector $q^\prime, k^\prime$ and $l^\prime$. Taking multiplicity into account, we have that $(q^\prime, l^\prime, k^\prime)$ has $|q^\prime|+| l^\prime|+| k^\prime|$ non zero components whose values are $1$. We denote them by
\begin{align*}
q_{j_{1}},\cdots, q_{j_{|q^\prime|}}, l_{j_{|q^\prime|+1}},\cdots,l_{j_{|q^\prime|+|l^\prime|}},k_{j_{|q^\prime|+|l^\prime|+1}},\cdots, k_{j_{|q^\prime|+|l^\prime|+|k^\prime|}}.
\end{align*}

 The estimates of $\{U,V\}^{q^\prime l^\prime k^\prime}_{21}$ can be obtained as follows.

Since $V$ has at most degree $2$ in $(z_j,\bar{z}_j)_{|j|>N}$, there are at least $|\tau|+|\mu+\gamma|-3$ components of $(q^\prime, l^\prime, k^\prime)$ coming from $(\tau,\mu,\gamma)$ with indexes bounded by $N$ or coming from $1,\cdots, n$. Therefore, the choice of $|\tau|+|\mu+\gamma|-3$ components of $(q^\prime, l^\prime, k^\prime)$ is less than $(2N+n)^{|\tau|+|\mu+\gamma|-3}(\le (2N+n)^{g-3})$. As for the remaining three components of $(\tau,\mu,\gamma)$, one position is $j$ with $|j|\le N$ and the choice is $2N$. Another position among the other two can be selected from the rest of the nonzero components $(q^\prime, l^\prime, k^\prime)$, the choice is $\left(\begin{array}{l}
1\\
|\kappa|+|\alpha|+|\beta|
\end{array}
\right)=|\kappa|+|\alpha|+|\beta|(\le h)$. The last one may be determined by the momentum conservation of $V$. Thus, we conclude that
\begin{align*}
\Vert\{U,V\}^{q^\prime l^\prime k^\prime}_{21}\Vert_{D(s)\times \Pi}\le &(2K)^nh(2N+n)^{g-3}2N\\
&\sup_{\kappa,\alpha,\beta}\Vert U^{\kappa\alpha\beta}\Vert_{D(s)\times \Pi}\cdot M_{\mu\gamma}\sup_{\tau,\mu,\gamma}\Vert V^{\tau\mu\gamma}\Vert_{D(s)\times \Pi}.
\end{align*}

The estimate of $\{U,V\}^{q^\prime l^\prime k^\prime}_{22}$ can be obtained as follows.

Since $V$ has degree at most $2$ in $(z_j,\bar{z}_j)_{|j|>N}$, there are at least $|\tau|+|\mu+\gamma|-2$ components of $(q^\prime, l^\prime, k^\prime)$ coming from $(\tau,\mu,\gamma)$ with indexes bounded by $N$ or coming from $1,\cdots, n$ and the choice of $|\tau|+|\mu+\gamma|-2$ components of $(q^\prime, l^\prime, k^\prime)$ is less than $(2N+n)^{|\tau|+|\mu+\gamma|-2}$. One position of the last two component can be selected from the rest of the nonzero components $(q^\prime, l^\prime, k^\prime)$, the choice is $\left(\begin{array}{l}
1\\
|\kappa|+|\alpha|+|\beta|
\end{array}
\right)=|\kappa|+|\alpha|+|\beta|$. The last one position may be determined by the momentum conservation of $V$. Thus, we conclude that
\begin{align*}
\Vert\{U,V\}^{q^\prime l^\prime k^\prime}_{21}\Vert_{D(s)\times \Pi}\le &(2K)^nh(2N+n)^{g-2}\\
&\sup_{\kappa,\alpha,\beta}\Vert U^{\kappa\alpha\beta}\Vert_{D(s)\times \Pi}\cdot M_{\mu\gamma}\sup_{\tau,\mu,\gamma}\Vert V^{\tau\mu\gamma}\Vert_{D(s)\times \Pi}.
\end{align*}

Similarly, one can prove
\begin{align*}
&\Vert\{U,V\}^{q^\prime l^\prime k^\prime}_{1}\Vert_{D(s-\sigma)\times \Pi}\\
 \le& (2K)^nh(2N+n)^{g-2} \Big(\sup_{\kappa,\alpha,\beta}\Vert \partial_xU^{\kappa\alpha\beta}\Vert_{D(s-\sigma)\times \Pi}\cdot \sup_{\tau,\mu,\gamma} \Vert V^{\tau\mu\gamma}\Vert_{D(s-\sigma)\times \Pi}\\
&+\sup_{\kappa,\alpha,\beta}\Vert U^{\kappa\alpha\beta}\Vert_{D(s-\sigma)\times \Pi}\cdot \sup_{\tau,\mu,\gamma} \Vert \partial_xV^{\tau\mu\gamma}\sup_{\tau,\mu,\gamma}\Vert_{D(s-\sigma)\times \Pi}\Big)\\
\le& (2K)^nh(2N+n)^{g-2}\cdot \frac{1}{\sigma}\sup_{\kappa,\alpha,\beta}\Vert U^{\kappa\alpha\beta}\Vert_{D(s)\times \Pi}\cdot  \sup_{\tau,\mu,\gamma}\Vert V^{\tau\mu\gamma}\Vert_{D(s)\times \Pi}.
%\le& (2K)^nh(2N+n)^{g-2} (\frac{1}{s_0-s+\sigma}\sup_{\kappa,\alpha,\beta}\Vert \partial_xU^{\kappa\alpha\beta}\Vert_{D(s_0)\times \Pi}\cdot \Vert V^{\tau\mu\gamma}\sup_{\tau,\mu,\gamma}\Vert_{D(s)\times \Pi}\\
%&+ \sup_{\kappa,\alpha,\beta}\Vert U^{\kappa\alpha\beta}\Vert_{D(s_0)\times \Pi}\cdot  K\Vert V^{\tau\mu\gamma}\sup_{\tau,\mu,\gamma}\Vert_{D(s)\times \Pi})
\end{align*}

%In conclusion, we have
%$$\Vert\{U,V\}^{q^\prime l^\prime k^\prime}\Vert_{D(s-\sigma)\times \Pi}$$
%$$\le CK^n(2N+h)^{g-2}h\frac{1}{s_0-s+\sigma}\sup_{\kappa,\alpha,\beta}\Vert U^{\kappa\alpha\beta}\Vert_{D(s_0)\times \Pi}\cdot M_{\mu\gamma}\sup_{\tau,\mu,\gamma}\Vert V^{\tau\mu\gamma}\Vert_{D(s)\times \Pi}.$$
Thus, one can easily obtains
 $$||| \{U,V\}|||_{D(s-\sigma,r-\sigma^\prime,r-\sigma^\prime)\times \Pi}\le C\frac{K^{n+1}(2N+n)^{g-2}}{\sigma\sigma^\prime r} ||| U|||_{D(s,r,r)\times \Pi}||| V|||^*_{D(s,r,r)\times \Pi}.$$

\end{proof}

The following is the Tame property which is useful for the long time stability estimate of KAM tori.

\begin{Lem}\label{guji}
Suppose the function
$$ W(x,y,z,\bar{z};\xi)=\sum_{\substack{\alpha\in\mathbb{N}^n,\mu,\gamma\in\mathbb{N}^{\mathbb{Z}_*}\\2|\alpha|+|\mu+\gamma|=h}} W^{\alpha\mu\gamma}(x;\xi)y^\alpha z^\mu\bar{z}^\gamma: D(s,r,r)\times\Pi\rightarrow \mathbb{R},$$
belongs to $\mathcal{H}$. Then, the following estimate hold
\begin{align}
&\Vert \partial_xW\Vert_{D(s-\sigma,r,r)\times \Pi}\le \frac{c_1^{h-2}}{e\sigma} |||W|||_{D(s,r,r)\times \Pi}r^{-h}\sum_{\alpha}|y|^{|\alpha|}\Vert z\Vert_2^{h-2|\alpha|},\\
&\Vert \partial_yW\Vert_{D(s-\sigma,r,r)\times \Pi}\le c_1^{h-2}h|||W|||_{D(s,r,r)\times \Pi}r^{-h}\sum_{\alpha}|y|^{|\alpha|-1}\Vert z\Vert_2^{h-2|\alpha|},
\end{align}
\begin{equation}
  \begin{split}
    \Vert ({\bf i} j \partial_{\bar{z}_j}W)_{j\in\mathbb{Z}_*}\Vert_{p-1}, \Vert (-{\bf i} j \partial_{z_j}W)_{j\in\mathbb{Z}_*}\Vert_{p-1} &\le  c_1^{h-2} h^{p+2} \frac{C_{J_n}}{e} |||W|||_{D(s,r,r)\times \Pi}r^{-h}\nonumber, \\
    & \cdot\sum_{\alpha}(\frac{1}{\sigma^{p-1}}|y|^{|\alpha|}\Vert z\Vert_1^{h-2|\alpha|-1}+\Vert z\Vert_s \Vert z\Vert_1^{h-2|\alpha|-2}),
  \end{split}
\end{equation}
%&\Vert X_W\Vert_{p-1,D(s-\sigma,r,r)\times \Pi}\le c_1^{h-2} h^{p+2} \frac{C_J}{e} |||W|||_{D(s,r,r)\times \Pi}r^{-h}\nonumber\\
%&\cdot\sum_{\alpha}(\frac{1}{\sigma}|y|^{|\alpha|}\Vert z\Vert_2^{h-2|\alpha|}+\frac{1}{\sigma^{p-1}}|y|^{|\alpha|}\Vert z\Vert_2^{h-2|\alpha|-1}+\Vert z\Vert_s \Vert z\Vert_2^{h-2|\alpha|-2}),\label{bush}
where $c_1=\sqrt{\sum_{j\in\mathbb{Z}^*}j^{-2}}$ and $C_{J_n}$ is a constant depending on the index set ${J_n}$. In addition, the following estimate hold
$$\Vert X_W\Vert_{p-1,D(s-\sigma,r,r)\times \Pi}\le c^{h-2} \frac{h^{p+2}}{\sigma^p} |||W|||_{D(s,r,r)\times \Pi}r^{-2},$$
where $c$ is a constant .

\end{Lem}

\begin{proof}
For $\partial_xW$, we have
$$\partial_xW(x,y,z,\bar{z};\xi)=\sum_{\substack{k\in\mathbb{Z}^n,\alpha\in\mathbb{N}^n,\mu,\gamma\in\mathbb{N}^{\mathbb{Z}_*}\\2|\alpha|+|\mu+\gamma|=h}} {\bf i} kW^{\alpha\mu\gamma}(k;\xi)e^{{\bf i}\langle k,x\rangle}y^\alpha z^\mu\bar{z}^\gamma.$$
Thus, we obtain
$$\Vert \partial_xW\Vert_{D(s-\sigma,r,r)\times \Pi}\le\sum_{k,\alpha} |k|e^{|k|(s-\sigma)} \sup_{\mu,\gamma} |W^{\alpha\mu\gamma}(k;\xi)||y|^{\alpha}\sum_{\substack{|\mu+\gamma|=h-2|\alpha|\\M(k,\mu,\gamma)=0}}| z^\mu\bar{z}^\gamma|.$$
By H\"older inequality, we have
$$\sum_{\substack{|\mu+\gamma|=h-2|\alpha|\\M(k,\mu,\gamma)=0}}| z^\mu\bar{z}^\gamma|\le\Bigg\Vert\Big(\sum_{\substack{|\mu+\gamma|=h-2\alpha\\j=M(k,\mu-e_j,\gamma)}} |z^{\mu-e_j}z^\gamma|\Big)_{j\in\mathbb{Z}^*}\Bigg\Vert_{l^2}\cdot \Vert z\Vert_{l^2},$$
It follows from the Young inequality
$$\Big\Vert(\sum_{k\in\mathbb{Z}^*}a_{j\pm k}b_k)_{j\in\mathbb{Z}^*}\Big\Vert_{l^2}\le \Vert a\Vert_{l^2}\Vert b\Vert_{l^1}$$
that, for each fixed $k$, the following estimate holds
\begin{align*}
\sum_{\substack{|\mu+\gamma|=h-2|\alpha|\\M(k,\mu,\gamma)=0}}| z^\mu\bar{z}^\gamma|\le  c_1^{h-2|\alpha|-2}\sum_{2|\alpha|+|\mu+\gamma|=h}\Vert z\Vert_1^{|\mu+\gamma|}.
\end{align*}
Thus, we get
\begin{align}
\Vert \partial_xW\Vert_{D(s-\sigma,r,r)\times \Pi}\le&\frac{c_1^{h-2}}{e\sigma}\sup_{\alpha,\mu,\gamma}\Vert W^{\alpha\mu\gamma}(x;\xi)\Vert\sum_{2|\alpha|+|\mu+\gamma|=h}|y|^\alpha\Vert z\Vert_1^{|\mu+\gamma|}\nonumber\\
= &\frac{c_1^{h-2}}{e\sigma} |||W|||_{D(s,r,r)\times \Pi}r^{-h}\sum_{\alpha}|y|^{|\alpha|}\Vert z\Vert_1^{h-2|\alpha|}\label{v1},
\end{align}

Similarly, for $\partial_yW$, one has
$$\partial_yW(x,y,z,\bar{z};\xi)=\sum_{\substack{k\in\mathbb{Z}^n,\alpha\in\mathbb{N}^n,\mu,\gamma\in\mathbb{N}^{\mathbb{Z}_*}\\2|\alpha|+|\mu+\gamma|=h}} \alpha W^{\alpha\mu\gamma}(k;\xi)e^{{\bf i}\langle k,x\rangle}y^{\alpha-1} z^\mu\bar{z}^\gamma.$$
Therefore, one can prove
\begin{align}\label{ve2}
\Vert \partial_yW\Vert_{D(s-\sigma,r,r)\times \Pi}\le&c_1^{h-2}h|||W|||_{D(s,r,r)\times \Pi}r^{-h}\sum_{\alpha}|y|^{|\alpha|-1}\Vert z\Vert_1^{h-2|\alpha|}.
\end{align}

For $({\bf i} j \partial_{\bar{z}_j}W,-{\bf i} j \partial_{z_j}W)_{j\in\mathbb{Z}_*}$, we have
\begin{align*}
&({\bf i} j \partial_{\bar{z}_j}W,-{\bf i} j \partial_{z_j}W)_{j\in\mathbb{Z}_*}\\
=&\Big({\bf i} j \gamma_j\sum_{\substack{k\in\mathbb{Z}^n,\alpha\in\mathbb{N}^n,\mu,\gamma\in\mathbb{N}^{\mathbb{Z}_*}\\2|\alpha|+|\mu+\gamma|=h}} W^{\alpha\mu\gamma}(k;\xi)e^{{\bf i}\langle k,x\rangle}y^{\alpha} z^\mu\bar{z}^{\gamma-e_j},\\
& -{\bf i} j \mu_j\sum_{\substack{k\in\mathbb{Z}^n,\alpha\in\mathbb{N}^n,\mu,\gamma\in\mathbb{N}^{\mathbb{Z}_*}\\2|\alpha|+|\mu+\gamma|=h}} W^{\alpha\mu\gamma}(k;\xi)e^{{\bf i}\langle k,x\rangle}y^{\alpha} z^{\mu-e_j}\bar{z}^{\gamma}\Big)_{j\in\mathbb{Z}_*}.
\end{align*}
Then, we have
\begin{align}
&\Vert ({\bf i} j \partial_{\bar{z}_j}W)_{j\in\mathbb{Z}_*}\Vert_{p-1}=\Vert ( \partial_{\bar{z}_j}W)_{j\in\mathbb{Z}_*}\Vert_{p}\nonumber\\
=&\left(\sum_{j\in\mathbb{Z}_*}\left(\sum_{k,2|\alpha|+|\mu+\gamma|=h}|\gamma_j||W^{\alpha\mu\gamma}(k;\xi)e^{{\bf i}\langle k,x\rangle}y^{\alpha} z^\mu\bar{z}^{\gamma-e_j}|\cdot |j|^p\right)^2\right)^\frac{1}{2}\nonumber\\
\le &\sum_{k,\alpha}\sup_{\mu,\gamma}|W^{\alpha\mu\gamma}(k;\xi)e^{|k|(s-\sigma)}|y|^\alpha \left(\sum_{j\in\mathbb{Z}_*}\left(\sum_{|\mu+\gamma|=h-2|\alpha|}|\gamma_j|| z^\mu\bar{z}^{\gamma-e_j}|\cdot |j|^p\right)^2\right)^\frac{1}{2}.
%= & \sum_{k,\alpha}\sup_{\mu,\gamma}|W^{\alpha\mu\gamma}(k;\xi)e^{|k|(s-\sigma)}|y|^\alpha\nonumber \\
%&\cdot\left(\sum_{j^\prime-\sum_{b}k_bj_b\in\mathbb{Z}_*}\left(\sum_{|\mu+\gamma|=h-2|\alpha|}|\gamma_{j^\prime-\sum_{b}k_bj_b}|| z^\mu\bar{z}^{\gamma-e_{j^\prime-\sum_{b}k_bj_b}}|\cdot \left|j^\prime-\sum_{b}k_bj_b\right|^p\right)^2\right)^\frac{1}{2}
\label{w11}
\end{align}
Since $j=M(k,\mu-e_j,\gamma)=\sum_{b}k_bj_b-\sum_{\beta\in \mathbb{Z}^*, \beta \neq j} \beta\mu_\beta +\sum_{\beta \in \mathbb{Z}^*}\beta \gamma_\beta$, we have
$$|j|^p\le h^p\big(C_{J_n}|k|^p+\sum_{\beta \in \mathbb{Z}^*}|\beta|^p\mu_\beta+\sum_{\beta\in \mathbb{Z}^*,\beta \neq j}|\beta|^p\gamma_\beta+(\gamma_j-1)|j|^p\big),$$
where $C_{J_n}$ is a constant depending on the index set $J_n$.
Therefore, one gets
\begin{align*}
&\left(\sum_{j\in\mathbb{Z}_*}\left(\sum_{|\mu+\gamma|=h-2|\alpha|}|\gamma_j|| z^\mu\bar{z}^{\gamma-e_j}|\cdot |j|^p\right)^2\right)^\frac{1}{2}\\
\le & h^p\left(\sum_{j\in\mathbb{Z}_*}\left(\sum_{|\mu+\gamma|=h-2|\alpha|}|\gamma_j|| z^\mu\bar{z}^{\gamma-e_j}|\right.\right.\\
&\left.\left.\cdot \Big(C_{J_n}|k|^p+\sum_{\beta \in \mathbb{Z}^*}|\beta|^p\mu_\beta+\sum_{\beta \in \mathbb{Z}^*, \beta \neq j}|\beta|^p\gamma_\beta +(\gamma_j-1)|j|^p\Big)\right)^2\right)^\frac{1}{2}\\
\le &h^p C_{J_n}|k|^p \left(\sum_{j\in\mathbb{Z}_*}\left(\sum_{|\mu+\gamma|=h-2|\alpha|}|\gamma_j| | z^\mu\bar{z}^{\gamma-e_j}|\ \right)^2\right)^\frac{1}{2}\\
&+ h^p\sum_{t\in \mathbb{Z}^*}\left(\sum_{j\in\mathbb{Z}_*}\left(\sum_{|\mu+\gamma|=h-2|\alpha|}|\gamma_j|| z^\mu\bar{z}^{\gamma-e_j}|\cdot |t|^p\mu_t\right)^2\right)^\frac{1}{2}\\
&+ h^p\sum_{t\in \mathbb{Z}^*,t\neq j}\left(\sum_{j\in\mathbb{Z}_*}\left(\sum_{|\mu+\gamma|=h-2|\alpha|}|\gamma_j|| z^\mu\bar{z}^{\gamma-e_j}|\cdot |t|^p\gamma_t\right)^2\right)^\frac{1}{2}\\
&+ h^p\left(\sum_{j\in\mathbb{Z}_*}\left(\sum_{|\mu+\gamma|=h-2|\alpha|}|\gamma_j|| z^\mu\bar{z}^{\gamma-e_j}|\cdot (\gamma_j-1)|j|^p\right)^2\right)^\frac{1}{2}.
\end{align*}
By Young's inequality, we have
\begin{align*}
&\left(\sum_{j\in\mathbb{Z}_*}\left(\sum_{|\mu+\gamma|=h-2|\alpha|}|\gamma_j|| z^\mu\bar{z}^{\gamma-e_j}|\cdot |j|^p\right)^2\right)^\frac{1}{2}\\
\le &h^{p+2} C_{J_n}|k|^p c_1^{h-2|\alpha|-1}\Vert z\Vert_1^{h-2|\alpha|-1}+h^{p+2}c_1^{h-2|\alpha|-2} \Vert z\Vert_s\Vert z\Vert_1^{h-2|\alpha|-2}.
\end{align*}
Substituting it into \eqref{w11}, one gets
\begin{align}\label{v3}
\Vert ({\bf i} j  \partial_{\bar{z}_j}W)_{j\in\mathbb{Z}_*}\Vert_{p-1}\le&  c_1^{h-2} h^{p+2} \frac{C_{J_n}}{e} |||W|||_{D(s,r,r)\times \Pi}r^{-h}\nonumber\\
&\cdot\sum_{\alpha}(\frac{1}{\sigma^{p-1}}|y|^{|\alpha|}\Vert z\Vert_1^{h-2|\alpha|-1}+\Vert z\Vert_s \Vert z\Vert_1^{h-2|\alpha|-2}).
\end{align}
Similar estimate holds for $\Vert (-{\bf i} j  \partial_{\bar{z}_j}W)_{j\in\mathbb{Z}_*}\Vert_{p-1}$.
%Combining \eqref{v1}, \eqref{ve2} and \eqref{v3} together, it's easy to obtain the final estimate \eqref{bush}

\end{proof}

%\begin{Lem}\label{huahua}
%If $p>\frac{1}{2}$, then $\Vert a*b\Vert_{p}\le c \Vert a\Vert_p\Vert b\Vert_p$ for $a,b\in l^p$ with a finite constant $c$.
%\end{Lem}
%The proof of this Lemma refer to \cite{}.

\begin{Lem}\label{youy}
Suppose the function
$$ W(x,y,z,\bar{z};\xi)=\sum_{\substack{\alpha\in\mathbb{N}^n,\mu,\gamma\in\mathbb{N}^{\mathbb{Z}_*}\\2|\alpha|+|\mu+\gamma|=h}} W^{\alpha\mu\gamma}(x;\xi)y^\alpha z^\mu\bar{z}^\gamma: D(s,r,r)\times\Pi\rightarrow \mathbb{R},$$
belong to $\mathcal{H}^*$. Then, the following estimate hold
$$\Vert X_W\Vert_{p,D(s-\sigma,r,r)\times \Pi}\le c^{h-2} \frac{h^{p+2}}{e\sigma^p} |||W|||^*_{D(s,r,r)\times \Pi}r^{-2},$$
where $c$ is a constant .
\end{Lem}

%\begin{proof}
%We only consider the estimate for $({\bf i} j W_{\bar{z}_j})_{j\in\mathbb{Z}_*}$. The other estimate can obtain easier than it. Actually, from the proof the the previous lemma, we have
%\begin{align}
%&\Vert ({\bf i} j W_{\bar{z}_j})_{j\in\mathbb{Z}_*}\Vert_{p-1}=\Vert ( W_{\bar{z}_j})_{j\in\mathbb{Z}_*}\Vert_{p}\nonumber\\
%=&\left(\sum_{j\in\mathbb{Z}_*}\left(\sum_{k,2|\alpha|+|\mu+\gamma|=h}|\gamma_j||W^{\alpha\mu\gamma}(k;\xi)e^{{\bf i}\langle k,x\rangle}y^{\alpha} z^\mu\bar{z}^{\gamma-e_j}|\cdot |j|^p\right)^2\right)^\frac{1}{2}\nonumber\\
%\le &\sum_{k,\alpha}\sup_{\mu,\gamma}|W^{\alpha\mu\gamma}(k;\xi)e^{|k|s}|y|^\alpha \left(\sum_{j\in\mathbb{Z}_*}\left(\sum_{|\mu+\gamma|=h-2|\alpha|}|\gamma_j|| z^\mu\bar{z}^{\gamma-e_j}|\cdot |j|^p\right)^2\right)^\frac{1}{2}
%\end{align}
%From Lemma \ref{huahua}, it is easy to obtain
%\begin{align*}
%\Vert ({\bf i} j W_{\bar{z}_j})_{j\in\mathbb{Z}_*}\Vert_{p-1}\le& \sum_{k,\alpha}\sup_{\mu,\gamma}|W^{\alpha\mu\gamma}(k;\xi)e^{|k|s}|y|^\alpha (c\Vert z\Vert_p)^{|\mu+\gamma|-1}\\
%\le& c^{h-2} |||W|||_{D(s,r,r)\times \Pi}r^{-2}.
%\end{align*}
%\end{proof}
The proof of this Lemma is similar to that of Lemma \ref{guji}. Here we omit it.

\section{Small denominator equation with large variable coefficient}\label{slarge}

Let $s_{0} > 0$ be constant and $\epsilon_{0} > 0$ de small constant. Define sequences $\left\{ {{s_m}} \right\}_{m = 0}^\infty$ and $\left\{ {{\varepsilon _m}} \right\}_{m = 0}^\infty$ as follows,

\begin{itemize}

  \item ${s_0} = {s_0}$,\ \ ${s_m} = {s_0}\left( {1 - \frac{{\sum\limits_{j = 1}^m {\frac{1}{{{j^2}}}} }}{{100\sum\limits_{j = 1}^\infty  {\frac{1}{{{j^2}}}} }}} \right),m = 1,2, \cdots $;
   \item $s_m^{(i)}=s_{m+1}+\left(1-\frac{i}{10}\right)(s_m-s_{m+1})$,\; $i=0,1,...,10$;
  \item  ${\varepsilon _0} = {\varepsilon _0}$, \ \  ${\varepsilon _m} = \varepsilon _0^{{{\left( {{4 \mathord{\left/
 {\vphantom {4 3}} \right.
 \kern-\nulldelimiterspace} 3}} \right)}^m}},m = 1,2, \cdots$ ;
\end{itemize}

Clearly, ${s_0} > {s_1} >  \cdots {s_m} > {s_{m + 1}} >  \cdots  > {{{s_0}} \mathord{\left/
 {\vphantom {{{s_0}} 2}} \right.
 \kern-\nulldelimiterspace} 2}$ and ${s_m} > s_m^{\left( 1 \right)} >  \cdots s_m^{\left( l \right)} >  \cdots  > {s_{m + 1}}$.

Define $D(s) = \left\{ {\theta  \in {{{{{C^d}} \mathord{\left/
 {\vphantom {{{C^d}} {\left( {2\pi Z} \right)}}} \right.
 \kern-\nulldelimiterspace} {\left( {2\pi \mathbb Z} \right)}}}^d}:\left| {\mathrm{Im}\theta } \right| \le s} \right\}$, thus $$ D({s_0})\supset D(s_1) \supset D(s_2) \supset  \cdots  \supset D(s_m) \supset  \cdots  \supset D(s_0/2),$$
 and
 $$D(s_m)\supset D( s_m^{\left( 1 \right)}) \supset D( s_m^{\left( 2 \right)})  \supset  \cdots  \supset D( s_m^{\left( l \right)})  \supset  \cdots  \supset D(s_{m + 1}).$$

 For any function $ f: D(s) \to \mathbb C$, we define ${\left\| f \right\|_s} = \mathop {\sup }\limits_{\theta  \in D(s)} \left| {f\left( \theta  \right)} \right|,\forall s \ge 0$. We call that a function $ f:D(s) \to \mathbb C$ is real analytic if it is analytic in $D(s)$ and it is real for real variable $\theta  \in D(s)$.

 Assume that we have a function $a\left( \theta  \right):D({s_m}) \to \mathbb C$ which obeys

\begin{itemize}
  \item

  \begin{equation}
  a\left( \theta  \right) = \sum\limits_{j = 0}^m {{a_j}\left( \theta  \right)},
  \label{1.1}
  \end{equation}

  \item ${a_j}\left( \theta  \right): D({s_j}) \to \mathbb C$ is real analytic, and
 \begin{equation}
\begin{array}{l}
{\left\| {{a_j}} \right\|_{{s_j}}} \le \varepsilon _j^{1 - }.
\label{1.2}
\end{array}
\end{equation}
\end{itemize}

Let $\omega  \in {\left[ {0,1} \right]^d}$ be a parameter vector and $\lambda  \ne 0$ be a constant which may be large. Assume also that $\left( {\omega ,\lambda } \right)$ is of Diophantine, i.e., there are constants $0 < \gamma  \ll 1$ and $\tau  > d$ such that
 \begin{equation}
\begin{array}{l}
\left| {\left\langle {k,\omega } \right\rangle } \right| \ge {{{\gamma _0}} \mathord{\left/
 {\vphantom {{{\gamma _0}} {{{\left| k \right|}^\tau },\forall k \in {Z^d}\backslash \left\{ 0 \right\}}}} \right.
 \kern-\nulldelimiterspace} {{{\left| k \right|}^\tau },\forall k \in {\mathbb Z^d}\backslash \left\{ 0 \right\}}},
 \label{1.3}
\end{array}
\end{equation}
and
 \begin{equation}
 \left| { \pm \left\langle {k,\omega } \right\rangle  + \lambda } \right| \ge {\gamma_0 \mathord{\left/
 {\vphantom {\gamma  {{{\left| k \right|}^\tau }}}} \right.
 \kern-\nulldelimiterspace} {{{\left| k \right|}^\tau }}},\forall k \in {\mathbb Z^d}.
 \label{1.4}
 \end{equation}
 where we use notation $|k|=1$, when $k = 0 \in {\mathbb Z^d}$ for convenience.

 \begin{Lem}\label{1 lemma}
 Let $A\left( \theta  \right):D(s) \to \mathbb C$ be real analytic for $s > 0$. Write ${D_\omega }: = \mathbf i \,\omega  \cdot {\partial _\theta }$, where ${\bf i}^{2}=-1$. Consider the differential equation
 \begin{equation}
 {D_\omega }{\rm X}\left( \theta  \right) = A\left( \theta  \right).
 \label{1.5}
 \end{equation}
Then \eqref{1.5} has unique solution ${\rm X} = {\rm X}\left( \theta  \right) = D_\omega ^{ - 1}A\left( \theta  \right): {D(s') \to \mathbb C}$ which is real analytic and
\begin{equation}
{\left\| {{\rm X}\left( \theta  \right)} \right\|_{{s^\prime}}} \le C\left( {d,\tau } \right)\frac{1}{{{\gamma _0}}}\frac{1}{{{{\left( {s - {s^\prime}} \right)}^{10\left( {d + \tau } \right)}}}}{\left\| A \right\|_s},\ \ \forall\ 0  < s'< s,
\label{1.6}
\end{equation}
where $C\left( {d,\tau } \right)$ is a constant depending on $d$ and $\tau$.
\end{Lem}
\begin{proof}
The proof is trivial in $KAM$ theory. We omit it.
\end{proof}
\bigskip

Let
\begin{equation}
\alpha \left( \theta  \right) = \sum\limits_{j = 1}^m {{\alpha _j}\left( \theta  \right)} ,
{\alpha _j}\left( \theta  \right) = D_\omega ^{ - 1}{a_j}\left( \theta  \right).
\label{1.7}
\end{equation}
By Lemma \ref{1.1}, ${\alpha _j}\left( \theta  \right):D(s_j^{\left( 1 \right)}) \to \mathbb C$ is real analytic with
\begin{equation}
{\left\| {{\alpha _j}\left( \theta  \right)} \right\|_{s_j^{\left( 1 \right)}}} \le C\left( {d,\tau } \right)\frac{1}{{{\gamma _0}}}{\left( {\frac{1}{{{s_j} - s_j^{\left( 1 \right)}}}} \right)^{10\left( {d + \tau } \right)}}{\left\| {{a_j}} \right\|_{{s_j}}} < \varepsilon _j^{1 - }.
\label{1.8}
\end{equation}
Thus
\begin{equation}
{\left\| \alpha  \right\|_{s_m^{\left( 1 \right)}}} \le \varepsilon _0^{1 - }.
\label{1.9}
\end{equation}

\begin{Lem}\label{2 lemma}
Taking any $\theta  = {\mathop{\rm Re}\nolimits}\ \theta  + \mathbf i \, {\mathop{\rm Im}\nolimits}\ \theta  \in D(s_j^{\left( 2 \right)})$ ,we have $\left| {{\mathop{\rm Im}\nolimits}\ {\alpha _j}\left( \theta  \right)} \right| < {\varepsilon_j ^{\frac{1}{{10}}}}\left| {{\mathop{\rm Im}\nolimits}\ \theta } \right|$.
\end{Lem}
\begin{proof}
Let $x = {\rm Re}\ \theta$, $y ={\rm Im}\ \theta$. Then
\begin{equation*}
\begin{aligned}
{\alpha _j}\left( \theta  \right) & = {\alpha _j}\left( {x + {\bf i}y} \right)\\
& = \sum_{p = 0}^{\infty } {\sum_{k \in \mathbb N^d, {|k|= p}}{\frac{{\partial _x^k{\alpha _j}\left( x \right)}}{{k!}} {\bf i}^{p}{y^k}}}\\
& = \sum_{p \in 2\mathbb N} {\sum_{k \in \mathbb N^d,|k|=p}{\frac{{\partial _x^k{\alpha _j}\left( x \right)}}{{k!}} {\bf i}^{p}{y^k}}}+\sum_{p \in { 2 \mathbb N+1  }} {\sum_{k \in \mathbb N^d,|k|=p}{\frac{{\partial _x^k{\alpha _j}\left( x \right)}}{{k!}} {\bf i}^{p}{y^k}}}
\end{aligned}
\end{equation*}
Then,
\begin{equation*}
\begin{aligned}
\left| {{\mathop{\rm Im}\nolimits}\ {\alpha _j}\left( \theta  \right)} \right| & < \left| {\sum\limits_{p \in 2\mathbb N + 1} {\left( { \pm 1} \right)\sum\limits_{\left| k \right| = p} {\frac{{\partial _x^k{\alpha _j}\left( x \right)}}{{k!}}{y^k}} } } \right|\\
 & \le \left( {\sum\limits_{p \in 2\mathbb N + 1} {\left( {\sum\limits_{\left| k \right| = p} 1 } \right){{\left( {\frac{1}{{s_j^{\left( 1 \right)}}}} \right)}^p}{{\left| y \right|}^p}} } \right){\left\| {{\alpha _j}} \right\|_{s_j^{\left( 1 \right)}}}\\
 & < C j^{10d^{2}}\frac{1}{{s_j^{(1)}}}\varepsilon _j^{1 - }\left| y \right|\\
 & \le \varepsilon _j^{\frac{1}{{10}}}\left| {{\mathop{\rm Im}\nolimits}\ \theta } \right|.
\end{aligned}
\end{equation*}
\end{proof}
This completes the proof of Lemma \ref{2 lemma}.

By Lemma \ref{2 lemma} and noting \eqref{1.7}, we have
\begin{equation}
\left| {{\mathop{\rm Im}\nolimits}\ \alpha \left( \theta  \right)} \right| \le \varepsilon _0^{\frac{1}{{20}}}\left| {{\mathop{\rm Im}\nolimits}\ \theta } \right|,\quad \forall \ \theta  \in D(s_m^{(2)}).
\label{1.10}
\end{equation}
Let $s=s_m,\, s'=s_{m}^{(2)}$. Moreover, by \eqref{1.10}, we can define the map $\tilde T : D(s') \to D(s)$ by
\begin{equation}
\theta  = \tilde T\varphi  = \varphi  + \alpha \left( \varphi  \right)\omega ,\forall \varphi  \in \mathbb T_{_{s'}}^d.
\label{1.11}
\end{equation}
By \eqref{1.9}, we can prove easily that there is a vector function $\hat \alpha :D(s'') \to D(s')$ with $s''=s_m^{(3)}$, which is real analytic such that the map $\tilde T$ is invertible, ${{\tilde T}^{ - 1}}:D(s'') \to D(s')$ is of the form
\begin{equation}
\varphi  = {{\tilde T}^{ - 1}}\left( \theta  \right) = \theta  + \hat \alpha \left( \theta  \right)\omega , \quad \forall \theta  \in D(s''),
\label{1.12}
\end{equation}
where $\hat \alpha \left( \theta  \right) = (\widehat{\alpha}(\theta):j=1,\cdots,d) \in {\mathbb C^d}$ , $\hat \alpha \left( \theta  \right)\omega  = (\alpha_{j}(\theta)\omega_{j}: j=1,\cdots,d)$, and
\begin{equation}
{\left|\left| {\hat \alpha }\right| \right|_{s''}} \le \varepsilon _0^{1 - },\quad s'' = s_m^{\left( 3 \right)}.
\label{1.13}
\end{equation}
Actually, the proof is as follows. Let $\varphi  = {\tilde T^{ - 1}}\theta  = \theta  + h\left( \theta  \right)$. Since $\theta  = \varphi  + \alpha \left( \varphi  \right)\omega $, we have $\theta  = \theta  + h\left( \theta  \right) + \alpha \left( {\theta  + h\left( \theta  \right)} \right)\omega $, by which we get a Picard sequence
\begin{equation}
\left\{ \begin{array}{l}
{h_0}\left( \theta  \right) =  - \alpha \left( \theta  \right)\omega, \\
{h_\nu }\left( \theta  \right) =  - \alpha \left( {\theta  + {h_{\nu  - 1}}\left( \theta  \right)} \right)\omega ,\quad \nu  = 1,2, \cdots.
\end{array} \right.
\label{1.14}
\end{equation}

By contraction mapping principle, as well as $\eqref{1.9}$ and $\eqref{1.10}$, there exists $$h\left( \theta  \right) = \mathop {\lim }\limits_{\nu  \to \infty } {h_\nu }\left( \theta  \right): D(s'') \to D(s')$$ which is real analytic. Let $\hat \alpha \left( \theta  \right) =  - \alpha \left( {\theta  + h\left( \theta  \right)} \right)$, we have
\begin{equation}
{\left\| {\hat \alpha } \right\|_{s''}} \le \varepsilon _0^ {1-} ,\quad \left| {{\mathop{\rm Im}\nolimits}\ \hat \alpha \left( \theta  \right)} \right| \le \varepsilon _0^ {1-} \left| {{\mathop{\rm Im}\nolimits}\ \theta } \right|,\quad \forall \theta  \in D(s'').
\label{1.15}
\end{equation}
By ${{\tilde T}^{ - 1}}$, we also induce a linear operator ${{\tilde T}^{ - 1}}:\mc L\left( D(s'')\right) \to \mc L\left( D(s') \right)$ by
\begin{equation}
\left( {{{\tilde T}^{ - 1}}f} \right)\left( \theta  \right) = f\left( {{{\tilde T}^{ - 1}}\theta } \right) = f\left( {\theta  + \hat \alpha \left( \theta  \right)\omega } \right) = f\left( \varphi  \right),
\label{1.16}
\end{equation}
where $ \mc L\left( D(s) \right)$ denotes the set of all real analytic functions defined in $D(s)$.

%%%%%%%%%%--------------------------------
\begin{Lem}\label{Improved Liu-Yuan lemma}
The homological equation
\begin{equation}
\left( {i\omega  \cdot {\partial _\varphi } + \lambda \left( {1 + a\left( \varphi  \right)} \right)} \right)x\left( \varphi  \right) = R\left( \varphi  \right),\varphi  \in D(s_m)
\label{1.27}
\end{equation}
has a unique solution $x=x(\varphi):D(s_{m+1}) \to \mathbb C$ which is real analytic and
\begin{equation}
\|x\|_{s_{m+1}} \leq C(d,\tau)\frac{1}{\gamma}(\frac{1}{s_{m}-s_{m+1}})^{20(d+\tau)} \|R\|_{s_{m}},
\label{1.28}
\end{equation}
where $(\omega,\lambda)$ satisfies the Diophantine \eqref{1.3}, \eqref{1.4} and $a=a(\theta)$ obeys \eqref{1.1},\eqref{1.2}, and $R=R{\varphi}: D(s_{m})\rightarrow \mathbb{C}$ is real analytic.
\end{Lem}
%%%%%%%%%---------------------------

\begin{proof}
Let $x=x(\varphi)$ be the solution of the homological equation \eqref{1.27}.
%\begin{equation*}\label{*}
%\left[ {i\omega  \cdot {\partial _\varphi } + \lambda \left( {1 + a\left( \varphi  \right)} \right)} \right]x\left( \varphi  %\right) = R\left( \varphi  \right)
%\eqno(*)
%\end{equation*}
%%%%%%%%%%-
And let $\tilde T^{-1} x=y$. Then $x(\varphi)=(\tilde T y)(\varphi)$ where $y=y(\theta): D(s'') \to D(s')$. Then
\begin{equation}
x\left( \varphi  \right) = y\left( {\tilde T\varphi } \right) = y\left( \theta  \right) = y\left( {\varphi  + \alpha \left( \varphi  \right)\omega } \right).
\label{1.17}
\end{equation}

By computation, we have
\begin{equation}
\begin{aligned}
\omega  \cdot {\partial _\varphi }x\left( \varphi  \right) & = \omega  \cdot {\partial _\varphi }y\left( {\varphi  + \alpha \left( \varphi  \right)\omega } \right)\\
 &= \sum\limits_{j = 1}^d {{\omega _j}{\partial _{{\varphi _1}}}y\left( {{\varphi _j} + \alpha \left( \varphi  \right){\omega _1}, \cdots ,{\varphi _d} + \alpha \left( \varphi  \right){\omega _d}} \right)} \\
& = \sum\limits_{j = 1}^d {{\omega _j}{\partial _{{\theta _j}}}y\left( \theta  \right)}  + \sum\limits_{i = 1}^d {\sum\limits_{j = 1}^d {{\omega _i}{\omega _j}{\partial _{{\theta _i}}}y\left( \theta  \right)} } {\partial _{{\varphi _j}}}\alpha \left( \varphi  \right)\\
& =( \omega  \cdot {\partial _\theta }y\left( \theta  \right)) \left( {1 + \omega  \cdot {\partial _\varphi }\alpha \left( \varphi  \right)} \right),
\end{aligned}
\label{1.18}
\end{equation}
where $\theta  = \tilde T\varphi  = \varphi  + \alpha \left( \varphi  \right)\omega $. And
\begin{equation}
\lambda \left( {1 + a\left( \varphi  \right)} \right)x\left( \varphi  \right) = \lambda \left( {1 + a\left( \varphi  \right)} \right)y\left( \theta  \right) = \lambda \left( {1 + a\left( {{{\tilde T}^{ - 1}}\theta } \right)} \right)y\left( \theta  \right).
\label{1.19}
\end{equation}

By \eqref{1.18} and \eqref{1.19}, the homological equation \eqref{1.27} can be rewritten as
\begin{equation}
\left( {\bf i}\omega  \cdot {\partial _\theta }y\left( \theta  \right)\right)\left( {1 + \omega  \cdot {\partial _\varphi }\alpha \left( \varphi  \right)} \right) + \lambda \left( {1 + a\left( \varphi  \right)} \right)y\left( \theta  \right) = R\left( \varphi  \right).
\label{1.20}
\end{equation}
By \eqref{1.7},$D_\omega ^{ - 1}a\left( \theta  \right) = \alpha \left( \theta  \right),\forall \theta  \in D(s'')$, that implies $\omega  \cdot {\partial _\varphi }\alpha \left( \varphi  \right) = a\left( \varphi  \right)$.
Let $\tilde R\left( \varphi  \right) = {\left( {1 + a\left( \varphi  \right)} \right)^{ - 1}}R\left( \varphi  \right)$. Then \eqref{1.20} reads
\begin{equation}
{\bf i}\,\omega  \cdot {\partial _\theta }y\left( \theta  \right) + \lambda y\left( \theta  \right) = \tilde R\left( \varphi  \right) = \tilde R\left( {{{\tilde T}^{ - 1}}\theta } \right): = {R^*}\left( \theta  \right),
\label{1.21}
\end{equation}
where we take $\theta  \in D(s'')$. By \eqref{1.2},
\begin{equation}
{\left\| {\tilde R} \right\|_{s''}} \le {\left( {1 + c{\varepsilon _0}} \right)^{ - 1}}{\left\| R \right\|_{{s_m}}} \le 2{\left\| R \right\|_{{s_m}}}.
\label{1.22}
\end{equation}
Moreover, by \eqref{1.15},
\begin{equation}
{\left\| {\tilde R\left( {{{\tilde T}^{ - 1}}\left( \theta  \right)} \right)} \right\|_{s''}} \le 2{\left\| R \right\|_{{s_m}}}.
\label{1.23}
\end{equation}
By Fourier expansion, the solution to \eqref{1.21} is
\begin{equation}
y\left( \theta  \right) = \sum\limits_{k \in {Z^d}} {\frac{1}{{\lambda  - \left\langle {k,\omega } \right\rangle }}} {{\hat R}^*}\left( k \right){e^{ik\theta }},\theta  \in D(s''),
\label{1.24}
\end{equation}
where $\widehat{R}^{*}(k)$ is k-Forier coefficient of $R^{*}$.
By the Diophantine Conditions \eqref{1.3},
\begin{equation}
\begin{aligned}
{\left\| y \right\|_{s_m^{(3)}}} & \leq C\left( {d,\tau } \right)\frac{1}{\gamma }\frac{1}{{{{\left( {s'' - s_{m}^{(3)}} \right)}^{10\left( {d + \tau } \right)}}}}{\left\| {{R^*}} \right\|_{s''}}\\
  & \leq C\left( {d,\tau } \right)\frac{1}{\gamma }\frac{1}{{{{\left( {s'' - s_{m}^{(3)}} \right)}^{10\left( {d + \tau } \right)}}}}{\left\| {{R^*}} \right\|_{{s_m}}},
\end{aligned}
\label{1.25}
\end{equation}

Returning to \eqref{1.17}, we have
\begin{equation}
{\left\| {x\left( \varphi  \right)} \right\|_{s_{m+1}}} \le C\left( {d,\tau } \right)\frac{1}{\gamma }\frac{1}{{{{\left( {s_{m} - s_{m+1}} \right)}^{10\left( {d + \tau } \right)}}}}{\left\| {{R^*}} \right\|_{{s_m}}}.
\label{1.26}
\end{equation}
This completes the proof.
\end{proof}

 \section{Acknowledgments}

 The authors would like to thank the anonymous referees for their helpful comments and suggestions.

\bibliographystyle{plain}
%\bibliography{Myreference}

\begin{thebibliography}{10}

\bibitem{arnold1963proof}
V.~I. Arnold.
\newblock Proof of a theorem of a. n. kolmogorov on the invariance of
  quasi-periodic motions under small perturbations of the
  hamiltonian(kolmogoroff theorem on invariance of quasi- periodic motions
  under small perturbations of hamiltonian).
\newblock {\em Russian Mathematical Surveys}, 18:9--36, 1963.

\bibitem{Arn64}
V.~I. Arnold.
\newblock Instability of dynamical systems with many degrees of freedom.
\newblock {\em Dokl. Akad. Nauk SSSR}, 156:9--12, 1964.

\bibitem{baldi2018time}
P.~Baldi, M.~Berti, E.~Haus, and R.~Montalto.
\newblock Time quasi-periodic gravity water waves in finite depth.
\newblock {\em Inventiones mathematicae}, 214(2):739--911, 2018.

\bibitem{baldi2016kam}
P.~Baldi, M.~Berti, and R.~Montalto.
\newblock Kam for autonomous quasi-linear perturbations of kdv.
\newblock In {\em Annales de l'Institut Henri Poincar{\'e} C, Analyse non
  lin{\'e}aire}, volume~33, pages 1589--1638. Elsevier, 2016.

\bibitem{baldi2019existence}
P.~Baldi and E.~Haus.
\newblock On the existence time for the kirchhoff equation with periodic
  boundary conditions.
\newblock {\em Nonlinearity}, 33(1):196, 2019.

\bibitem{Bam03}
D.~Bambusi.
\newblock Birkhoff normal form for some nonlinear {PDE}s.
\newblock {\em Comm. Math. Phys.}, 234(2):253--285, 2003.

\bibitem{BG}
D.~Bambusi and B.~Gr\'{e}bert.
\newblock Birkhoff normal form for partial differential equations with tame
  modulus.
\newblock {\em Duke Math. J.}, 135(3):507--567, 2006.

\bibitem{BFG88}
G.~Benettin, J.~Fr\"{o}hlich, and A.~Giorgilli.
\newblock A {N}ekhoroshev-type theorem for {H}amiltonian systems with
  infinitely many degrees of freedom.
\newblock {\em Comm. Math. Phys.}, 119(1):95--108, 1988.

\bibitem{bernier2020long}
J.~Bernier, E.~Faou, and B.~Gr{\'e}bert.
\newblock Long time behavior of the solutions of nlw on the-dimensional torus.
\newblock In {\em Forum of Mathematics, Sigma}, volume~8, page e12. Cambridge
  University Press, 2020.

\bibitem{BBP13}
M.~Berti, L.~Biasco, and M.~Procesi.
\newblock K{AM} theory for the {H}amiltonian derivative wave equation.
\newblock {\em Ann. Sci. \'{E}c. Norm. Sup\'{e}r. (4)}, 46(2):301--373 (2013),
  2013.

\bibitem{berti2022birkhoff}
M.~Berti, R.~Feola, and F.~Pusateri.
\newblock Birkhoff normal form and long time existence for periodic gravity
  water waves.
\newblock {\em Communications on Pure and Applied Mathematics}, 2022.

\bibitem{biasco2020abstract}
L.~Biasco, J.~E. Massetti, and M.~Procesi.
\newblock An abstract birkhoff normal form theorem and exponential type
  stability of the 1d nls.
\newblock {\em Communications in Mathematical Physics}, 375(3):2089--2153,
  2020.

\bibitem{Bou96-GAFA}
J.~Bourgain.
\newblock Construction of approximative and almost periodic solutions of
  perturbed linear {S}chr\"{o}dinger and wave equations.
\newblock {\em Geom. Funct. Anal.}, 6(2):201--230, 1996.

\bibitem{CK}
J.~Colliander, M.~Keel, G.~Staffilani, H.~Takaoka, and T.~Tao.
\newblock Transfer of energy to high frequencies in the cubic defocusing
  nonlinear {S}chr\"{o}dinger equation.
\newblock {\em Invent. Math.}, 181(1):39--113, 2010.

\bibitem{cong2015long}
H.~Cong, M.~Gao, and J.~Liu.
\newblock Long time stability of kam tori for nonlinear wave equation.
\newblock {\em Journal of Differential Equations}, 258(8):2823--2846, 2015.

\bibitem{CLY}
H.~Cong, J.~Liu, and X.~Yuan.
\newblock Stability of {KAM} tori for nonlinear {S}chr\"{o}dinger equation.
\newblock {\em Mem. Amer. Math. Soc.}, 239(1134):vii+85, 2016.

\bibitem{guardia2015growth}
M.~Guardia and V.~Kaloshin.
\newblock Growth of sobolev norms in the cubic defocusing nonlinear
  schr{\"o}dinger equation.
\newblock {\em Journal of the European Mathematical Society}, 17(1):71--149,
  2015.

\bibitem{he2021long}
X.~He, J.~Shi, and X.~Yuan.
\newblock Long time stability of kam tori for the nonlinear schr$\backslash$"
  odinger equation.
\newblock {\em arXiv preprint arXiv:2112.03464}, 2021.

\bibitem{cong}
Cong Hongzi.
\newblock The existence of full dimensional kam tori for nonlinear schrödinger
  equation.
\newblock {\em Mathematische Annalen}, pages 1--49, 2023.

\bibitem{imekraz2016long}
R.~Imekraz.
\newblock Long time existence for the semi-linear beam equation on irrational
  tori of dimension two.
\newblock {\em Nonlinearity}, 29(10):3067, 2016.

\bibitem{ionescu2019long}
A.D. Ionescu and F.~Pusateri.
\newblock Long-time existence for multi-dimensional periodic water waves.
\newblock {\em Geometric and Functional Analysis}, 29(3):811--870, 2019.

\bibitem{kappeler2021stability}
T.~Kappeler and R.~Montalto.
\newblock On the stability of periodic multi-solitons of the kdv equation.
\newblock {\em Communications in Mathematical Physics}, 385(3):1871--1956,
  2021.

\bibitem{kolmogorov1954conservation}
K.~N. Kolmogorov.
\newblock On conservation of conditionally periodic motions for a small change
  in hamilton's function.
\newblock In {\em Dokl. akad. nauk Sssr}, volume~98, pages 527--530, 1954.

\bibitem{K}
S.~B. Kuksin.
\newblock {\em Nearly integrable infinite-dimensional {H}amiltonian systems},
  volume 1556 of {\em Lecture Notes in Mathematics}.
\newblock Springer-Verlag, Berlin, 1993.

\bibitem{kuksin1997small}
S.~B. Kuksin.
\newblock On small-denominators equations with large variable coefficients.
\newblock {\em Zeitschrift f{\"u}r angewandte Mathematik und Physik ZAMP},
  48(2):262--271, 1997.

\bibitem{Kuk00}
S.~B. Kuksin.
\newblock {\em Analysis of {H}amiltonian {PDE}s}, volume~19 of {\em Oxford
  Lecture Series in Mathematics and its Applications}.
\newblock Oxford University Press, Oxford, 2000.

\bibitem{KP}
S.~B. Kuksin and J.~P\"{o}schel.
\newblock Invariant {C}antor manifolds of quasi-periodic oscillations for a
  nonlinear {S}chr\"{o}dinger equation.
\newblock {\em Ann. of Math. (2)}, 143(1):149--179, 1996.

\bibitem{liu2010spectrum}
J.~Liu and X.~Yuan.
\newblock Spectrum for quantum duffing oscillator and small-divisor equation
  with large-variable coefficient.
\newblock {\em Communications on pure and applied mathematics},
  63(9):1145--1172, 2010.

\bibitem{LY11}
J.~Liu and X.~Yuan.
\newblock A {KAM} theorem for {H}amiltonian partial differential equations with
  unbounded perturbations.
\newblock {\em Comm. Math. Phys.}, 307(3):629--673, 2011.

\bibitem{LiuYuan2014}
J.~Liu and X.~Yuan.
\newblock K{AM} for the derivative nonlinear {S}chr\"{o}dinger equation with
  periodic boundary conditions.
\newblock {\em J. Differential Equations}, 256(4):1627--1652, 2014.

\bibitem{maspero2018long}
A.~Maspero and M.~Procesi.
\newblock Long time stability of small finite gap solutions of the cubic
  nonlinear schr{\"o}dinger equation on t2.
\newblock {\em Journal of Differential Equations}, 265(7):3212--3309, 2018.

\bibitem{moser1962}
J.~M{\"o}ser.
\newblock On invariant curves of area-preserving mappings of an annulus.
\newblock {\em Nachr. Akad. Wiss. G{\"o}ttingen, II}, pages 1--20, 1962.

\bibitem{Neh}
N.~N. Nekhoro\v{s}ev.
\newblock An exponential estimate of the time of stability of nearly integrable
  {H}amiltonian systems.
\newblock {\em Uspehi Mat. Nauk}, 32(6(198)):5--66, 287, 1977.

\bibitem{Pos93}
J.~P\"{o}schel.
\newblock Nekhoroshev estimates for quasi-convex {H}amiltonian systems.
\newblock {\em Math. Z.}, 213(2):187--216, 1993.

\bibitem{P1}
J.~P\"{o}schel.
\newblock A {KAM}-theorem for some nonlinear partial differential equations.
\newblock {\em Ann. Scuola Norm. Sup. Pisa Cl. Sci. (4)}, 23(1):119--148, 1996.

\bibitem{P2}
J.~P\"{o}schel.
\newblock Quasi-periodic solutions for a nonlinear wave equation.
\newblock {\em Comment. Math. Helv.}, 71(2):269--296, 1996.

\bibitem{P3}
J.~P{\"o}schel.
\newblock On the construction of almost periodic solutions for a nonlinear
  schr{\"o}dinger equation.
\newblock {\em Ergodic Theory and Dynamical Systems}, 22(5):1537--1549, 2002.

\bibitem{W}
C.~E. Wayne.
\newblock Periodic and quasi-periodic solutions of nonlinear wave equations via
  {KAM} theory.
\newblock {\em Comm. Math. Phys.}, 127(3):479--528, 1990.

\bibitem{yuan}
X.~Yuan.
\newblock K{AM} theorem with normal frequencies of finite limit points for some
  shallow water equations.
\newblock {\em Comm. Pure Appl. Math.}, 74(6):1193--1281, 2021.

\bibitem{yuan2014long}
X.~Yuan and J.~Zhang.
\newblock Long time stability of hamiltonian partial differential equations.
\newblock {\em SIAM Journal on Mathematical Analysis}, 46(5):3176--3222, 2014.

\bibitem{yuan2016averaging}
X.~Yuan and J.~Zhang.
\newblock Averaging principle for the kdv equation with a small initial value.
\newblock {\em Nonlinearity}, 29(2):603, 2016.

\bibitem{yuan2013reduction}
X.~Yuan and K.~Zhang.
\newblock A reduction theorem for time dependent schr{\"o}dinger operator with
  finite differentiable unbounded perturbation.
\newblock {\em Journal of Mathematical Physics}, 54(5):052701, 2013.

\end{thebibliography}

\end{document}